\documentclass[12pt]{amsart}
\usepackage{mathpazo}
\usepackage{avant}
\usepackage{amsmath,amssymb,tensor}
\usepackage{mathrsfs}
\usepackage{tikz-cd}
\usepackage{booktabs,longtable,array}
\usepackage[backend=biber,style=alphabetic,maxbibnames=99]{biblatex}
\usepackage[hidelinks]{hyperref}
\hypersetup{
  pdftitle={Coefficient Riemann-Hilbert for an Affine Singular Foliation},
  pdfauthor={Qingyun Zeng}
}
\AtBeginBibliography{%
  \setlength{\emergencystretch}{3em}%
  \raggedright
  \interlinepenalty=10000
}

\newtheorem{thm}{Theorem}[section]
\newtheorem{cor}[thm]{Corollary}
\newtheorem{lem}[thm]{Lemma}
\newtheorem{prop}[thm]{Proposition}

\newtheorem{maintheorem}{Theorem}

\theoremstyle{definition}
\newtheorem{defn}[thm]{Definition}
\newtheorem{example}[thm]{Example}

\theoremstyle{remark}
\newtheorem{rem}[thm]{Remark}

\numberwithin{equation}{section}

\newcommand{\R}{\mathbb R}

\newcommand{\Z}{\mathbb Z}

\newcommand{\CB}{\mathcal B}

\newcommand{\CF}{\mathcal F}
\newcommand{\CG}{\mathcal G}

\newcommand{\E}{\mathcal E}
\newcommand{\li}{L_\infty}
\newcommand{\ai}{A_\infty}

\newcommand{\Wone}{W_1}

\newcommand{\Ggeom}{\mathcal G_M}

\newcommand{\NQ}{NQ}
\newcommand{\Coh}{\mathsf{Coh}}

\newcommand{\Perf}{\mathsf{Perf}}
\newcommand{\Rep}{\mathsf{Rep}}
\newcommand{\Cat}{\mathsf{Cat}}
\newcommand{\Aff}{\mathsf{Aff}}

\newcommand{\y}{\mathsf y}

\newcommand{\weak}{\xrightarrow{\sim}}

\DeclareMathOperator{\Map}{Map}
\DeclareMathOperator{\Hom}{Hom}
\DeclareMathOperator{\End}{End}

\DeclareMathOperator{\Int}{Int}

\DeclareMathOperator{\id}{Id}

\DeclareMathOperator{\pr}{pr}
\DeclareMathOperator{\im}{im}

\DeclareMathOperator{\Tot}{Tot}

\title[Affine coefficient Riemann--Hilbert]
{Coefficient Riemann--Hilbert for an Affine Singular Foliation}
\author{Qingyun Zeng}
\date{October 30, 2026; revision 6}
\keywords{singular foliation, Lie infinity-algebroid, higher integration,
  cohesive module, representation up to homotopy,
  Riemann--Hilbert correspondence}

\begin{document}

\begin{abstract}
For the affine foliation $\langle x\partial_x,y\partial_x\rangle$,
we prove that smooth coefficient integration is a strictly unital
$A_\infty$ quasi-equivalence on every ordinary open subset of
$\mathbb R^2$, for globally bounded, degreewise finite-projective
real coefficient bodies, allowing internal cohomology of
nonconstant rank. A horizontal-hull restriction of a finite Lie
$2$-groupoid model admits a smooth simplicial section of the
comparison from full Maurer--Cartan integration, with a homotopy
relative to all unit-generated simplices. Strict differentiation
and two relative Hom models yield a real-linear chain map whose
internal-degree diagonal is the identity. Its inverse is a finite
sum acting directly on the possibly rank-jumping body complex.
Global acyclicity of the differentiation kernels gives
invertible counits that are the identity on the body.
Infinitely flat and nongeneric rank-jumping examples illustrate
the nonregular theorem. Separately, regular two-term coefficients
admit a residue classification and an exact ordinary-quotient
descent criterion. Nonzero residue obstructs pullback through any
strict map to an ordinary groupoid nerve with represented objects
in the stated raw coefficient theory. For regular vertical
families over $S\times\mathbb R^2$, we also prove a
finite-projective bundle-valued theorem natural under smooth
parameter pullback; its two-term invariants are bundle arrows,
rather than pointwise rank profiles.
\end{abstract}

\maketitle
\clearpage
\tableofcontents

\section{Introduction}\label{sec:introduction}

\subsection{A singular foliation and its coefficients}

Let
\[
 M=\R^2,\qquad
 \CF=\langle x\partial_x,y\partial_x\rangle.
\]
The foliation has rank one away from the origin and rank zero at the
origin. Its module of generators is not locally free, but has the
global resolution
\[
 0\longrightarrow C^\infty(M)
 \xrightarrow{(-y,x)}C^\infty(M)^2
 \xrightarrow{(x\partial_x,y\partial_x)}\CF(M)\longrightarrow0.
\]
The associated geometric Lie $2$-algebroid $E_{\mathrm{aff}}$ has CE algebra
\[
 B=C^\infty(M)[\theta_1,\theta_2,\eta],
 \qquad |\theta_i|=1,\quad|\eta|=2,
\]
with differential
\begin{equation}\label{eq:intro-affine-ce}
 \begin{aligned}
 Qf&=f_x(x\theta_1+y\theta_2),&
 Q\theta_1&=-y\eta,\\
 Q\theta_2&=\theta_1\theta_2+x\eta,&
 Q\eta&=0.
 \end{aligned}
\end{equation}
The degree-two generator records a relation between vector-field
generators. At the origin it contributes genuine higher isotropy.
This paper proves coefficient Riemann--Hilbert on every ordinary
open $O\subseteq M$, including opens containing the origin.
The coefficient bodies are globally bounded and degreewise finite
projective over $C^\infty(O)$; their cohomology ranks may jump.
We then determine how higher isotropy is retained by integration
and which regular two-term coefficients can pass to ordinary quotients.

Geometric resolutions and their universal Lie infinity-algebroid
structures provide a general way to organize such relations
\cite{LGLS2020}. Their homotopy theory is not the same question as
constructing the ordinary holonomy groupoid
\cite{AndroulidakisSkandalis2009}. Here we retain full
Maurer--Cartan simplices and work on a fixed smooth base.
For a geometric algebroid $E$, full smooth integration is
\begin{equation}\label{eq:introduction-integration}
 \Int(E)_n(U)=
 \Map_{\NQ}\bigl(U\times T[1]\Delta^n,E[1]\bigr),
\end{equation}
with differential in the simplex direction and joint smoothness
through all faces. A simplex includes its base map and anchor
equations, not just pointwise Lie-algebra-valued forms.

\subsection{The finite model and the Riemann--Hilbert theorem}

Let $G=\Aff_+(1)$ act by $(a,b)(x,y)=(ax+by,y)$.
The crossed module defined in Section~\ref{sec:affine-integration}
has a strict Lie $2$-groupoid nerve
\[
 K_n=M\times G^n\times\R^{\binom n2}.
\]
This model appeared in the companion transport construction
\cite{ZengHigherHolonomy}. We recall its crossed-module structure
and horn properties. The comparison with full MC integration
and the coefficient results are the subject of the present paper.
For an open $O\subseteq M$, put $E_O=E_{\mathrm{aff}}|_O$ and
$B_O=B_{E_O}$. The appropriate finite-model restriction is
\begin{equation}\label{eq:intro-horizontal-hull}
 K[O]_n=
 \left\{k\in K_n:
 [\min_i x_i,\max_i x_i]\times\{y\}\subset O\right\},
\end{equation}
where $(x_i,y)$ are its vertices. The horizontal-hull condition
keeps the vertices in a single horizontal component of $O$;
vertex membership alone does not ensure this. Likewise
$\Int(E_O)$ uses full MC simplices whose entire base image
lies in $O$.

Here $\y$ denotes the represented smooth sheaf, and $c(\y O)$
is the constant simplicial subsheaf of unit-generated simplices.
Write $\Coh(B_O)$ for cohesive coefficients with real, globally
bounded, degreewise finite-projective bodies over $C^\infty(O)$.
Regularity means locally constant ranks of internal cohomology.
The representation notation refers to the concrete normalized
smooth cochain category of Section~\ref{sec:representations},
with the same coefficient-body conditions. Throughout,
representations retain their specified vertex complexes and
normalized cochain operations.

\begin{maintheorem}[Ordinary-open full affine comparison]
\label{main:affine-integration}
For every ordinary open $O\subseteq\R^2$, $K[O]$ is a
geometric Lie $2$-groupoid and there is an actual
vertex-preserving comparison
\[
 \Phi_O:\Int(E_O)\longrightarrow\y K[O].
\]
It is a trivial Kan fibration on every smooth parameter manifold.
It admits a smooth simplicial section $\Psi_O$, with
$\Phi_O\Psi_O=1$, and a homotopy from $\Psi_O\Phi_O$ to the identity,
over $\y K[O]$ and fixed on the section image. In particular
the homotopy is relative to the whole $c(\y O)$.
Actual pullback gives a dg quasi-equivalence
\[
 \Phi_O^*:\Rep_{\mathrm{sm}}^{\mathrm{perf}}(\y K[O])
 \longrightarrow
 \Rep_{\mathrm{sm}}^{\mathrm{perf}}(\Int(E_O)).
\]
This also restricts to the regular full subcategories.
\end{maintheorem}

The map records edge holonomies and
$h_{ijk}=-\int_{[ijk]}\beta$, where $\beta$ is the degree-two
MC field. Boundary extension, a Whitney--radial identity and a
rank-uniform flat correction prove the comparison. The smooth
simplicial retraction supplies more than an abstract weak
equivalence: it justifies a change of the actual raw coefficient
target. For a general $O$, smooth charts on the convex admissible
arrow domains keep the whole base simplex in $O$. This
range-preserving construction is distinct from the whole-plane
section. It neither divides by a rank coordinate nor
removes the origin. For different opens the sections are chosen
independently; strict restriction-naturality is an additional
compatibility condition.

\begin{maintheorem}[Ordinary-open coefficient RH]\label{main:affine-rh}
For every ordinary open $O\subseteq\R^2$, the actual functors
\[
 \begin{aligned}
 I_O=\Psi_O^*\mathsf{IN}_{E_O}&:
 \Coh(B_O)\longrightarrow
 \Rep_{\mathrm{sm}}^{\mathrm{perf}}(\y K[O]),\\
 \mathsf{IN}_{E_O}&:
 \Coh(B_O)\longrightarrow
 \Rep_{\mathrm{sm}}^{\mathrm{perf}}(\Int(E_O))
 \end{aligned}
\]
are smooth strictly unital $\ai$ quasi-equivalences, without a
regularity assumption on internal cohomology. Every integer Hom
degree is covered. On the finite model, actual differentiation
$F$ supplies a preimage $FR$ on the same body as $R$ and an actual
closed invertible morphism
\[
 q_R:I_O(FR)\longrightarrow R,\qquad (q_R)_0=1,\quad F(q_R)=1.
\]
The full-MC comparison likewise has actual invertible
identity-body counits. The preimages and counits are chosen
objectwise.
\end{maintheorem}

The proof works directly on the body complex, including where
its cohomology ranks jump.
Strict differentiation is constructed on the actual normalized
cochains and on globally finite-projective bundles. After a
closed rectification of edge transport, two prescribed
upper-decalage Hom models give a comparison $T$ whose internal
diagonal is the identity. If $[\ell,u]$ bounds the internal Hom
degrees, then
\[
 T=1+N,\qquad N^{u-\ell+1}=0,\qquad
 T^{-1}=\sum_{j=0}^{u-\ell}(-N)^j.
\]
A weighted normal-ideal estimate and an all-degree scalar
calibration prove the diagonal statement, including every power
of the even generator $\eta$. This is a real-linear chain inverse.
Internal degree is not a filtration by subcomplexes when the
body differential is nonzero. The inverse is constructed in
real vector spaces, and the two Hom differentials are related
by conjugation.

Actual rectification extends the Hom comparison to all pairs.
Global graded surjectivity then gives acyclicity of the actual
smooth differentiation kernels. Solving one equation in that
kernel constructs $q_R$; its inverse is finite because its
positive-arity correction lowers internal degree.
The whole-unit-relative prism passes the conclusion to full MC.
Infinitely flat and nongeneric rank-jumping coefficients make
the absence of a cohomology-regularity hypothesis concrete.

The whole-plane scalar formula
$H^\bullet(B)=C^\infty(\R_y)[\eta]/(y\eta)$ does not extend by
restriction to arbitrary nonrectangular opens. Regular
minimalization and finite MC lifting instead provide the tools
for the regular classification and vertical-family results;
they are distinct from the direct body-complex argument of
Theorem~\ref{main:affine-rh}.

\subsection{A concrete higher-isotropy application}

The two-term coefficient sector makes the significance of the
higher class explicit. Here ``two-term'' means internal cohomology
in degrees $0,1$, allowing acyclic terms in the original complex.

\begin{maintheorem}[Two-term coefficients and ordinary quotients]
\label{main:two-term-applications}
The $H^0$-isomorphism classes of regular two-term cohesive
coefficients, and of their corresponding raw representations
on $\y K$ and full MC, are classified by
\[
 (d_0,d_1,r),\qquad 0\le r\le\min(d_0,d_1),
\]
where $d_i$ are the internal cohomology dimensions and $r$ is
the rank of a residue map $N_0:V^1\to V^0$.
A representative has $D(v_1)=N_0(v_1)\eta$ and integrates on $K$
to $r_2=-hN_0$.

There is an ordinary groupoid $\overline{\CG}$ in smooth sheaves,
with objects $\y M$, and a strict map
$q:\y K\to N(\overline{\CG})$ fixing objects.
Such a representation descends along $q$ exactly when $r=0$.
For $r>0$ it cannot be a pullback, up to $H^0$-isomorphism,
through any strict map to an ordinary groupoid nerve with
represented objects, using the bounded finite-projective raw
coefficient theory. The same obstruction holds on full MC.
\end{maintheorem}

Section~\ref{sec:applications} proves the normal forms by a global
smooth gauge and an exact residue calculation. The residue matrix
is considered up to graded frame changes, not as a canonical
numeric matrix. The obstruction restricts to the exact
$B^2\Z$ at the origin. Ordinary groupoid representations can have
higher components, but normalization makes their pullbacks
constant on this particular simplicial subobject.
The quotient is constructed in smooth sheaves, without a
representability hypothesis on its arrow sheaf. Comparison with
Androulidakis--Skandalis (AS) holonomy or refined higher holonomy
requires additional geometric input. The obstruction allows
nonregular and nonminimal coefficients on the ordinary target.

\subsection{Vertical families and nontrivial coefficient bundles}

We also study regular coefficients for the vertical algebroid
$E_S$ on $S\times\R^2$, for every finite-dimensional Hausdorff
second-countable smooth manifold $S$, including noncompact
and disconnected bases. Its anchor has no $S$-component:
this is not an inverse-image algebroid adding $TS$.
Theorem~\ref{thm:family-rh} proves actual regular RH for
globally bounded finite-projective coefficient bundles,
without requiring those bundles to be trivial.
Here regularity is a standing hypothesis on the fiber cohomology
of the body complex over $S\times\R^2$. Smooth parameter dependence
of the simplices in Theorem~\ref{main:affine-integration} and
regularity of a coefficient bundle are distinct conditions.
Corollary~\ref{cor:family-parameter-naturality} gives compatibility
with all smooth parameter pullbacks and a formal equivalence
of the corresponding $\Cat_\infty$-valued parameter diagrams.

This requires more than tensoring the scalar comparison with
an arbitrary function space. A normalized global projector
dictionary first makes raw finite perturbation valid for
nontrivial vertex bundles. The reduction then produces a bundle
over $S$, and finite-projective layers over $C^\infty(S)$ permit
global Hom comparison and MC lifting. The actual counit is
given in \eqref{eq:family-counit} and is chosen objectwise.

The two-term family invariant is a bundle arrow
$N:H^1\to H^0$ over $S$, up to bundle isomorphism
(Theorem~\ref{thm:family-two-term}). The M\"obius-line example
detects coefficient topology, while $t$ and $t^3$ detect
smooth vanishing order despite identical rank profiles.
The rank of $N$ may vary.
Corollary~\ref{cor:family-ordinary-descent} proves the exact
ordinary-descent criterion $N\equiv0$. The family theorem varies
the parameter space and coefficient bundles while keeping the
affine algebroid fixed.

\subsection{Relation to earlier work and organization}

The general fixed-base theory is developed in the companion
\emph{Global $\li$-Algebroid Resolutions: Homotopy Invariance and
Smooth Coefficient Integration} \cite{ZengGlobalResolutions2026}.
We use its pointed polynomial comparison, coherent
global-resolution results, coefficient invariance and smooth
coefficient-integration construction. Section~\ref{sec:preliminaries}
records the hypotheses needed here. The global-resolution
machinery compares other admissible global resolutions with the
selected affine model. The explicit RH proof for that model is
given here independently of those comparisons.

Nuiten supplies the fixed-base semi-model and polynomial
mapping-space inputs \cite{Nuiten2019}. For general weak maps,
we use the connective mapping-space comparison. Its connective
form is essential at the singular foliation, where the anchor
is not transitive.
Chuang--Holstein--Lazarev establish the general smooth and
normalized-cochain homotopy-invariance principles
\cite[Theorem~4.2, Corollary~3.1 and Remark~4.1;
Theorem~5.2, Proposition~6.2 and Corollary~6.1]{CHL2021}.
The companion gives the precise jointly smooth, finite-projective
and vertex-relative forms used here.

The universal integration kernel is due to Arias Abad--Sch\"atz
\cite{AriasAbadSchaetz2013}. The smooth strictly unital functor
used here is constructed in \cite[Theorem~9.2]{ZengGlobalResolutions2026}.
Section~\ref{sec:coefficient-integration} retains its units,
vertex behavior, matrix signs, actual Taylor formulas and
fixed-base naturality. The affine task is to prove that
\emph{this functor} is a quasi-equivalence and to compute its
residue-sensitive representations.

An actual higher Riemann--Hilbert correspondence for full
de Rham geometry on compact closed manifolds was established
by Block--Smith \cite[Theorem~4.1, p.~16]{BlockSmith2014}.
Ordinary-groupoid loss of higher coefficient information also
has a direct antecedent: Arias Abad--Sch\"atz give a two-term
$TS^2$ coefficient whose integration does not descend to the
pair groupoid
\cite[Section~5.6 and Proposition~5.4, pp.~39--41]{AriasAbadSchaetz2013}.
These references use the preprint versions specified
in the bibliography. Our contribution here is the actual-map
reconstruction for this singular affine model, its explicit
residue and bundle-arrow classifications, and the obstruction
for strict maps to ordinary groupoid nerves with represented objects.

A complementary regular-foliation correspondence is developed in
\cite{ZengRegularRH}. On a compact base it compares globally bounded
finite-rank regular cohesive modules with global plot-smooth objects
of the same regularity, using sheafified raw Hom presheaves.
The affine equivalence here
instead uses actual raw Hom complexes and permits cohomology ranks
to jump. These different targets must be kept distinct when
comparing the respective descent statements.

Higher Van Est theory provides an important antecedent
\cite{AnguloCueca2025,CabreraDelHoyo2026}. The last-horn criterion
of \cite[Lemma~4.3.3 and Theorem~4.3.4, pp.~28--29]{CabreraDelHoyo2026}
applies to this $K$ in every scalar degree.
That theorem concerns a differentiation map into its intrinsic
CE quotient. Here we identify the affine quotient, construct
actual coefficient differentiation, and calibrate the prescribed
relative comparison in every internal degree. The reverse
functor is the specified $\Psi_O^*\mathsf{IN}$, rather than an
arbitrarily chosen inverse to a scalar Van Est map.

After the conventions and cited interfaces,
Section~\ref{sec:affine-integration} introduces the affine crossed
module and the full-MC comparison map; its longer global lifting
and retraction proofs are in Appendix~\ref{app:affine-global-lifting}.
Sections~\ref{sec:representations}--\ref{sec:coefficient-integration}
fix the raw coefficient target, prove its actual affine pullback
comparison, and record the integration formulas.
Section~\ref{sec:affine-open-geometry} constructs the
range-preserving ordinary-open comparison. The central argument in
Sections~\ref{sec:affine-differentiation}--\ref{sec:affine-nonregular-rh}
is strict differentiation, the prescribed relative Hom comparison,
and Theorem~\ref{thm:affine-open-rh}.
Section~\ref{sec:nonregular-examples} illustrates its new scope
with rank-jumping coefficients and a nonrectangular open.
Section~\ref{sec:applications} gives the point-parameter residue
classification and ordinary non-descent theorem.
Section~\ref{sec:affine-families} proves the vertical-family
extension and its bundle-valued examples.
Section~\ref{sec:programme} formulates questions for further work.
Appendices~\ref{sec:affine-hom} and~\ref{sec:affine-reconstruction}
retain the regular reduction, scalar-parameter and finite
MC-lifting tools used by the classification and family arguments,
separately from the nonregular proof.

\section{Conventions and comparison results}
\label{sec:preliminaries}

Manifolds are finite-dimensional, Hausdorff and second-countable.
Write $A=C^\infty(M)$. Algebroid complexes are
homological, with $\partial:E_i\to E_{i-1}$, and we use Nuiten's
convention $l_1=-\partial$, $|l_n|=n-2$. The geometric CE generator
dual to $E_i$ has cohomological degree $i+1$. Coefficient complexes
are cohomological.

\subsection{Geometric objects and supplied global resolutions}
\label{sec:polynomial}

Let $\Ggeom$ be the category of split geometric $\li$-algebroids
over $M$ whose homological section complexes are connective and
degreewise finitely generated projective over $A$. Such complexes
are underlying-$A$-cofibrant. Morphisms are fixed-base nonlinear
infinity-morphisms, and $\Wone$ consists of those whose first Taylor
component is a quasi-isomorphism. Their Taylor components factor
$A$-linearly through the corresponding symmetric powers and satisfy
the anchor condition separately. The CE differential and anchored
brackets obey their respective Leibniz rules over $A$. The geometric
dictionary and underlying-cofibrancy statement are given in
\cite[Section~3, Lemmas~3.1 and 3.6]{ZengGlobalResolutions2026}.

For $E\in\Ggeom$, write $B_E$ for the graded algebra of global
geometric CE sections with its differential $Q$. Positive grading
forces $Qf$ for $f\in A$ to be linear in the degree-one generators.
Thus the ordinary-body anchor is fiber-linear, also at singular rank.
Infinity-morphisms may nevertheless be nonlinear in positive
coordinates. Throughout, the base is an ordinary smooth manifold.

Write $T_A=\operatorname{Der}_{\R}(A)$.
A geometric resolution of a foliation module
$\mathfrak f\hookrightarrow T_A$ is a strict augmentation
$p_E:E\weak\mathfrak f$ from such an $E$ whose underlying global
section complex resolves $\mathfrak f$ in degree zero.
Injectivity here is injectivity of the module anchor; it need not
survive passage to a residue-field fiber. All resolution-comparison
statements concern global resolutions with these hypotheses.
The module $\mathfrak f$ may be nonprojective. Global finite
projectivity of the geometric terms is an explicit hypothesis;
local finite rank or exactness of residue-field fibers is not a
substitute. Resolution comparisons take place among these
supplied global objects.

\subsection{Smooth sheaves and full Maurer--Cartan simplices}

Use the ordinary site of Euclidean opens with smooth maps and
open covers, including the empty cover of the empty object.
Write $\E$ for its sheaf topos and $\y M(U)=C^\infty(U,M)$.
Euclidean opens are dense, so these sheaves can also be evaluated
on smooth manifolds. Simplicial smooth sheaves are objects of
$s\E$. Their tensor with a simplicial set is the sheafified
copower; simplex labels may be locally constant on disconnected
tests. A finite simplicial set has finitely many nondegenerate
simplices, and $\Lambda^0[1]=\{0\}$.

An $\NQ$-map is a map of the graded manifolds that commutes
with their cohomological vector fields. Full smooth integration is
\[
 \Int(E)_n(U)=
 \Map_{\NQ}\bigl(U\times T[1]\Delta^n,E[1]\bigr).
\]
The differential acts only in the simplex direction. For an
infinite positive grading this means degreewise geometric CE
maps. A simplex contains its entire base map, anchor equation and
higher flatness equations. All data are jointly smooth through
the faces, not an algebraic tensor with $C^\infty(U)$.
Gluing in smooth parameters for a fixed algebroid gives a
simplicial sheaf \cite[Lemma~4.1]{ZengGlobalResolutions2026}.

For $S_0=\y M$, write $c(\y M)\subseteq S$ for the whole
unit-generated subsheaf of totally degenerate simplices.
A homotopy relative to $c(\y M)$ is constant on
$c(\y M)\times\Delta[1]$, not merely at vertex endpoints.
For integration these are constant base maps with zero
positive-degree fields.

Write $|S|$ for open-cover infinity-sheafification of the
objectwise realization, and
$\operatorname{Mon}^{\mathrm{sh}}(E)=|\Int(E)|$.
We also consider the hypercompletion of this smooth infinity-sheaf.
The associated derived fixed-object Cech groupoid has levels
$\y M\times^h_{|S|}\cdots\times^h_{|S|}\y M$.
The raw coefficient category of Section~\ref{sec:representations}
is defined on $S$ itself, with its specified vertices and simplices.

\subsection{Cohesive coefficients}
\label{sec:cohesive}

Let $\CB=\bigoplus_{r\ge0}\CB^r$ be a sheaf of commutative
cochain dg real algebras with $\CB^0=C^\infty_M$.
Each $\CB^r$ is a finite-rank bundle sheaf, multiplication is
smooth, and $d_B$ is a local differential operator commuting
with differentiation and compact integration in external smooth
parameters. The geometric CE sheaves above satisfy these
hypotheses. Put $B=\bigoplus_r\Gamma(M,\CB^r)$.

Following \cite{Block2010}, an object of $\Coh(B)$ is a globally
bounded graded finite-projective $A$-module $V$ with a flat
degree-one right superconnection on $V\otimes_A B$:
\[
 D(vb)=D(v)b+(-1)^{|v|}v\,d_Bb,\qquad D^2=0.
\]
A homogeneous degree-$m$ morphism is signlessly right-linear,
$\phi(vb)=\phi(v)b$, and has differential and composition
\[
 d\phi=D_W\phi-(-1)^m\phi D_V,\qquad
 (\psi\phi)(v)=\psi(\phi(v)).
\]
The internal differential is $D^{[0]}$.
The regular full subcategory $\Coh^{\mathrm{reg}}(B)$ requires
locally constant fiber-cohomology ranks. It is shift-closed;
a cone belongs to it when its fiber-cohomology ranks are locally
constant. All coefficient ranks and
amplitudes are globally bounded. In a fixed Hom degree $m$,
only finitely many terms with $j+r=i+m$ occur. The degrees
of $B$ may be unbounded.
These are the conventions of
\cite[Section~6]{ZengGlobalResolutions2026}. In interval
formulas $dt$ is on the left.

The notation $\Rep_{\mathrm{sm}}^{\mathrm{perf}}$ refers to the actual
bounded finite-rank vertex-complex category of
Section~\ref{sec:representations}. Its definition retains the
specified vertex complexes and normalized cochains, rather than
passing to a derived-perfect replacement or localizing over stack
presentations.

\subsection{Fixed-base and global-resolution comparisons}
\label{sec:coherent}

The following results compare other admissible global resolutions
with the selected affine resolution. The explicit affine comparison
and its RH proof proceed independently of these results.
All categories are taken in compatible universes.

Let $\varepsilon_E:B_E\to A$ be the graded zero-section
augmentation. A polynomial parameter map
$H:B_{E'}\to B_E\otimes_\R\Omega^\bullet_{\mathrm{poly}}(\Delta^n)$
is pointed when
\[
 (\varepsilon_E\otimes\id)H=\iota\varepsilon_{E'},
\]
where $\iota$ is constant extension.
Fixing base functions alone does not impose this condition:
positive generators must have no pure parameter-form terms.
This is \cite[Definition~3.5 and Lemma~3.6]{ZengGlobalResolutions2026}.

For general weak maps, we use the connective mapping-space
comparison in \cite[Remarks~5.27--5.28]{Nuiten2019}, as specified in
\cite[Section~3]{ZengGlobalResolutions2026}. It identifies the
polynomial model with the localization of connective strict-anchor
algebroids and strict maps. Its target-fibrancy condition is
positive-degree surjectivity, automatic over the discrete base.
In contrast, the unbounded fibrancy condition requires a surjective
anchor, which a singular resolution need not have.

The next statement collects
\cite[Theorems~3.9, 4.5 and 6.4; Corollary~8.4]{ZengGlobalResolutions2026}.
\begin{thm}[Fixed-base comparison]
\label{main:geometric-invariance}
Under the preceding geometric hypotheses and connective
mapping-space input, a given $F:E\rightsquigarrow E'$ in $\Wone$
has a geometric inverse up to finite chains of pointed polynomial
homotopies. Its cohesive pullback is a dg quasi-equivalence,
including on regular coefficients. The actual map $\Int(F)$
is an objectwise simplicial homotopy equivalence, and its raw
representation pullback is a dg quasi-equivalence, also on the
regular full subcategories.
\end{thm}
The inverse homotopies can be made jointly smooth. They preserve
the entire base simplex and are relative to $c(\y M)$.
This is stronger than an abstract equivalence of smooth stacks.

The following resolution-specific route does not require that
general connective localization input.
For a degree-zero algebroid with injective module anchor
$\mathfrak f\hookrightarrow T_A$, let $\mathscr R_{\mathfrak f}$
be the nonempty simplicial category of its supplied strict
global geometric resolutions, with all actual reduced pointed
polynomial mapping simplices.

The resolution mapping spaces satisfy the following uniqueness
result \cite[Theorem~5.2]{ZengGlobalResolutions2026}.
\begin{thm}[Coherent resolution uniqueness]
\label{thm:coherent-resolutions}
For $E,E'\in\mathscr R_{\mathfrak f}$, the mapping space is a
weakly contractible Kan complex. The homotopy-coherent nerve
$N_{\mathrm{hc}}\mathscr R_{\mathfrak f}$ is contractible.
Every mapping simplex automatically commutes with the resolution
augmentations.
\end{thm}

Integration extends to these mapping spaces by
\cite[Theorem~5.3 and equation~(5.2)]{ZengGlobalResolutions2026}.
\begin{thm}[Enriched full-MC integration]
\label{thm:enriched-integration}
Full smooth integration sends pointed polynomial mapping
simplices to smooth simplicial mapping simplices, compatibly
with identities and composition. For a mapping $n$-simplex $H$,
an MC $m$-simplex $\sigma$ and $\tau:[m]\to[n]$, its value is
\[
 H_{\NQ}\circ(\sigma,T[1]\lambda_\tau),\qquad
 (\lambda_\tau(t))_j=\sum_{\tau(i)=j}t_i.
\]
It preserves the actual base simplex and fixes the whole
$c(\y M)$ throughout every mapping simplex.
\end{thm}

Consequently, one obtains the following framed comparison
\cite[Corollary~5.5]{ZengGlobalResolutions2026}.
\begin{cor}[Framed monodromy comparison]
\label{cor:coherent-monodromy}
For either smooth infinity-target, the integration diagram
on $N_{\mathrm{hc}}\mathscr R_{\mathfrak f}$ is equivalent to
a constant diagram. After selecting a resolution and framing
the component there by the identity, the space of coherent
comparisons is contractible.
\end{cor}

For actual comparisons between these global resolutions, the
cohesive and raw representation pullbacks are dg
quasi-equivalences, including on regular coefficients, by
\cite[Corollaries~6.3 and 8.5]{ZengGlobalResolutions2026}.
Thus the supplied global resolutions give equivalent cohesive
and raw representation categories.

\paragraph{Citation versions.}
Numbered citations to \cite{Nuiten2019} follow arXiv:1712.03441v2;
those to \cite{CHL2021} use the published article. Page references
to \cite{NSS2015} use arXiv:1207.0249v1.
For \cite{AriasAbadSchaetz2013} the numbered source is
arXiv:1011.4693v2. The higher Van Est references
\cite{CabreraDelHoyo2026,AnguloCueca2025} use their respective
version~2 preprints. Citations to \cite{ZengGlobalResolutions2026}
use its October 30, 2026 revision 8.

\section{A finite Lie 2-groupoid model for the affine resolution}
\label{sec:affine-integration}

We now compare full smooth MC integration with a finite-dimensional
model, without truncating the source simplicial sheaf. The example
is the resolution of
\[
 M=\R^2,\qquad
 \CF=\langle x\partial_x,y\partial_x\rangle
\]
used in the transport construction of
\cite{ZengHigherHolonomy}.
The point-base Morita comparison of
\cite[arXiv v4, Theorem~2.15, p.~11]{ShengZhu2012}
concerns the 2-truncated integration of a finite-dimensional
strict Lie 2-algebra.
Here the comparison allows nonconstant base maps in the MC
simplices, arbitrary smooth parameter families, and all simplicial degrees.

\subsection{The resolution and the strict target}

Put $A=C^\infty(\R^2)$. The resolution has
\[
 0\longrightarrow A\varepsilon\xrightarrow{d_2}
 Ae_1\oplus Ae_2\xrightarrow{\rho}\CF(M)\longrightarrow0,
 \qquad
 d_2\varepsilon=-ye_1+xe_2,
\]
with $\rho(e_1)=x\partial_x$ and $\rho(e_2)=y\partial_x$.
Its CE differential is
\begin{equation}\label{eq:affine-mc-ce}
 \begin{aligned}
 Qf&=(xf_x)\theta_1+(yf_x)\theta_2,
 &Q\theta_1&=-y\eta,\\
 Q\theta_2&=\theta_1\theta_2+x\eta,
 &Q\eta&=0,
 \end{aligned}
 \qquad |\theta_i|=1,\quad|\eta|=2.
\end{equation}
We denote this geometric algebroid by $E_{\mathrm{aff}}$.

Let
\[
 \begin{gathered}
 G=\Aff_+(1)
 =\left\{\begin{pmatrix}a&b\\0&1\end{pmatrix}:a>0\right\},
 \\
 (a,b)(x,y)=(ax+by,y),
 \qquad
 Z_m=\begin{pmatrix}-y_m&x_m\\0&0\end{pmatrix}.
 \end{gathered}
\]
The additive group bundle $M\times\R$ maps to $G\ltimes M$ by
\[
 \partial(m,h)=(\exp(hZ_m),m).
\]
The $G$-action transports $m$ and leaves the scalar $h$ unchanged.

\begin{prop}[Affine resolution and crossed module]\label{prop:affine-source-data}
The displayed complex is exact, and
\eqref{eq:affine-mc-ce} squares to zero.
The displayed group-bundle map is a crossed module.
Its strict Lie 2-groupoid has a Duskin nerve
\[
 K_n\cong M\times G^n\times\R^{\binom{n}{2}}.
\]
The horn maps are surjective submersions, with unique fillers
above degree two. Its smooth second isotropy is $\R$ at the
origin and zero elsewhere.
\end{prop}

\begin{proof}
If $xf+yg=0$, restriction to $x=0$ and continuity give
$g(0,y)=0$. Hadamard division gives $g=xh$; density of
$x\ne0$ gives $f=-yh$ and also proves injectivity of $d_2$.
The same integral division works globally on $\R^2$.
The graded Leibniz rule verifies $Q^2=0$ on the smooth
coordinate functions and the three positive generators;
for example
\[
 Q^2\theta_2
 =-y\eta\theta_2-x\theta_1\eta
       +(x\theta_1+y\theta_2)\eta=0.
\]

Matrix multiplication gives $Z_mm=0$ and
$gZ_mg^{-1}=Z_{gm}$. Moreover
\begin{equation}\label{eq:affine-entire-boundary}
 \exp(hZ_m)=
 \begin{pmatrix}
 e^{-hy_m}&x_mh\varphi(hy_m)\\0&1
 \end{pmatrix},
 \qquad
 \varphi(t)=\int_0^1e^{-st}\,ds.
\end{equation}
Thus the boundary is smooth, fixes $m$, and is equivariant.
The Peiffer identity follows from the trivial action on
the scalar and the commutativity of the coefficient group.

Fix the triangle convention
\begin{equation}\label{eq:affine-triangle-convention}
 g_{ik}=g_{jk}g_{ij}\exp(h_{ijk}Z_{m_i}).
\end{equation}
The root edges $a_i=g_{0i}$ and rooted triangle parameters
$h_{0ij}$ determine every other coordinate:
\begin{equation}\label{eq:affine-root-coordinates}
 m_i=a_i m_0,\qquad
 g_{ij}=a_j\exp(-h_{0ij}Z_{m_0})a_i^{-1},\qquad
 h_{ijk}=h_{0ij}+h_{0jk}-h_{0ik}.
\end{equation}
The tetrahedral identity is $\delta h=0$.
Degree-one horns have the source and target submersions
of the action groupoid. In degree two, each missing edge
is solved by multiplication and inverse for every chosen
$h$; these coordinate changes identify the horn map
with a projection with fiber $\R$.
The tetrahedral identity determines a missing degree-three
filler, and its cocycle identity gives unique higher
fillers. All formulas are smooth, also after applying
$C^\infty(U,-)$.
Finally, $\exp(hZ_m)=1$ forces $h=0$ when $y_m\ne0$,
or when $y_m=0,x_m\ne0$. At $m=0$ it imposes no
condition on $h$.
\end{proof}

The nerve has all degrees; an $n$-simplex is determined
by its compatible 2-skeleton, but not every 2-skeleton
extends. The scalar isotropy just described is a smooth
group of 2-automorphisms, not the ordinary homotopy group
of the topological manifold $\R$.

\subsection{The full MC equations and comparison map}

For a smooth parameter manifold $U$, a full MC simplex
has a base map $\sigma=(x,y):U\times\Delta^n\to M$ and
relative forms $\alpha_1,\alpha_2,\beta$ of degrees $1,1,2$.
All are jointly smooth through the faces, and
\begin{equation}\label{eq:affine-full-mc}
 \begin{gathered}
 dy=0,\qquad dx=x\alpha_1+y\alpha_2,\\
 d\alpha_1=-y\beta,\qquad
 d\alpha_2=\alpha_1\wedge\alpha_2+x\beta,\qquad
 d\beta=0.
 \end{gathered}
\end{equation}
Here and below $d$ differentiates only in $\Delta^n$.
Write
\[
 \mathsf A=\begin{pmatrix}\alpha_1&\alpha_2\\0&0\end{pmatrix},
 \qquad
 Z_\sigma=\begin{pmatrix}-y&x\\0&0\end{pmatrix}.
\]
The equations are equivalently
\begin{equation}\label{eq:affine-fake-curvature}
 d\sigma=\mathsf A\sigma,\qquad
 d\mathsf A-\mathsf A\wedge\mathsf A=Z_\sigma\beta,\qquad
 d\beta=0,
\end{equation}
and imply
$dZ_\sigma=\mathsf AZ_\sigma-Z_\sigma\mathsf A$.

For an oriented edge $[ij]$, let $g_{ij}$ be the
holonomy defined by the left matrix equation
\[
 P'(t)=\mathsf A(\partial_t)P(t),\qquad P(0)=1.
\]
For every oriented triangle set
\begin{equation}\label{eq:affine-surface-coordinate}
 h_{ijk}=-\int_{[ijk]}\beta.
\end{equation}

\begin{thm}[Full affine integration comparison]\label{thm:affine-integration}
These edge holonomies and triangle coordinates define a
natural vertex-preserving map of simplicial smooth sheaves
\[
 \Phi:\Int(E_{\mathrm{aff}})\longrightarrow\y K.
\]
For every smooth parameter manifold $U$, $\Phi(U)$ is a
trivial Kan fibration of simplicial sets. In particular,
it is an objectwise weak equivalence of full simplicial
sheaves, including families passing through the origin
and the locus $y=0$.
\end{thm}

Appendix~\ref{app:affine-global-lifting} proves the theorem
by constructing lifts against every boundary inclusion and
then supplies an actual smooth simplicial section and
whole-unit-relative homotopy. The comparison itself is
intrinsic; the auxiliary choices in a lift need not commute
with all simplex operators. Section~\ref{sec:affine-open-geometry}
gives the range-preserving ordinary-open construction.

\section{A concrete coefficient target on a simplicial smooth sheaf}
\label{sec:representations}

\subsection{Smooth cochains with two endpoints}

We use a raw dg category of representations on a simplicial
smooth sheaf. The construction is a standard
twisting-cochain form of homotopy-coherent
representations, with explicit totalization signs;
compare \cite[Definition~2.2 and Remark~2.4]{AriasAbadSchaetz2013}.
The plot-smooth family condition and the uniform
boundedness imposed here are additional requirements.
We give the definitions and signs below; the dg identities follow
from
\cite[Section~7]{ZengGlobalResolutions2026}.

Let $S$ be a simplicial smooth sheaf with $S_0=\y M$.
Let $V,W$ be globally bounded finite-rank graded
bundles on $M$, with finite-projective global
section modules. Write
\[
 C_{\mathrm{sm}}^{p,q}(S;V,W)
\]
for plot-smooth sections of
$\Hom^q(v_p^*V,v_0^*W)$ on $S_p$, for $p\ge0$.
Precisely, such a section assigns to each
$\sigma\in S_p(T)$ a smooth degree-$q$ bundle
map from $(v_p\sigma)^*V$ to $(v_0\sigma)^*W$
on the test manifold $T$, naturally under
smooth pullback. Transport acts from the last
vertex to the first.

The total degree is $p+q$. Set
\[
 C^m_{\mathrm{sm}}(S;V,W)
 =\prod_{p\ge0}C^{p,m-p}_{\mathrm{sm}}(S;V,W).
\]
This product is finite in each total degree:
the internal degrees of $\Hom(V,W)$ range
over a fixed finite interval.

For $\alpha$ of bidegree $(p,q)$ and
$\beta$ of bidegree $(r,s)$ with composable
coefficient bundles, define
\begin{equation}\label{eq:cochain-product}
 (\alpha\star\beta)(\sigma)
 =(-1)^{qr}
 \alpha(\sigma[0,\ldots,p])\,
 \beta(\sigma[p,\ldots,p+r]).
\end{equation}
Here juxtaposition is composition, with
the right factor applied first. Define the
inner-face differential by
\begin{equation}\label{eq:inner-differential}
 (\delta\alpha)(\sigma)
 =\sum_{i=1}^{p}(-1)^i\alpha(d_i\sigma),
 \qquad \sigma\in S_{p+1}.
\end{equation}
Inner faces retain the two endpoints, so
every term has the correct bundle type.

The product and inner differential satisfy the following identities
\cite[Lemma~7.1]{ZengGlobalResolutions2026}.
\begin{lem}\label{lem:cochain-identities}
The product is associative, $\delta^2=0$,
and
\[
 \delta(\alpha\star\beta)
 =(\delta\alpha)\star\beta
  +(-1)^{p+q}\alpha\star(\delta\beta).
\]
\end{lem}

\subsection{Representations and their Hom differential}

\begin{defn}\label{def:smooth-representations}
An object of $\Rep_{\mathrm{sm}}^{\mathrm{perf}}(S)$
is $(V,r_V)$, where
\[
 r_V=\sum_{p\ge0}r_{V,p},\qquad
 r_{V,p}\in C^{p,1-p}_{\mathrm{sm}}(S;V,V),
\]
satisfies
\[
 \delta r_V+r_V\star r_V=0.
\]
We impose strict unital normalization
\[
 r_{V,1}(s_0x)=\id_{V_x},\qquad
 r_{V,p}(s_i\sigma)=0\quad(p\ge2),
\]
and require every edge transport to be
a fiberwise quasi-isomorphism. A
degree-$m$ morphism $\Phi$ has components
$\Phi_p\in C^{p,m-p}_{\mathrm{sm}}(S;V,W)$,
vanishing on degeneracies for $p\ge1$.
Its differential is
\begin{equation}\label{eq:representation-differential}
 D\Phi=\delta\Phi+r_W\star\Phi
       -(-1)^m\Phi\star r_V.
\end{equation}
Composition is $\star$.
\end{defn}

Here $r_{V,0}$ is the internal
differential. The degree-one
Maurer--Cartan equation says that
$r_{V,1}$ is a chain map. In degree two
it reads
\[
 -r_{V,1}(d_1\sigma)
 +r_{V,0}r_{V,2}
 +r_{V,1}(\sigma[0,1])r_{V,1}(\sigma[1,2])
 +r_{V,2}r_{V,0}=0.
\]
Thus $r_{V,2}$ supplies the indicated
chain homotopy. More precisely, set
$\varepsilon_p=(-1)^{p(p-1)/2}$ and
$F_p=\varepsilon_p r_{V,p}$.
For $k\ge1$ the identity
\[
 (-1)^{(1-p)(k-p)}
 \varepsilon_{k-1}\varepsilon_p\varepsilon_{k-p}
 =(-1)^{p+1}
\]
converts the degree-$k$ structure equation to
\cite[Definition~2.2, equation~(1)]{AriasAbadSchaetz2013},
whose front/back cup has no internal Koszul sign.
For $k=0$ both equations say $r_{V,0}^2=0$,
and $\varepsilon_1=1$ preserves the strict unit.
This object rescaling is not the morphism dictionary.
Its all-degree extension, including the differential
and composition, is
Lemma~\ref{lem:in-sign-dictionary}.

If $S$ is Kan on point-valued simplices, inverse edges
and two-simplices give chain-homotopy inverses of the
edge maps, hence the required fiberwise
quasi-isomorphisms. They do not imply that the
original edge maps are degreewise invertible.
The edge quasi-isomorphism condition is included
explicitly so that the definition does not assume
raw Kan fibrancy for every integration used below.

These operations give the following dg structure
\cite[Proposition~7.3]{ZengGlobalResolutions2026}.
\begin{prop}\label{prop:representation-dg-category}
Definition~\ref{def:smooth-representations}
defines a dg category. The identity has
only its degree-zero component $\id_V$.
\end{prop}

The inner differential $\delta$ alone need not preserve
normalization. It is the full twisted differential
\eqref{eq:representation-differential}, with its unit terms,
that does so.

At simplicial degree zero,
$(D\Phi)_0=r_{W,0}\Phi_0-(-1)^m\Phi_0r_{V,0}$.
When $S$ is the nerve of an ordinary groupoid and the coefficient
bundles are concentrated in degree zero, these objects are
ordinary contravariant groupoid representations.

Let $\Rep_{\mathrm{sm}}^{\mathrm{perf,reg}}(S)$
be the full subcategory for which
$r_{V,0}$ has locally constant
fiber-cohomology ranks. A cone belongs
to this subcategory when its fiber-cohomology
ranks are locally constant.

Strict simplicial maps induce actual raw dg pullback functors
by precomposition and pullback of the vertex bundles.
The inverse-homotopy hypotheses of
Section~\ref{sec:representation-invariance} ensure that the pullback
is a quasi-equivalence.
For an open $O\subseteq M$, $\Int(E|_O)$ uses simplices whose
entire base image lies in $O$, not merely their vertices.
The affine RH theorem uses this raw category before
Hom-sheafification. The latter construction, described in
\cite[Proposition~7.4]{ZengGlobalResolutions2026}, sheafifies the
morphism complexes while retaining the existing global objects.

\section{Representation invariance and the affine coefficient comparison}
\label{sec:representation-invariance}

We use the actual raw dg category of
Section~\ref{sec:representations} on the simplicial
smooth sheaf itself.
Homotopy invariance for normalized-cochain
Maurer--Cartan and twisted-module categories
is established in
\cite[Theorem~5.2, Proposition~6.2 and
Corollary~6.1]{CHL2021}. These results give
dg quasi-equivalences.
Their scalar coefficient models and strong
interval conventions differ from the
plot-smooth bundle category used here.
The general plot-smooth prism proof is
\cite[Section~8]{ZengGlobalResolutions2026}.
We recall its relative-homotopy hypotheses
and prove the consequence for the specified
affine $\Phi,\Psi,H$.

\subsection{The normalized prism}

Let $X_0=Y_0=\y M$ and let
$H:X\times\Delta[1]\to Y$ have endpoints
$f=H(-,0)$ and $g=H(-,1)$, both the
identity on the vertex sheaf.
A homotopy \emph{relative to $c(\y M)$}
is constant on the whole subsheaf
$c(\y M)\times\Delta[1]$, not merely at
its degree-zero endpoints.
For $\sigma\in X_n$ put
\[
 h_i(\sigma)=H_{n+1}(s_i\sigma,\tau_i),
 \qquad 0\le i\le n,
\]
where $\tau_i:[n+1]\to[1]$ is zero
through $i$ and one thereafter.
For $\alpha\in C_{\mathrm{sm}}^{p,q}(Y;V,W)$
with $p\ge1$ define
\begin{equation}\label{eq:representation-prism}
 (\mathcal P_H\alpha)(\sigma)
   =\sum_{i=0}^{p-1}(-1)^i\alpha(h_i\sigma),
 \qquad \sigma\in X_{p-1}.
\end{equation}
Set $\mathcal P_H=0$ on simplicial degree zero.
It lowers total degree by one, without an
additional internal-degree sign.
The first endpoint of every prism simplex
is at time zero and the last at time one.
Since $f_0=g_0=\id_M$, the bundle types
in the sum coincide.

The prism satisfies the following identities
\cite[Lemma~8.1]{ZengGlobalResolutions2026}.
\begin{lem}[Normalized prism identities]
\label{lem:normalized-prism}
The operator is natural on smooth plots,
preserves positive-degree normalization, and
satisfies
\begin{align}
 \delta\mathcal P_H+\mathcal P_H\delta&=0,
 \label{eq:inner-prism}\\
 \mathcal P_H(\alpha\star\beta)
  &=(\mathcal P_H\alpha)\star g^*\beta
   +(-1)^{|\alpha|}f^*\alpha\star\mathcal P_H\beta.
 \label{eq:product-prism}
\end{align}
\end{lem}

\subsection{Finite gauges and all-degree naturality}

For a whole-unit-relative homotopy, the prism also gives the
following gauges
\cite[Lemma~8.2 and equation~(8.4)]{ZengGlobalResolutions2026}.
\begin{lem}[Pointed gauges]\label{lem:prism-gauge}
Assume now that $H$ is relative to $c(\y M)$.
For a representation $R=(V,r)$ set
$U_R=\mathcal P_H r$. This is a normalized,
closed degree-zero morphism
\[
 U_R:g^*R\longrightarrow f^*R,\qquad
 (U_R)_0=\id_V,
\]
with an actual closed inverse.
For a homogeneous morphism
$\alpha:R_V\to R_W$ of any total
degree $m\in\Z$, one has
\begin{equation}\label{eq:twisted-prism}
 D\mathcal P_H\alpha+\mathcal P_H D\alpha
   =U_W\star g^*\alpha-f^*\alpha\star U_V.
\end{equation}
Here the first differential is on
$\Hom_X(g^*R_V,f^*R_W)$; the differential
inside $\mathcal P_H D$ is on the
original Hom complex over $Y$.
\end{lem}

The hypothesis on the whole $c(\y M)$ makes the vertical
edge $H(s_0x,\tau_0)$ degenerate, so $(U_R)_0=1$.
Equality of endpoint object maps alone does not justify it.
For $V$ supported in $[a,b]$, write $U_R=1+u$.
The inverse is the finite formula
\cite[equation~(8.5)]{ZengGlobalResolutions2026}
\begin{equation}\label{eq:finite-prism-inverse}
 U_R^{-1}=\sum_{j=0}^{b-a}(-u)^{\star j}
\end{equation}
with the reverse representation type. Positive components
$u_p\in C^{p,-p}$ lower internal degree, so
$u^{\star N}=0$ for $N>b-a$.
The twisted identity applies to nonclosed morphisms and
negative total degrees as well.

\subsection{Raw homotopy invariance}

These identities yield homotopy invariance
\cite[Theorem~8.3]{ZengGlobalResolutions2026}.
\begin{thm}[Raw representation invariance]
\label{thm:raw-representation-invariance}
Let $f:X\to Y$ and $g:Y\to X$ fix the
common vertex sheaf $\y M$. Suppose there
are simplicial homotopies
\[
 H_Y:fg\Longrightarrow\id_Y,\qquad
 H_X:gf\Longrightarrow\id_X
\]
relative to $c(\y M)$ on both sides.
Then the actual pullback dg functor
\[
 f^*:\Rep_{\mathrm{sm}}^{\mathrm{perf}}(Y)
       \longrightarrow
       \Rep_{\mathrm{sm}}^{\mathrm{perf}}(X)
\]
is a dg quasi-equivalence. Its maps on
all Hom complexes are chain homotopy
equivalences. The result restricts to
the regular full subcategories and
requires no Kan hypothesis on $X$ or $Y$.
\end{thm}

\subsection{The actual affine pullback}

\begin{cor}[Global affine coefficient comparison]
\label{cor:affine-representation-invariance}
For the affine comparison of
Theorem~\ref{thm:affine-sdr}, pullback gives
a dg quasi-equivalence
\[
 \Phi^*:\Rep_{\mathrm{sm}}^{\mathrm{perf}}(\y K)
       \longrightarrow
       \Rep_{\mathrm{sm}}^{\mathrm{perf}}
                   (\Int(E_{\mathrm{aff}})),
\]
also on the regular full subcategories.
For objects $R,Q$ on $\y K$, write
\[
 i=\Phi^*,\qquad p=\Psi^*
\]
on their Hom complexes. Then
\begin{equation}\label{eq:affine-hom-retraction}
 pi=\id,\qquad
 D\mathcal P_H+\mathcal P_H D=\id-ip.
\end{equation}
Thus the Hom complex on full MC integration
retracts onto the image of the Hom complex
on the finite nerve.
\end{cor}

\begin{proof}
Because $H$ is over $T=\y K$, its prism
simplices satisfy
\[
 \Phi(h_i\sigma)=s_i\Phi(\sigma).
\]
Strict units and normalization give
$U_{\Phi^*R}=\id$: the zero component
is the unit, and every positive one
uses a degenerate structure component
of degree at least two.
Also $\mathcal P_H(\Phi^*\gamma)=0$
for a normalized morphism $\gamma$.

For an arbitrary cochain $\alpha$
between $\Phi^*R$ and $\Phi^*Q$,
the twisted prism identity therefore
becomes
\[
 D\mathcal P_H\alpha+\mathcal P_H D\alpha
        =\alpha-(\Psi\Phi)^*\alpha.
\]
The objects are fixed by $(\Psi\Phi)^*$,
since $\Phi\Psi=\id_T$, so this is an
identity on one Hom complex, not only
on cochains already in the image of $i$.
The equality $\Psi^*\Phi^*=\id$ gives
\eqref{eq:affine-hom-retraction}.
For any object $R'$ on $\Int(E_{\mathrm{aff}})$, the gauge
$U_{R'}:R'\to\Phi^*\Psi^*R'$ is closed
and invertible. This proves essential
surjectivity and completes the comparison.
Regularity follows as in the theorem.
\end{proof}

This proof uses specified inverse data, not
weak equivalence alone. The whole-plane
$\Psi,H$ need not restrict to an arbitrary
open. For each $O$, the range-preserving
$\Psi_O,H_O$ satisfy the same whole-unit-relative
hypotheses, so the same prism argument applies
to $\Phi_O$. Thus the comparison acts on the
raw coefficient categories. The subsequent
differentiation and reconstruction argument
identifies them with cohesive coefficients.

For fixed-base weak maps, use
\cite[Corollary~8.4]{ZengGlobalResolutions2026}
with the connective input of
Theorem~\ref{main:geometric-invariance}.
For the supplied global resolutions, the direct route is
\cite[Corollary~8.5]{ZengGlobalResolutions2026}.
Those inverse homotopies preserve the entire base simplex,
so their restriction and Hom-sheafification assertions apply
directly. This differs from the open-by-open choices of
affine $\Psi_O,H_O$. Hom-sheafification changes
the morphism complexes while retaining the existing global objects.
Together with Corollary~\ref{cor:affine-representation-invariance},
these results connect the raw representation categories of
the admissible global affine resolutions to that of $K$.
Effective descent instead concerns reconstructing objects from
local data; Section~\ref{sec:programme} discusses this and the
related coherence and holonomy questions.

\section{Smooth coefficient integration}
\label{sec:coefficient-integration}

The general smooth construction is
\cite[Theorem~9.2]{ZengGlobalResolutions2026}, using the universal
higher de Rham kernel of Arias Abad--Sch\"atz
\cite{AriasAbadSchaetz2013}. We recall the actual Taylor formulas
and signs needed below: these specify the functor whose affine
quasi-equivalence we will prove. The cited companion establishes
convergence, gauge equivariance and smooth parameter dependence.

\begin{thm}[Coefficient integration]\label{thm:coefficient-integration}
For every geometric fixed-base $\li$-algebroid $E$ in
Section~\ref{sec:polynomial}, there is a smooth strictly unital
$\ai$-functor
\[
 \mathsf{IN}_E:\Coh(B_E)\longrightarrow
 \Rep_{\mathrm{sm}}^{\mathrm{perf}}(\Int(E)).
\]
It retains the underlying graded bundle and its internal
differential, preserves the regular full subcategories, and
commutes with open restriction and fixed-base infinity-morphisms.
\end{thm}

This is \cite[Theorem~9.2]{ZengGlobalResolutions2026}.
It uses globally bounded finite-projective coefficients and
allows unbounded degrees of $E$. The construction depends on the
CE algebra and full MC simplices, independently of Nuiten's
connective mapping-space comparison. For the affine algebroid,
Theorem~\ref{thm:affine-open-rh} will establish quasi-equivalence.

\subsection{Full scalar extension and totalization}

For $\sigma\in\Int(E)_p(U)$, let
$x_\sigma:U\times\Delta^p\to M$ be its base map and put
$\Omega_\sigma=\Omega^\bullet(U\times\Delta^p/U)$.
The full MC equations give a dg algebra map
$\sigma^*:B_E\to\Omega_\sigma$.
Scalar extension of $M_V=V\otimes_A B_E$ uses
\[
 D_\sigma(m\otimes\omega)
 =D_Vm\otimes\omega+(-1)^{|m|}m\otimes d_\Delta\omega.
\]
It is a flat relative superconnection on the actual
$x_\sigma^*V$ and preserves Hom differentials and composition
\cite[Lemma~9.3]{ZengGlobalResolutions2026}.
Thus scalar extension retains the base map and anchor derivatives.

We use the forms-first braiding
\[
 \mathfrak b(v\otimes\omega)=(-1)^{|v||\omega|}\omega\otimes v.
\]
It gives
\[
 D^L(\omega v)=d_\Delta\omega\,v+(-1)^{|\omega|}\omega D^L(v),
 \qquad
 \phi^L(\omega v)=(-1)^{m|\omega|}\omega\phi^L(v)
\]
for a degree-$m$ homogeneous morphism. The right-module
convention remains the one in Section~\ref{sec:cohesive}.

This is \cite[Lemma~9.4 and equation~(9.1)]{ZengGlobalResolutions2026}.
\begin{lem}[All-degree target dictionary]\label{lem:in-sign-dictionary}
Let $\alpha$ have total degree $m$ in the raw target.
The transformation
\begin{equation}\label{eq:in-sign-dictionary}
 (T_m\alpha)_p=(-1)^{pm+p(p+1)/2}\alpha_p
\end{equation}
identifies its dg category with the convention of
\cite[Definition~2.2 and Remark~2.4, pp.~6--7]{AriasAbadSchaetz2013},
where object structures transform by
$F_p=(-1)^{p(p-1)/2}r_p$.
It preserves units, normalization and edge quasi-isomorphisms.
\end{lem}

\subsection{The kernel, matrix signs and twisting formulas}

Use $|sa|=|a|-1$, with suspended dg coderivation
\begin{equation}\label{eq:in-bar-convention}
 b_1(sa)=-s(da),\qquad
 b_2(sa\otimes sb)=(-1)^{|a|}s(ab).
\end{equation}
The unsuspended Taylor convention is
$\overline F_k=sF_k(s^{-1})^{\otimes k}$, including graded
tensor signs. These are
\cite[equation~(9.2)]{ZengGlobalResolutions2026}.
Both coderivation terms are the negatives of
\cite[Remark~3.1]{AriasAbadSchaetz2013}; negating both
source and target coderivations leaves the bar-map equation
unchanged.

For a smooth manifold $N$, write $\Omega=\Omega^\bullet(N)$
for its de Rham algebra and
$C_{\mathrm{norm}}=C^\bullet_{\mathrm{norm}}(S_{\mathrm{sm}}(N))$
for normalized smooth singular cochains. The scalar kernel below
is the cited natural $\ai$-map between these algebras. In the
present application it is evaluated on the simplex-direction
forms in $\Omega_\sigma$, with the external parameter held fixed.
Jointly smooth dependence on that parameter is supplied by the
stated all-derivative convergence and gluing results. No differential
in the external parameter is included in $d_\Delta$.

Let $\psi$ be the kernel of
\cite[Definition~3.23 and Theorem~3.25]{AriasAbadSchaetz2013}
from $(\Omega,-d,\wedge)$ to ordinary smooth singular cochains.
For scalar forms put $\tau(a)=(-1)^{|a|}a$ and fix the
uncurved $\ai$-morphism
\[
 \Gamma=\psi\tau:
 (\Omega,d,\wedge)\longrightarrow
 (C_{\mathrm{norm}},d_{\mathrm{full}},\cup).
\]
Thus
$\overline\Gamma_N(sa_1,\ldots,sa_N)
=(-1)^{\sum_i|a_i|}\psi_N(sa_1,\ldots,sa_N)$, with
$\overline\Gamma_N=s\Gamma_N(s^{-1})^{\otimes N}$.
The completed kernel is natural under smooth pullback.
For a $p$-form its unary term is
\[
 \Gamma_1(a)(\sigma)=\int_{\Delta^p}\sigma^*a.
\]
In form degree zero this is vertex evaluation.
Every higher $\Gamma_N$, $N>1$, vanishes if any input has
ordinary form degree zero. Its outputs are normalized
\cite[Proposition~3.26]{AriasAbadSchaetz2013}.

For matrices identify
$s(\Omega\otimes\End V)=(s\Omega)\otimes\End V$.
On elementary inputs $s\omega_a\otimes T_a$ with constant
graded matrices, the suspended extension collects the scalar
forms with sign
\begin{equation}\label{eq:in-matrix-sign}
 (-1)^{\sum_{a<b}|T_a|(|\omega_b|-1)},
\end{equation}
and multiplies the matrices in their original order.
The corresponding unsuspended formula is
\begin{equation}\label{eq:in-unsuspended-matrix-sign}
 \Gamma_N^V(\omega_1\otimes T_1,\ldots,\omega_N\otimes T_N)
 =(-1)^{\sum_{a<b}|T_a||\omega_b|}
   \Gamma_N(\omega_1,\ldots,\omega_N)\otimes T_1\cdots T_N.
\end{equation}
These are exactly
\cite[equations~(9.3)--(9.4)]{ZengGlobalResolutions2026},
the forms-first matrix extension of
\cite[Proposition~A.17]{AriasAbadSchaetz2013}.
For composable morphisms use the corresponding block of
the finite direct sum of coefficient bundles. This matrix
extension is understood in every formula below.

In a simplexwise frame write $D_i^L=d+A_i$.
For $\phi_a:V_{a-1}\to V_a$ of degree $m_a$, the actual
twisted Taylor expression is
\cite[equation~(9.5)]{ZengGlobalResolutions2026}
\begin{equation}\label{eq:in-twisted-taylor}
 \begin{split}
 \mathcal I_k(\phi_k,\ldots,\phi_1)
 =\sum_{j_0,\ldots,j_k\ge0}(-1)^{\chi}
 \Gamma_{k+\sum j_i}
 (A_k^{\otimes j_k},\phi_k,A_{k-1}^{\otimes j_{k-1}},
       \ldots,\phi_1,A_0^{\otimes j_0}),\\
 \chi=\sum_{a=1}^k(m_a-1)\sum_{\ell<a}j_\ell,
 \qquad k\ge1.
 \end{split}
\end{equation}
The suspended expression inserts $(sA_i)^{\otimes j_i}$
around $s\phi_k,\ldots,s\phi_1$ without an additional sign.
Only the diagonal $A_i$ are MC elements: the off-diagonal
morphism inputs need not be closed or MC.

For an object the formula is
\cite[equation~(9.6)]{ZengGlobalResolutions2026}
\begin{equation}\label{eq:in-object-formula}
 a_V=\sum_{N\ge1}\Gamma_N(A,\ldots,A),\qquad
 r_V=e+a_V,\qquad e_1=\id,\quad e_p=0\ (p\ne1).
\end{equation}
Here $e$ is a one-cochain in the temporary constant frame,
not a globally specified map between distinct base fibers.
The unit shift is governed by
\cite[equation~(9.7)]{ZengGlobalResolutions2026}:
\begin{equation}\label{eq:in-full-inner}
 d_{\mathrm{full}}\alpha
 =\delta\alpha+e_W\star\alpha
                -(-1)^{|\alpha|}\alpha\star e_V.
\end{equation}
The endpoint terms complete the inner differential to the
full coboundary, with final-face sign $(-1)^{p+1}$ on a
$p$-cochain. The displayed outputs are already in our
ordinary-cochain totalization. \emph{No further factor $T$
is applied to them.}

\subsection{Smoothness, degree control and naturality}

For fixed output simplex degree, Taylor arity, homogeneous
inputs and a compact subset of a local smooth family-frame
chart, the object and Taylor series converge locally uniformly
with every external-parameter derivative
\cite[Lemma~9.5]{ZengGlobalResolutions2026}.
In a nonzero higher kernel word every input has positive
ordinary form degree and
\begin{equation}\label{eq:in-surplus-degree}
 \sum_{\mathrm{inputs}}(\deg_{\mathrm{form}}-1)=p-1.
\end{equation}
This is \cite[equation~(9.8)]{ZengGlobalResolutions2026}.
One-form connection insertions can be arbitrarily numerous;
the general twisting series may therefore be infinite.
Structure inputs of form degree at least two lower internal
degree, whereas the distinguished morphism slots have their
specified arbitrary internal degrees. Matrix products
retain the input order in
\eqref{eq:in-unsuspended-matrix-sign}; this is the filtration
control used in finite affine reconstruction.

The outputs are intrinsic smooth two-ended bundle maps,
independent of temporary family frames
\cite[Lemma~9.6]{ZengGlobalResolutions2026}.
With endpoint identifications understood, their evaluation is
\cite[equation~(9.10)]{ZengGlobalResolutions2026}:
\begin{equation}\label{eq:in-taylor-definition}
 \bigl(\mathsf{IN}_{E,k}(\phi_k,\ldots,\phi_1)\bigr)_p(\sigma)
 =\mathcal I_k(\phi^L_{k,\sigma},\ldots,\phi^L_{1,\sigma})
                    (\id_{\Delta^p}).
\end{equation}
The $k$th Taylor map has degree $1-k$; this $p$-component
has internal degree $\sum_a|\phi_a|+1-k-p$.
Bounded coefficient amplitudes give finite cochain support
in each total degree, without changing the normal-convergence
assertion for the kernel words.

At vertices, $r_0=D^{[0]}$, the unary morphism is $\phi^{[0]}$,
and higher Taylor maps vanish. Strict units mean
$\mathsf{IN}_{E,1}(1)=1$ and
$\mathsf{IN}_{E,k}(\ldots,1,\ldots)=0$ for $k>1$.
Edges are chain isomorphisms of the vertex complexes.
For a fixed-base infinity-morphism $F:E\rightsquigarrow E'$,
\[
 \mathsf{IN}_E F^*=\Int(F)^*\mathsf{IN}_{E'}
\]
in every arity, under the canonical scalar-extension
identifications. Open restriction uses simplices whose
entire base image lies in the open. These are properties
of the actual functor in
\cite[Theorem~9.2]{ZengGlobalResolutions2026}.

\subsection{Postcomposition and the affine higher-isotropy test}

\begin{cor}\label{cor:in-postcomposition}
For the affine resolution the actual section in
Theorem~\ref{thm:affine-sdr} gives a global functor
\[
 \mathsf{IN}_{E_{\mathrm{aff}}}^{K}
 =\Psi^*\mathsf{IN}_{E_{\mathrm{aff}}}:
 \Coh(B_{E_{\mathrm{aff}}})\longrightarrow
 \Rep_{\mathrm{sm}}^{\mathrm{perf}}(\y K).
\]
It is a smooth strictly unital $\ai$-functor and
preserves regular coefficients.
\end{cor}

\begin{proof}
Postcomposition uses the already constructed dg functor
$\Psi^*$, with unchanged vertex complexes.
The global raw target change is justified by
Corollary~\ref{cor:affine-representation-invariance}.
\end{proof}

For each ordinary open $O$, postcomposition with the separately
constructed $\Psi_O$ likewise gives
$I_O=\Psi_O^*\mathsf{IN}_{E_O}$ on $\y K[O]$.
The sections are chosen separately for each open; strict
restriction-naturality is a property of the intrinsic forward
functor $\mathsf{IN}$ itself.

\begin{prop}[A coefficient detecting the affine two-form]
\label{prop:in-affine-probe}
Let $V=\R v_0\oplus\R v_1$ over $\R^2$ with
$|v_0|=0$, $|v_1|=1$, and
\[
 D(v_0)=0,\qquad D(v_1)=v_0\eta.
\]
Set $N(v_1)=v_0$, $N(v_0)=0$.
Its integrated representation has
\[
 r_0=0,\qquad r_1=1,\qquad
 r_2(\sigma)=\left(\int_{[012]}\beta\right)N,\qquad
 r_p=0\quad(p\ge3).
\]
On $K$ the two-component is exactly $-h_{012}N$,
including at the singular origin.
\end{prop}

\begin{proof}
The coefficient is flat because $Q\eta=0$ and $N^2=0$.
After full scalar extension and braiding its form is
$\beta\otimes N$. Every word with two or more such
inputs vanishes by $N^2=0$. The actual linear kernel
$\Gamma_1=\int$ therefore fixes the displayed sign
and normalization; the identity one-cochain supplies $r_1$.
The degree-three object equation is
\[
 (-h_{123}+h_{023}-h_{013}+h_{012})N=0,
\]
where the first and last terms come from the two
unit products, not the inner differential alone.
Degree four uses $N^2=0$. Since $\Phi\Psi=1$ and
$h=-\int\beta$, the formula on $K$ follows.
At the origin the exponential boundary $e^{hZ_0}$
is the identity, but this two-component still varies
linearly with the retained scalar $h$.
\end{proof}

MC alone would also allow $r_2=-c hN$ for any scalar $c$,
including zero. It is the linear integration kernel,
not that equation alone, which determines $c=1$.
For these affine functors, the all-open Hom comparison and
nonregular reconstruction establish quasi-equivalence in
Theorem~\ref{thm:affine-open-rh}. General-algebroid reconstruction
and effective descent are discussed in Section~\ref{sec:programme}.

\section{Range-preserving integration on arbitrary opens}
\label{sec:affine-open-geometry}

Let $O\subseteq\R^2$ be an ordinary open set and write
$E_O=E_{\mathrm{aff}}|_O$.  We construct the finite model and
the inverse-homotopy data separately on each $O$.  The interior
base maps of the global section in Theorem~\ref{thm:affine-sdr}
need not stay in a prescribed open set, so the construction
must preserve the range in $O$.

When the raw coefficient notation is applied to a simplicial
manifold, we suppress its levelwise Yoneda embedding:
$\Rep_{\mathrm{sm}}^{\mathrm{perf}}(K[O])$ means
$\Rep_{\mathrm{sm}}^{\mathrm{perf}}(\y K[O])$.

Write $\mathsf G=\Aff_+(1)$, with identity $e=(1,0)$, and use
the root coordinates of \eqref{eq:affine-root-coordinates}.
For a simplex of $K_n$, denote its vertices by
$m_i=(x_i,y_i)$.  All $y_i$ agree.  Define
\begin{equation}\label{eq:open-hull-model}
 K[O]_n=
 \left\{
 k\in K_n:
 [\min_i x_i,\max_i x_i]\times\{y_0\}\subseteq O
 \right\}.
\end{equation}
The notation refers to the whole closed horizontal hull, not
merely to membership of the vertices in $O$.

\begin{prop}[The open finite model]
\label{prop:open-finite-model}
The spaces $K[O]_n$ form a simplicial open submanifold of $K$.
They constitute a Lie $2$-groupoid with the same triangle and
tetrahedral formulas as $K$.  For $m=(x,y)\in O$, its source
fibre is
\begin{equation}\label{eq:open-arrow-domain}
 D_m=
 \left\{
 (a,b):a>0,\
 [\min\{x,ax+by\},\max\{x,ax+by\}]\times\{y\}\subseteq O
 \right\}.
\end{equation}
The total space
$D=\{(m,g):m\in O,\ g\in D_m\}$ is open, and every $D_m$
is convex and contains $e$.
\end{prop}

\begin{proof}
A compact horizontal segment contained in $O$ remains contained
in $O$ after a sufficiently small perturbation of its endpoints
and its height.  Thus \eqref{eq:open-hull-model} defines an open
subset of $K_n$.  Faces delete vertices and degeneracies repeat
them, so these operations preserve the hull condition.

For fixed $m=(x,y)$, let $I_y(m)$ be the component of
$O\cap(\R\times\{y\})$ containing $m$, identified with an open
interval in the $x$-coordinate.  Then
\[
 D_m=\{(a,b):a>0,\ ax+by\in I_y(m)\}.
\]
This is an intersection of an open half-plane with the inverse
image of an interval under an affine function.  It is convex
and contains $(1,0)$.  Openness of the total space follows from
compactness of the parameter interval in the map
\[
 (m,g,t)\longmapsto (1-t)m+t(gm),\qquad 0\le t\le1.
\]

The degree-one source and target maps are restrictions of
submersions to an open arrow manifold; they remain surjective
submersions because units are present.  For a horn of degree
at least two, its one-skeleton connects all its vertices.
The retained edges put those vertices in one horizontal
component of $O$.  Every filler in $K$ has the same vertices,
and hence belongs to $K[O]$.  The original degree-two horn
fibres are $\R$, and the higher horn fillers are unique by
Proposition~\ref{prop:affine-source-data}.  Restriction to the
open hull condition preserves these assertions.
\end{proof}

\subsection{A jointly smooth complete chart of the source fibres}

The following chart is used both in the open integration
section and in the relative de Rham contraction below.  Its
dependence on $O$ is fixed once and for all.

\begin{lem}[Complete based fibre charts]
\label{lem:open-based-chart}
There are jointly smooth diffeomorphisms
\[
 \psi_m:D_m\longrightarrow\R^2,\qquad \psi_m(e)=0.
\]
They are defined for every $m\in O$, including the origin when
it belongs to $O$.
\end{lem}

\begin{proof}
Translate the group coordinates by writing $z=(a-1,b)$.
Then $D$ is an open subset of $\R^4_{m,z}$ containing the zero
section over $O$.

Choose a countable family of nonnegative functions
$f_j\in C_c^\infty(D)$ whose positive sets cover $D$, and
extend them by zero to $\R^4$.  Choose $c_j>0$ such that
\[
 c_j\max_{|\alpha|\le j}\|\partial^\alpha f_j\|_\infty
 \le 2^{-j}.
\]
For every fixed multiindex $\alpha$, the differentiated tail
of $\sum_j c_jf_j$ converges uniformly.  Thus
\[
 F=\sum_{j=1}^{\infty}c_jf_j
\]
is smooth with all derivatives, is positive exactly on $D$,
and vanishes outside $D$.  Set
\[
 \widetilde\rho=\frac{F}{1+F}.
\]
Then $0\le\widetilde\rho<1$ and its positive locus is $D$.

Choose a locally finite partition of unity $\{\lambda_i\}$
on $O$ and constants $\varepsilon_i>0$ such that
\[
 \{m\}\times\overline B_{4\varepsilon_i}(0)\subseteq D
 \quad\text{for }m\in\operatorname{supp}\lambda_i.
\]
Put $\delta(m)=\sum_i\lambda_i(m)\varepsilon_i$.  At each
$m$, only finitely many terms are nonzero.  If
$\varepsilon_i$ is maximal among the active terms, then
$3\delta(m)\le3\varepsilon_i<4\varepsilon_i$; consequently
\[
 \{m\}\times\overline B_{3\delta(m)}(0)\subseteq D.
\]
The function $\delta$ is smooth and strictly positive.

Choose a smooth cutoff $\chi(m,z)$, depending on
$\|z\|^2/\delta(m)^2$, equal to one for
$\|z\|\le2\delta(m)$ and zero for
$\|z\|\ge3\delta(m)$.  Define
\[
 \rho(m,z)=\chi(m,z)+(1-\chi(m,z))\widetilde\rho(m,z),
 \qquad X(m,z)=\rho(m,z)z.
\]
The vector field is vertical, smooth on $O\times\R^2$,
positive radially on $D\setminus(O\times\{0\})$, zero outside $D$,
and equal to $z$ near the zero section.  The estimate
$\|X(m,z)\|\le\|z\|$ makes its flow $\varphi_t$ complete
for positive and negative time.

Each ray meets $D_m$ in an interval starting at zero.
On that interval the radial equation is
\[
 \dot r=\rho(m,r\omega)r.
\]
Backward flow tends to zero.  Forward flow tends to the
endpoint of the interval, finite or infinite: a positive
interior limit would contradict positivity of $\rho$.
A boundary point cannot be reached in finite time, since
it is stationary and uniqueness would then force the
original trajectory to be stationary.

For $v\ne0$, define
\[
 \Xi(m,v)=
 \varphi_{\log(\|v\|/\delta(m))}
 \left(m,\delta(m)\frac{v}{\|v\|}\right),
 \qquad
 \Xi(m,0)=(m,0).
\]
When $\|v\|\le2\delta(m)$, this formula is exactly
$\Xi(m,v)=(m,v)$, because the relevant trajectory remains
where $\rho=1$.  Thus it is smooth at the zero section.
On each ray it is strictly increasing and maps
$(0,\infty)$ onto the entire radial interval of $D_m$.
Its radial derivative is positive, and its angular
derivative has nonzero tangential part.  Hence it is a
bijective local diffeomorphism over $O$, with smooth
inverse.  Translating back from $z$ to $(a,b)$ gives
$\psi_m$.
\end{proof}

The chart determines the selected edge path
\begin{equation}\label{eq:open-selected-edge}
 \gamma_{m,g}(t)=\psi_m^{-1}(t\psi_m(g)),\qquad 0\le t\le1.
\end{equation}
It stays in $D_m$, depends jointly smoothly on $(m,g,t)$,
and extends smoothly past both endpoints.  The chart is chosen
separately for each open.

\subsection{Lifting with the current root}

We use the smooth face extension, Whitney--radial primitive,
and ordered split-locus extension of
Lemmas~\ref{lem:smooth-face-extension},
\ref{lem:degree-two-radial}, and
\ref{lem:split-locus-extension}.  These are linear operations
on the displayed form or function residuals.  The following
argument specifies how they interact with the moving root
and the open condition.

For a parameter sheaf $P$, a smooth family of simplex
$q$-forms means a morphism from $P$ to the sheaf
\[
 U\longmapsto\Omega^q(U\times\Delta^n/U),
\]
defined on plots. The sheaf $P$ may be nonrepresentable.
The ordered extension of such vector-valued residuals is the
construction of Lemma~\ref{lem:split-locus-extension}.

\begin{lem}[Range-preserving split-locus lifting]
\label{lem:open-split-lifting}
Let $P$ be a smooth parameter sheaf and let
$k:P\to K[O]_n$ have compatible prescribed full-MC faces
over its faces.  Suppose compatible full lifts are also
prescribed on finitely many split loci
$P_1,\ldots,P_N\subseteq P$, with retractions $r_j$
satisfying
\[
 r_j(P_i)\subseteq P_i\cap P_j\qquad(i<j).
\]
Then the data admit a jointly smooth full-MC lift over $k$
whose entire base simplex lies in $O$.  The parameter
retractions need not preserve $k$ or its root.
\end{lem}

\begin{proof}
For $n=0$, the lift is forced.  Suppose first that $n\ge2$,
and let $m_0=m_0(k)$.  Put $c=-h(k)$.  The tetrahedral
identity says $\delta c=0$.  Let $\beta_\partial$ be the
prescribed boundary two-form.  The baseline
\[
 \beta^0=
 W(c)+d\,\operatorname{Ext}_n
 \bigl(D_\partial(\beta_\partial-W(c)|_\partial)\bigr)
\]
is closed, has the prescribed boundary, and has triangle
periods $c$, by Lemma~\ref{lem:closed-period-extension}.

On $P_j$, let $\beta_j$ be the prescribed two-form and set
\[
 \gamma_j=\beta_j-\beta^0|_{P_j},\qquad
 \eta_j=D_n\gamma_j.
\]
The form $\gamma_j$ is closed and has zero triangle
periods, so $d\eta_j=\gamma_j$.  Face naturality of $D_n$
makes $\eta_j$ vanish on the whole boundary.  Extend the
compatible $\eta_j$ by the ordered formula of
Lemma~\ref{lem:split-locus-extension}, within the linear
space of one-forms vanishing on the boundary, and denote
the extension by $\eta$.  Then
\[
 \beta=\beta^0+d\eta
\]
recovers all prescribed two-forms and faces.  Its triangle
periods remain $c$: every edge is contained in the boundary,
so the exact correction has zero triangle integrals.

Set $u=H_{v_0}\beta$.  Then $du=\beta$ and $u$ vanishes
on every root edge.  Flatten each prescribed face by
\[
 \mathsf A_i'=\mathsf A_i-Z_{\sigma_i}u|_{F_i}.
\]
Lemma~\ref{lem:affine-flat-correction} gives compatible
flat potentials on faces containing $v_0$, normalized to
one there.  On the opposite face use the potential with
value $a_1=g_{01}(k)$ at $v_1$.  Its compatibility at $v_2$
is the exact calculation
\[
 \int_{[12]}u=-h_{012},\qquad
 g_0(v_2)=g_{12}a_1\exp(h_{012}Z_{m_0})=a_2.
\]
Uniqueness of the flat potential then gives agreement on
all face intersections, including the opposite-face
intersection in degree two.

The resulting boundary potential $g_\partial$ satisfies
$\sigma_\partial=g_\partial m_0$.  Every prescribed face
has connected base image in a horizontal component of $O$
containing its vertices.  Those vertices lie in the
component of $m_0$ by the hull condition.  Therefore
$g_\partial$ takes values in $D_{m_0}$.

Extend the chart-valued boundary function, using the
current root:
\[
 g^0=
 \psi_{m_0}^{-1}
 \bigl(\operatorname{Ext}_n(\psi_{m_0}(g_\partial))\bigr).
\]
On $P_j$, flatten the prescribed lift with its recovered
$u$ and take its normalized potential $g_j$.
It agrees with $g_\partial$ on the boundary and lies in
$D_{m_0}$.  Thus
\[
 w_j=\psi_{m_0}(g_j)-\psi_{m_0}(g^0)|_{P_j}
\]
is a compatible vector-valued residual vanishing on the
boundary.  Extend these residuals by the ordered split-locus
formula, obtaining $w$, and set
\[
 g=\psi_{m_0}^{-1}\bigl(\psi_{m_0}(g^0)+w\bigr).
\]
The retractions act on the vector residuals only.  The
inverse chart in this formula always uses the current
$m_0(k)$, even if a retraction changes the target simplex.
Since the chart has range all of $\R^2$, $g$ remains in
$D_{m_0}$.

Define
\[
 \sigma=gm_0,\qquad
 \mathsf A=dg\,g^{-1}+Z_\sigma u.
\]
Reverse flat correction gives all the full-MC equations
\eqref{eq:affine-full-mc}.  The base lies in $O$ because
$g\in D_{m_0}$.  Root-edge holonomies are the prescribed
$a_i$, and triangle periods are $-h_{ijk}$.  The root
coordinate formulas therefore give $\Phi(\sigma,\mathsf A,
\beta)=k$ exactly.  The construction recovers every face
and split-locus lift.

For $n=1$, use $\beta=u=0$ and the baseline
$g^0(t)=\psi_{m_0}^{-1}(t\psi_{m_0}(a_1))$.
Each prescribed interval lift has a normalized flat
potential.  Extend its chart-valued difference from $g^0$
within the functions vanishing at the two endpoints, and
reconstruct by $\sigma=gm_0$, $\mathsf A=dg\,g^{-1}$.
This gives the same conclusion.

All constructions use smooth pullbacks, finite linear
operations, the fixed smooth charts, and compact-interval
integration.  They commute with pullback of external plots
and preserve joint smoothness through the closed faces.
\end{proof}

\begin{thm}[The actual open-wise comparison]
\label{thm:affine-open-comparison}
For every ordinary open $O$, the holonomy and period map
restricts to
\[
 \Phi_O:\Int(E_O)\longrightarrow\y K[O].
\]
For every smooth parameter manifold $U$, $\Phi_O(U)$ is a
trivial Kan fibration.  There are a smooth simplicial section
$\Psi_O$ and a smooth simplicial homotopy
\[
 H_O:\Int(E_O)\times\Delta[1]\longrightarrow\Int(E_O),
 \qquad
 H_O(-,0)=\Psi_O\Phi_O,\quad H_O(-,1)=\id,
\]
such that
\[
 \Phi_O\Psi_O=1,\qquad
 \Phi_OH_O=\Phi_O\pr,\qquad
 H_O|_{\Psi_O(K[O])}=\pr.
\]
In particular, $H_O$ is relative to the whole
$c(\y O)$.
\end{thm}

\begin{proof}
For a full-MC simplex in $O$, the equation $dy=0$ puts its
connected base image in one horizontal component of $O$.
Its vertices therefore satisfy the hull condition, so
$\Phi$ restricts to $\Phi_O$.

Lemma~\ref{lem:open-split-lifting}, without additional
split loci, solves every relative boundary problem.
In degree zero the map is the identity.  Thus each
$\Phi_O(U)$ has the right lifting property against all
boundary inclusions and is a trivial Kan fibration.

For the smooth section, induct on simplicial degree.
At degree $n$, the prescribed faces are
$\Psi_{O,n-1}d_i$.  On the degeneracy locus
$B_j=\im s_j$, prescribe
\[
 s_j\Psi_{O,n-1}d_{j+1}.
\]
The retractions $s_jd_{j+1}$ preserve earlier degeneracy
loci.  On $B_i\cap B_j$, $i<j$, both prescriptions are
the common double degeneracy
\[
 s_is_{j-1}\Psi_{O,n-2}
 =
 s_js_i\Psi_{O,n-2}
\]
of the uniquely determined lower simplex.  These are the
identities of Lemma~\ref{lem:ordered-degeneracy-loci}.
The split-locus lifting lemma therefore supplies
$\Psi_{O,n}$ with all required faces and degeneracies.

Construct the homotopy in the same way, degreewise for
each label $\tau:[n]\to[1]$.  Constant labels give the
two endpoints.  For a nonconstant label prescribe the
already constructed faces and first impose the fixed
locus $\Psi_O(K[O]_n)$, retracted onto by $\Psi_O\Phi_O$.
For every $j$ with $\tau(j)=\tau(j+1)$, impose the
degeneracy locus with prescription
\[
 s_jH_{O,n-1}(d_{j+1}\sigma,d_{j+1}\tau).
\]
On a double-degeneracy intersection, the two reduced
labels agree because they factor through the same
surjection
$\varsigma_{j-1}\varsigma_i=\varsigma_i\varsigma_j$.
The two prescriptions are consequently equal by the
simplicial identities.  On the fixed-locus intersections,
they agree by the lower-degree fixedness condition.
Each $s_jd_{j+1}$ preserves the fixed locus, so the fixed locus
can be placed first in the ordered split-locus construction.
Lemma~\ref{lem:open-split-lifting} completes the induction.

Every lift is taken over the current target and remains
inside $O$.  Labels on disconnected test manifolds are
locally constant, so the constructions glue on such tests.
Finally, $\Psi_O$ is the identity in degree zero and
preserves degeneracies; hence its image contains the
whole unit-generated subsheaf, which $H_O$ fixes.
\end{proof}

\begin{cor}[Raw coefficients on the two actual sources]
\label{cor:affine-open-raw-comparison}
The actual pullback dg functors along $\Phi_O$ and $\Psi_O$
are quasi-inverse dg quasi-equivalences on globally bounded
finite-projective raw coefficients.  Their maps on all Hom
complexes are chain-homotopy equivalences.  The body cohomology
ranks may vary.
\end{cor}

\begin{proof}
Apply Theorem~\ref{thm:raw-representation-invariance} to
the inverse-homotopy data of
Theorem~\ref{thm:affine-open-comparison}.  Its hypothesis is
the homotopy relative to the entire $c(\y O)$, not only
agreement of the endpoint vertex maps.
\end{proof}

\begin{rem}\label{rem:open-choice-scope}
The canonical maps $\Phi_O$ are compatible with restriction.
The charts $\psi_m$, sections $\Psi_O$, and homotopies $H_O$
are chosen separately on each $O$; strict compatibility with
inclusions of opens is an additional condition on such choices.
The empty-open assertions are vacuous.
\end{rem}

\section{Actual differentiation, its integration retraction,
and edge rectification}
\label{sec:affine-differentiation}

Fix an ordinary open $O$, and put
\[
 A=C^\infty(O),\qquad G=K[O],\qquad
 B_O=A[\theta_1,\theta_2,\eta].
\]
The differential on $B_O$ is the restriction of
\eqref{eq:intro-affine-ce}.  The coefficient categories
retain the globally bounded finite-projective hypotheses
of Sections~\ref{sec:cohesive} and \ref{sec:representations}.
Write
\[
 I=\Psi_O^*\mathsf{IN}_{E_O}.
\]
We first construct differentiation on the actual raw
category, including every Hom degree.

\subsection{The scalar normalized quotient}

For $n>0$, let
\[
 S_k=s_{k-1}(G_{n-1})\subseteq G_n,\qquad 1\le k\le n,
\]
and write $I(S_k)$ for the ideal of smooth functions on $G_n$
vanishing on $S_k$. Define
\begin{equation}\label{eq:open-differentiating-ideal}
 C_N^n=\bigcap_{k=1}^n I(S_k),\qquad
 J_n=\sum_{k=1}^n\bigl(C_N^n\cap I(S_k)^2\bigr),
 \qquad
 \widehat J=J+\delta_{\mathrm{full}}J.
\end{equation}
Here $J=\bigoplus_{n\ge0}J_n$, and $\delta_{\mathrm{full}}$
is the full simplicial coboundary. In degree zero,
$C_N^0=A$ and $J_0=0$. The product is the ordinary
front/back cup product.

\begin{prop}[The affine scalar differentiation]
\label{prop:open-scalar-differentiation}
Write $\mathfrak a=a_{01}$ and $b=b_{01}$ for the dilation
and translation coordinates of an edge $g_{01}$, and
$h_{\mathrm{geom}}=h_{012}$ for the triangle coordinate
of \eqref{eq:affine-surface-coordinate}.
The normalized DGA quotient of
\cite[Theorem~3.3.5 and Definition~3.3.6]{CabreraDelHoyo2026}
gives a unital DGA epimorphism
\[
 \operatorname{ve}:C_N^\bullet(G)\longrightarrow B_O
\]
with
\begin{equation}\label{eq:open-ve-generators}
 \operatorname{ve}(\log\mathfrak a)=\theta_1,\qquad
 \operatorname{ve}(b)=\theta_2,\qquad
 \operatorname{ve}(h_{\mathrm{geom}})=-\eta.
\end{equation}
Its target differential is exactly
\[
 Qf=f_x(x\theta_1+y\theta_2),\quad
 Q\theta_1=-y\eta,\quad
 Q\theta_2=\theta_1\theta_2+x\eta,\quad Q\eta=0.
\]
\end{prop}

\begin{proof}
The normalized infinitesimal bundles have ranks two in
degree one and one in degree two.  In degree two, fixing
the two positive faces fixes both root edges, leaving
the $h$-coordinate.  Above degree two the positive horn
fillers are unique.  The quotient theorem therefore makes
the graded algebra semi-free on the three generators in
\eqref{eq:open-ve-generators}.

For $f\in A$,
\[
 \delta_{\mathrm{full}}f(g,m)=f(gm)-f(m).
\]
Its first normal derivative is
$f_x(x\log\mathfrak a+yb)$, giving the formula for $Qf$.

Put $\ell=\log\mathfrak a$.  The first matrix entry of
\[
 g_{02}=g_{12}g_{01}\exp(hZ_{m_0})
\]
gives
\[
 \ell_{02}=\ell_{12}+\ell_{01}-yh,
 \qquad \delta_{\mathrm{full}}\ell=yh.
\]
Hence $Q\theta_1=-y\eta$.

For the second entry, with
$\varphi(t)=\int_0^1e^{-st}\,ds$, one obtains
\[
 \delta_{\mathrm{full}}b
 =(1-a_{12})b_{01}
 -a_{12}a_{01}xh\varphi(hy).
\]
Modulo $J_2$, the first term is $-b_{01}\ell_{12}$:
the difference contains $\ell_{12}^2b_{01}$ times a
smooth function.  For the second term use
\[
 a_{12}a_{01}\varphi(hy)-1
 =(a_{12}-1)a_{01}\varphi(hy)
 +(a_{01}-1)\varphi(hy)+\varphi(hy)-1.
\]
After multiplication by $h$, each summand is normalized
and lies in one of the required squared degeneracy ideals.
The last assertion follows from
$\varphi(t)-1=t\psi(t)$ with $\psi$ smooth; it is not
division by $y$.  Thus
\[
 \delta_{\mathrm{full}}b
 \equiv -b_{01}\ell_{12}-xh\pmod{J_2}.
\]
Graded commutativity of the full quotient by $\widehat J$
gives
\[
 Q\theta_2=-\theta_2\theta_1+x\eta
 =\theta_1\theta_2+x\eta.
\]
Finally, the tetrahedral equation is
$\delta_{\mathrm{full}}h=0$, so $Q\eta=0$.
All calculations are local near units and hence apply
on every $O$ by the germ-locality of the quotient
\cite[Lemma~3.1.5]{CabreraDelHoyo2026}.
\end{proof}

The quotient is by $J+\delta_{\mathrm{full}}J$, not merely
by adjacent-degeneracy squares.  The additional
allowed-transposition relations are part of its DGA
structure.

\begin{prop}[Scalar Van Est on every open]
\label{prop:open-scalar-van-est}
The map $\operatorname{ve}$ is a quasi-isomorphism in
every degree.
\end{prop}

\begin{proof}
The degree-one last-horn fibres are target fibres, which
inversion identifies with the convex domains $D_m$ of
\eqref{eq:open-arrow-domain}.  A degree-two last horn
specifies $g_{12}$ and $g_{02}$.  For every $h\in\R$,
the unique remaining edge is
\[
 g_{01}=g_{12}^{-1}g_{02}\exp(-hZ_{m_0}).
\]
The vertices already lie in one horizontal component,
so this filler remains in $K[O]$.  The fibre is therefore
$\R$, including at the origin.  Higher last-horn fibres
are points.  The last-horn criterion and higher Van Est
theorem of
\cite[Lemma~4.3.3 and Theorem~4.3.4]{CabreraDelHoyo2026}
apply in every finite degree.
\end{proof}

\subsection{Matrices, arbitrary Hom degrees, and projective bodies}

In a local homogeneous frame let $e$ be the temporary
identity one-cochain, with $e_1=1$ and all other components
zero.  It is not normalized.  Directly from the raw product,
\[
 \delta_{\mathrm{in}}e+e\star e=0,
\]
and, for a cochain of bidegree $(p,q)$,
\begin{equation}\label{eq:open-unit-shift}
 \delta_{\mathrm{full}}\phi
 =
 \delta_{\mathrm{in}}\phi
 +e\star\phi-(-1)^{p+q}\phi\star e.
\end{equation}
The final endpoint coefficient is
$-(-1)^{p+q}(-1)^q=(-1)^{p+1}$.

\begin{thm}[Actual strict dg differentiation]
\label{thm:open-dg-differentiation}
There is a restriction-natural strict dg functor
\[
 F=\operatorname{Diff}_O:
 \Rep_{\mathrm{sm}}^{\mathrm{perf}}(G)
 \longrightarrow \Coh(B_O).
\]
It retains the original graded bundle and its body
differential.  In a local homogeneous frame its formulas are
\begin{equation}\label{eq:open-diff-formula}
 a_R=r_R-e,\qquad
 FR=Q+\mu_R,\quad \mu_R=\operatorname{ve}(a_R),
 \qquad
 F(\phi)=\operatorname{ve}(\phi),
\end{equation}
where differentiation is entrywise.
These formulas apply to nonclosed morphisms in every
integer degree and retain every mixed and terminal
Maurer--Cartan equation.
\end{thm}

\begin{proof}
The matrix product on normalized cochains is
\[
 (cX)(dY)=(-1)^{|X||d|}(c\smile d)XY.
\]
It is precisely the forms-first version of the raw product.
Writing $r_R=e+a_R$, equation \eqref{eq:open-unit-shift}
turns the raw MC equation into
\[
 \delta_{\mathrm{full}}a_R+a_R\star a_R=0.
\]
For a homogeneous morphism of total degree $m$,
\[
 D\phi
 =
 \delta_{\mathrm{full}}\phi
 +a_S\star\phi-(-1)^m\phi\star a_R.
\]
Apply the scalar DGA map entrywise.  It gives
\[
 Q\mu_R+\mu_R^2=0
\]
and
\[
 F(D\phi)
 =
 QF(\phi)+\mu_SF(\phi)-(-1)^mF(\phi)\mu_R.
\]
This is the cohesive Hom differential in forms-first order.
The fixed braiding of
Section~\ref{sec:coefficient-integration} returns the
signlessly right-linear convention.  Products and identity
zero-cochains are preserved strictly.

If old coordinates equal $g$ times new coordinates, the
actual endpoint frame transformation is
\[
 a_R^g=g^{-1}a_Rg+g^{-1}\delta_{\mathrm{full}}g.
\]
Differentiation gives
\begin{equation}\label{eq:open-diff-frame-law}
 \mu_R^g=g^{-1}\mu_Rg+g^{-1}Qg.
\end{equation}
Thus the local operators glue to superconnections on the
original bundle.

For finite homogeneous presentations
$V_R=\im\Pi_R$ and $V_S=\im\Pi_S$, raw Hom cochains are
endpoint corners.  Degree-zero multiplication is cup
multiplication at the corresponding endpoint, so
\begin{equation}\label{eq:open-diff-corner}
 F\bigl(\Pi_S(v_0)\phi\Pi_R(v_p)\bigr)
 =
 \Pi_SF(\phi)\Pi_R.
\end{equation}
This proves global finite-projective compatibility directly
from the endpoint projectors.

The degree-zero component of the scalar quotient is the
identity, so body differentials are unchanged.  The single
identity $Q\mu_R+\mu_R^2=0$ includes all curvature
coefficients, including the last possible one allowed by
the body amplitude.  The construction therefore retains the
full superconnection on the given body.
\end{proof}

\begin{lem}[Global graded Hom lifting]
\label{lem:open-graded-hom-lifting}
For every ordered pair of actual raw objects,
\[
 F_{R,S}:
 \Hom^\bullet_{\mathrm{raw}}(R,S)
 \longrightarrow
 \Hom^\bullet_{\Coh(B_O)}(FR,FS)
\]
is degreewise surjective.
\end{lem}

\begin{proof}
Lift a scalar ordered monomial by
\[
 f\theta_1^{\varepsilon_1}\theta_2^{\varepsilon_2}\eta^k
 \longmapsto
 f\smile\ell^{\smile\varepsilon_1}
 \smile b^{\smile\varepsilon_2}
 \smile(-h_{\mathrm{geom}})^{\smile k},
 \qquad \varepsilon_i\in\{0,1\}.
\]
Proposition~\ref{prop:open-scalar-differentiation} shows
that this is a graded lift.  In finite homogeneous
presentations, lift every ambient matrix entry and then
apply both endpoint projectors.
Equation~\eqref{eq:open-diff-corner} gives the prescribed
cohesive corner.  The sums are finite in each total
degree because the internal Hom degrees are bounded.
The lifts are global, smooth, and normalized.  This construction
provides a section at the level of graded vector spaces.
\end{proof}

\subsection{Literal differentiation of the actual integration}

The geometric weighted function estimate of
Lemma~\ref{lem:open-weighted-functions} will annihilate higher
normal-weight terms in the retraction calculation.

\begin{lem}[The two-form of the selected open section]
\label{lem:open-section-whitney}
For the section constructed in
Theorem~\ref{thm:affine-open-comparison},
\[
 \beta_{\Psi_O(k)}=W(-h(k))
\]
in every simplicial degree.  Its degree-one potential is
the path \eqref{eq:open-selected-edge}.
\end{lem}

\begin{proof}
Whitney forms commute with every ordinal simplex map.
Indeed, if $\lambda$ is such a map, then
$\lambda^*t_j=\sum_{\lambda(i)=j}t_i$; substitution in
the alternating formula for the Whitney form gives
\[
 \lambda^*\omega_{j_0\ldots j_r}
 =
 \sum_{\substack{i_0<\cdots<i_r\\
                 \lambda(i_a)=j_a}}
 \omega_{i_0\ldots i_r}.
\]
Repeated image indices have zero alternating coefficient.
This proves naturality for degeneracies as well as faces.

Induct on the section degree.  In degrees zero and one
both two-forms vanish.  If the assertion holds on the
boundary of an $n$-simplex, the two-form baseline in
Lemma~\ref{lem:open-split-lifting} is already $W(-h)$:
its boundary residual is zero.  On every prescribed
degeneracy locus, the two-form is also $W(-h)$ by
naturality.  Thus all split-locus two-form residuals are
zero, and the linear ordered extension leaves them zero.
The subsequent potential extension does not change
$\beta$.

In degree one, the sole degeneracy prescription is the
identity edge.  Its potential and the radial baseline
are both constantly $e$, so its residual is zero.
The selected potential is therefore exactly
\eqref{eq:open-selected-edge}.
\end{proof}

\begin{thm}[Literal integration retraction]
\label{thm:open-literal-retraction}
For every ordinary open $O$,
\[
 \boxed{FI=\id}
\]
as strictly unital $\ai$-functors.  This is literal on
objects and unary morphisms, on the original bundles;
all higher Taylor maps of $FI$ vanish.
\end{thm}

\begin{proof}
Use root group coordinates
$u_i=(\log\mathfrak a_i,b_i)$ of weight one and rooted
triangle coordinates $H_{ij}$ of weight two, with the
equivalent adjacent coordinates of
Lemma~\ref{lem:open-weighted-functions}.
Weights concern smooth normal Taylor germs at the whole
unit section; they do not require analyticity.

By Lemma~\ref{lem:open-section-whitney},
$\beta=W(-h)$ has weight at least two.  Its simplex-radial
primitive $u=H_{v_0}\beta$ has the same bound.
At a totally degenerate unit, the selected MC simplex is
constant.  Its normalized flat potential is therefore
identically one.  Smoothness gives
$g-1,d_\Delta g$ weight at least one.  The reconstruction
formula
\[
 \mathsf A= d_\Delta g\,g^{-1}+Z_\sigma u
\]
then shows that both MC one-forms have weight at least
one.  Consequently the pullback of a homogeneous scalar
CE element of degree $d$ has normal weight at least $d$.

Write $\Gamma_N^\Psi$ for the actual scalar kernel pulled
back along the chosen full-MC section $\Psi_O$.
For this kernel, an output in
simplicial degree $n$ satisfies
\[
 n=\sum_i d_i+1-N,
\]
where $d_i$ are the ordinary input form degrees.
For $N>1$, every term with a zero-form input vanishes.
In every remaining term with $N>1$, all $d_i$ are positive;
fixed simplex pullbacks, products, and compact integrations then give
normal weight at least
\[
 \sum_i d_i=n+N-1\ge n+1.
\]
The output is normalized.  Lemma~\ref{lem:open-weighted-functions}
therefore gives
\[
 \operatorname{ve}\Gamma_N^\Psi=0\qquad(N>1).
\]
There is no nonzero higher output in degree zero.

Let $F_N=\operatorname{ve}\Gamma_N^\Psi$.
Vertex evaluation gives $F_1(f)=f$.  On the selected
edge,
\[
 \int\alpha_1=\log\mathfrak a,\qquad
 \int\alpha_2=b-\int b(t)\,d\log\mathfrak a(t).
\]
The last integral has weight at least two and vanishes
under degree-one differentiation.  Thus
$F_1(\theta_i)=\theta_i$.  The exact period
$\int\beta=-h$ gives $F_1(\eta)=\eta$.
Since all higher $F_N$ vanish, the $\ai$ identities make
$F_1$ a DGA map.  Its values on $A,\theta_1,\theta_2,\eta$
give $F_1=1$, including all powers of the even generator.

Now use the actual matrix extension
\eqref{eq:in-unsuspended-matrix-sign} and the twisted
formula \eqref{eq:in-twisted-taylor}.  In every word of
length greater than one, either a zero-form slot kills
the word or the preceding normal-weight argument kills
its differentiated scalar entries.  For an object,
only the bare unary term of $r-e$ survives, recovering
the entire original superconnection, including its
body differential.  For a unary morphism, only the
uninserted unary word survives.  For a higher Taylor
map, every word has length at least two and vanishes
after differentiation.  These formulas use the original
target totalization throughout.

The general kernel convergence statement in
Section~\ref{sec:coefficient-integration} is local uniform
convergence with every smooth parameter derivative.
The quotient in degree $n$ uses a finite normal jet,
so these derivatives can be taken termwise in the
actual infinite connection series.  Equivalently, the
ordered-time estimates have a factorial denominator,
with only polynomial growth from the finitely many
distinguished slots and derivative allocations.
Thus the vanishing applies to the actual sums, not
only to individual formal words.

Finally, choose a homogeneous frame near a unit.
Compactness of the simplex and joint smoothness of the
chosen $\Psi_O$ allow the normal neighborhood to be
shrunk so that every selected base simplex lies in that
frame patch.  This uses the restriction of $\Psi_O$,
not a separately chosen section for the smaller patch.
The endpoint frame laws and
\eqref{eq:open-diff-corner} glue the identities on the
original finite-projective bundles.
\end{proof}

\subsection{Actual invertible-edge rectification}

The original raw edges need only be fibrewise
quasi-isomorphisms.  The following replacement makes
them actual chain isomorphisms by a closed raw gauge.

\begin{thm}[Closed edge rectification]
\label{thm:open-edge-rectification}
Every raw object $R$ on $G$ admits a normalized closed
isomorphism
\[
 E_R:R^{\mathrm e}\longrightarrow R,\qquad (E_R)_0=1,
\]
where the edges of $R^{\mathrm e}$ are the invertible
last-to-first ODE transports of the connection
differentiated from $r_{R,1}$ along
\eqref{eq:open-selected-edge}. The inverse of $E_R$ is given
by a finite series. No regularity assumption is required.
\end{thm}

\begin{proof}
Fix an edge $g:m\to gm$ and put
$\gamma(t)=\psi_m^{-1}(t\psi_m(g))$, $m_t=\gamma(t)m$.
For small $\varepsilon$, the increment
\[
 c_{t,\varepsilon}
 =\gamma(t+\varepsilon)\gamma(t)^{-1}
\]
and the $h=0$ triangle formed with $\gamma(t)$ and
$\gamma(t+\varepsilon)$ are admissible: their vertices
lie in one horizontal component of $O$.

In local frames along the path, set
\[
 A_t=\left.\partial_\varepsilon\right|_0
 r_{R,1}(c_{t,\varepsilon},m_t).
\]
Under a frame change $u(t)$ it transforms as
$u^{-1}A_tu+u^{-1}\dot u$.  Thus the equation
\[
 \partial_tT_{s,t}=T_{s,t}A_t,\qquad T_{s,s}=1
\]
defines intrinsic smooth invertible transport
$T_{s,t}:V_{m_t}\to V_{m_s}$.
Differentiating the edge chain-map equation gives
\[
 \dot d_t=d_tA_t-A_td_t,\qquad
 d_sT_{s,t}=T_{s,t}d_t,
 \quad d_t=r_{R,0}(m_t).
\]

Write $R_t=r_{R,1}(\gamma(t),m)$ and define
\[
 K_t=\left.\partial_\varepsilon\right|_0
 r_{R,2}
 \bigl(\gamma(t),c_{t,\varepsilon},\gamma(t+\varepsilon);0\bigr).
\]
At $\varepsilon=0$ the triangle is degenerate and its
$r_2$ is zero.  The derivative is therefore independent
of a temporary identification of the moving endpoint.
The triangle MC equation yields
\[
 \dot R_t-R_tA_t=d_mK_t+K_td_t.
\]
The bundle-valued integral
\[
 H_t=\int_0^tK_sT_{s,t}\,ds:V_{m_t}\to V_m
\]
satisfies $\dot H_t=K_t+H_tA_t$.
Using the formula for $\dot d_t$, both
$R_t-T_{0,t}$ and $d_mH_t+H_td_t$ solve the same
inhomogeneous ODE with zero initial value.  Hence
\[
 R_t-T_{0,t}=d_mH_t+H_td_t.
\]

Put $H(g,m)=H_1$ and define $E=1+H$, with $H$ only
in simplicial degree one and internal degree minus one.
At a unit, $H=0$.  If the body width is $w$, then
$H^{\star(w+1)}=0$ and
\[
 E^{-1}=\sum_{j=0}^{w}(-H)^{\star j}.
\]
Define
\[
 r^{\mathrm e}
 =
 E^{-1}\star(\delta_{\mathrm{in}}E+r_R\star E).
\]
The raw derivation identities prove its MC equation and
the closedness of $E:R^{\mathrm e}\to R$.
Normalization follows from
\[
 r^{\mathrm e}=r_R+E^{-1}\star D_RE:
\]
Here $D_R$ is the differential on $\Hom_{\mathrm{raw}}(R,R)$.
It preserves normalized cochains, which are closed under the
required products. The body is unchanged. In degree one,
\[
 r^{\mathrm e}_1
 =r_{R,1}-d_mH-Hd_{gm}
 =T_{0,1}.
\]
Thus the new edges are actual chain isomorphisms.
The finite inverse is closed by the Hom Leibniz rule.
All operations are intrinsic, smooth in parameters,
and compatible with the original projective bundles.
\end{proof}

\begin{rem}\label{rem:open-rectification-differentiation}
Strict dg differentiation sends $E_R$ to a closed
invertible cohesive map, with
$F(E_R^{-1})=F(E_R)^{-1}$.  Its body is the identity,
but $F(E_R)$ itself need not be the identity.
These actual differentiated gauges will be retained
when comparison is transferred back from rectified
objects.
\end{rem}

\section{Relative Hom transfer and normal-symbol calibration}
\label{sec:affine-relative-comparison}

Fix an ordinary open $O$, put $A=C^\infty(O)$ and
$G=K[O]$, and let $R,S$ be raw coefficients. Write
$\partial_R=r_{R,0}$ and $\partial_S=r_{S,0}$ for their body
differentials. If their bodies occupy $[a_R,b_R]$ and $[a_S,b_S]$, set
\[
 \ell=a_S-b_R,\qquad u=b_S-a_R.
\]
The construction below retains the body differential.
Internal degree will control triangular operators, not
a filtration by subcomplexes.

\subsection{A finite perturbation lemma}

\begin{lem}[Finite perturbation with a retained small differential]
\label{lem:open-finite-perturbation}
Suppose
\[
 (C_0,b)\mathrel{\mathop{\rightleftarrows}^{i}_{p}}(C_1,d)
\]
has a degree-minus-one map $h$ satisfying
\[
 pi=1,\quad di=ib,\quad pd=bp,\quad
 dh+hd=1-ip,\quad h^2=hi=ph=0.
\]
Let $(d+\Delta)^2=0$, and suppose $h\Delta$ and
$\Delta h$ are locally nilpotent with finite bounds in
each graded component.  Put
\[
 \Sigma=\Delta(1+h\Delta)^{-1}
       =(1+\Delta h)^{-1}\Delta,
\]
\[
 i'=i-h\Sigma i,\quad
 p'=p-p\Sigma h,\quad
 h'=h-h\Sigma h,\quad
 b'=b+p\Sigma i.
\]
Then these formulas define a strong deformation retract (SDR)
of cochain complexes, with the same side conditions, between
$(C_0,b')$ and $(C_1,d+\Delta)$.
\end{lem}

\begin{proof}
All inverses are finite sums.  Set
$X=1+h\Delta$ and $Y=1+\Delta h$.
The formulas give
$i'=X^{-1}i$, $p'=pY^{-1}$, and
$h'=X^{-1}h=hY^{-1}$.
Using $d\Delta+\Delta d+\Delta^2=0$ and
$dh+hd=1-ip$, direct multiplication gives
\[
 X(d+\Delta)i'=ib',
 \qquad
 p'(d+\Delta)Y=b'p.
\]
The side conditions give $p'i'=1$.
Moreover, both
\[
 X\bigl((d+\Delta)h'+h'(d+\Delta)\bigr)Y
 \quad\text{and}\quad
 X(1-i'p')Y
\]
expand to
\[
 1-ip+\Delta h+h\Delta+h\Delta^2h.
\]
Thus the chain and homotopy identities hold.
From $(d+\Delta)i'=i'b'$ and $p'i'=1$ one obtains
$b'^2=0$.  Expanding the finite formulas and using
$h^2=hi=ph=0$ gives
$h'^2=h'i'=p'h'=0$.
No vanishing of $b$, and no module-linearity of
$\Delta$, is required.
\end{proof}

\subsection{The pullback-coefficient relative complex}

Let
\[
 P_p=G_{p+1},\qquad
 \pi_p=d_{p+1}:P_p\to G_p,\qquad
 t_p=v_{p+1}:P_p\to O.
\]
Give $P$ the upper-decalage faces and degeneracies, indexed
$0,\ldots,p$ in degree $p$. The last degeneracy of $G$ is not
one of these operators. The map $\pi$ forgets the appended
terminal vertex, whereas $t$ records it.

Define
\[
 \mathscr D_\pi^{p,r,q}(R,S)=
 \Omega^r\left(
 P_p/G_p;\
 \pi_p^*\underline{\Hom}^q(v_p^*V_R,v_0^*V_S)
 \right),
 \qquad |\Phi|=p+r+q.
\]
All cochains are normalized with respect to the genuine
degeneracies of $P$.  The front and back maps retain the
appended terminal vertex.  For composable tridegrees,
define
\begin{equation}\label{eq:open-relative-product}
 \Phi\star_\pi\Psi
 =
 (-1)^{(r+q)p'+qr'}
 \Phi_{\mathrm{front}}\wedge\Psi_{\mathrm{back}}.
\end{equation}
The wedge includes ordinary matrix composition in the
displayed order; the exponent supplies its forms-first
Koszul sign.

\begin{prop}[The actual relative Hom differential]
\label{prop:open-relative-differential}
With $d_v=(-1)^pd_\pi$, the operator
\begin{equation}\label{eq:open-relative-differential}
 \mathbb D_\pi\Phi
 =
 \delta_{\mathrm{in}}\Phi+d_v\Phi
 +\pi^*r_S\star_\pi\Phi
 -(-1)^{|\Phi|}\Phi\star_\pi\pi^*r_R
\end{equation}
squares to zero and preserves normalization.
Its body-plus-relative part is
\begin{equation}\label{eq:open-column-unperturbed}
 \mathbb D_{\pi,0}
 =
 (-1)^p\bigl(d_\pi+(-1)^r d_H\bigr),
 \qquad
 d_HX=\partial_SX-(-1)^qX\partial_R.
\end{equation}
The perturbation
$\Delta_\pi=\mathbb D_\pi-\mathbb D_{\pi,0}$
preserves relative degree, raises horizontal degree,
and does not increase internal degree.
\end{prop}

\begin{proof}
The coefficient bundle is pulled from $G_p$, so its
relative derivative is canonical.  In a projective
presentation the two endpoint projectors have zero
$d_\pi$-derivative.  This assertion does not apply to
terminal projectors pulled back by $t$.

The product phase in \eqref{eq:open-relative-product}
is bilinear in the two degree triples; its two
associativity exponents therefore agree.
The usual inner-face cancellation gives the total
Leibniz rule for $\delta_{\mathrm{in}}$.
The ordinary relative Leibniz rule gives the same
total rule for $d_v$: differentiating the front form
changes the product exponent by $p'$, and
differentiating the back form changes it by $q$.
Genuine face pullbacks commute with $d_\pi$, so
$\delta_{\mathrm{in}}d_v+d_v\delta_{\mathrm{in}}=0$.

The coefficients $\pi^*r_R,\pi^*r_S$ have zero relative
derivative.  The square of
\eqref{eq:open-relative-differential} is therefore the
left and right action of their pulled-back MC
curvatures, and is zero.  At a repeated endpoint,
the unit terms complete the inner differential to
the full normalized cancellation.  In particular the
last unit coefficient is
\[
 -(-1)^{p+r+q}(-1)^{r+q}=(-1)^{p+1}.
\]
Higher structure components vanish on their degeneracies.
Thus the full operator preserves normalization.

The body signs are $(-1)^{p+r}$ at the target and
$-(-1)^{p+r+q}$ at the source.  They give
\eqref{eq:open-column-unperturbed}.  The body coefficients
are constant along each $\pi_p$-fibre, so this operator
squares to zero.  The remaining terms have degree
$(1,0,0)$ or $(k,0,1-k)$ with $k\ge1$, proving the
last assertion.
\end{proof}

\subsection{The prescribed normalized column SDR}

In root coordinates the fibre of $\pi_p$ is
\[
 D_{m_0}\times\R^p,\qquad
 g=a_{p+1},\quad z_i=h_{0,i,p+1}.
\]
Use the charts of Lemma~\ref{lem:open-based-chart} and
the fixed contraction
\begin{equation}\label{eq:open-column-radial}
 \varrho_s(g,z)
 =
 \bigl(\psi_{m_0}^{-1}(s\psi_{m_0}(g)),sz\bigr).
\end{equation}
It contracts to
\[
 \kappa_p(k)=(k,e,0,\ldots,0),
\]
which appends the root.  This is not the last-degeneracy
section.

Let $k_p$ be its relative Poincare operator.  In the
based fibre coordinates $w=(\psi_{m_0}(g),z)$, for
an $r$-form with $r\ge1$,
\[
 (k_p\omega)_w(v_1,\ldots,v_{r-1})
 =
 \int_0^1s^{r-1}
 \omega_{sw}(w,v_1,\ldots,v_{r-1})\,ds,
\]
and $k_p=0$ on relative zero-forms.

\begin{prop}[The column SDR]
\label{prop:open-column-sdr}
Set
\[
 J_{\pi,0}=\pi^*,\qquad
 e_{\pi,0}=\kappa^*,\qquad
 H_{\pi,0}=(-1)^pk_p,
\]
where $e_{\pi,0}$ kills positive relative degree.  Put
\[
 \Sigma_\pi
 =
 \Delta_\pi(1+H_{\pi,0}\Delta_\pi)^{-1},
\]
\begin{equation}\label{eq:open-column-transfer}
 \begin{split}
 J_\pi&=\pi^*,\\
 r_\pi&=e_{\pi,0}-e_{\pi,0}\Sigma_\pi H_{\pi,0},\\
 H_\pi&=H_{\pi,0}-H_{\pi,0}\Sigma_\pi H_{\pi,0}.
 \end{split}
\end{equation}
Then
\[
 r_\pi J_\pi=1,\qquad
 J_\pi r_\pi
 =
 1-\mathbb D_\pi H_\pi-H_\pi\mathbb D_\pi,
\]
with all SDR side conditions.  The small differential
is literally the original raw Hom differential.
The construction applies to the original, possibly
noninvertible edge maps.
\end{prop}

\begin{proof}
Cartan's homotopy identity gives
\[
 d_\pi k_p+k_pd_\pi=1-\pi_p^*\kappa_p^*.
\]
The operator kills pullbacks and satisfies
$\kappa_p^*k_p=0$.  Its square is zero: the iterated
integral inserts two collinear radial vectors into an
alternating form.

A genuine $P$-degeneracy repeats a prefix vertex.
It preserves $m_0$ and $g$ and either inserts an
initial zero in the $z$-list or repeats one of its
entries.  It therefore commutes with
\eqref{eq:open-column-radial}.  Thus all column maps
preserve normalization.  They act as scalars on the
pulled-back coefficient bundle and commute with the
plain $d_H$.  The body cross terms in the homotopy
identity are
\[
 (-1)^{r-1}d_Hk_p+(-1)^rk_pd_H=0.
\]
Consequently these maps give an SDR for
$\mathbb D_{\pi,0}$ over the raw body-only differential
$(-1)^pd_H$.

The perturbation preserves relative degree and the
homotopy lowers it.  On a relative $r$-form,
\[
 (H_{\pi,0}\Delta_\pi)^{r+1}
 =
 (\Delta_\pi H_{\pi,0})^{r+1}=0.
\]
Lemma~\ref{lem:open-finite-perturbation} applies.
Moreover $H_{\pi,0}\Delta_\pi J_{\pi,0}=0$, since its
argument has relative degree zero.  Thus the perturbed
inclusion is still $\pi^*$.  The simplicial map $\pi$
preserves the original products and twisting data, so
\[
 \Delta_\pi J_{\pi,0}
 =
 J_{\pi,0}\bigl(D_{\mathrm{raw}}-(-1)^pd_H\bigr).
\]
The transferred differential is exactly
$D_{\mathrm{raw}}$.

The formulas are finite in every relative degree.
Smoothness follows from compact integration in the
jointly smooth charts.  They act intrinsically on
the pullback bundles, so they glue on finite-projective
bodies.
\end{proof}

\subsection{Full terminal conjugation}

Assume now that the edges of $R,S$ are actual smooth
chain isomorphisms.  For a raw object $R$, let $B_R$
be the terminal-body representation on $P$: its fibre
is $V_R|_{m_*}$, its body is $\partial_R(m_*)$, its
edges are identities, and its higher structures vanish.

\begin{lem}[The upper-decalage gauge]
\label{lem:open-upper-gauge}
The maps
\[
 U_p^R=(-1)^pr_{R,p+1}:B_R\longrightarrow\pi^*R
\]
form a normalized closed degree-zero morphism.
Their closed inverse is given by a finite geometric series.
\end{lem}

\begin{proof}
The component has type $V_{m_{p+1}}\to V_{m_0}$ and
internal degree $-p$.  For $p\ge1$, normalization
follows from that of $r_{p+1}$.
Expanding its Hom differential gives
\[
 \begin{split}
 (DU)_p={}&
 \sum_{i=1}^{p-1}(-1)^{p+i-1}r_p(d_i\sigma)\\
 &+\sum_{a=0}^{p}(-1)^{a(p-a)}
 r_a(\sigma[0,\ldots,a])
 r_{p-a+1}(\sigma[a,\ldots,p+1])\\
 &-r_p(d_p\sigma)
 +(-1)^{p+1}r_{p+1}(\sigma)\partial_R(m_{p+1}).
 \end{split}
\]
The source identity edge supplies the term
$-r_p(d_p\sigma)$.  This is
$(-1)^{p+1}$ times the original MC equation in degree
$p+1$, and vanishes.  At $p=0$ the assertion is the
edge chain-map equation.

Let $X=U_0$ and $v=X^{-1}\star U-1$.
If the body width is $w$, then $v$ has only positive
horizontal components and $v^{\star(w+1)}=0$.
Therefore
\[
 U^{-1}
 =
 \left(\sum_{j=0}^{w}(-v)^{\star j}\right)\star X^{-1}.
\]
The inverse is two-sided in the graded cochain algebra,
has the reverse endpoint type, and is closed by the
Hom Leibniz rule.  The zero-cochain $X$ alone need not
be closed for the full raw differential.
\end{proof}

Define
\[
 \mathfrak A(\chi)
 =
 (U^S)^{-1}\star_\pi\chi\star_\pi U^R,
 \qquad
 \mathfrak A^{-1}(\zeta)
 =
 U^S\star_\pi\zeta\star_\pi(U^R)^{-1}.
\]
The terminal complex has coefficients
$t_p^*\underline{\Hom}(V_R,V_S)$.  Define its differential
by full conjugation:
\begin{equation}\label{eq:open-terminal-differential}
 \mathbb D_t
 =
 \mathfrak A\mathbb D_\pi\mathfrak A^{-1}
 =
 \mathbb D_{t,0}+\nabla_U,
\end{equation}
where
\[
 \mathbb D_{t,0}
 =
 \delta_{\mathrm{full}}+(-1)^{p+r}d_{H,t},
 \qquad
 \nabla_U=\mathfrak A d_v\mathfrak A^{-1}.
\]
Closedness of the full $U$ proves this split.  It
retains all body and higher-curvature terms.
The derivative is applied only after returning by
$\mathfrak A^{-1}$ to the bundle with its canonical
relative derivative.

Each component of $\mathfrak A$ or $\mathfrak A^{-1}$
adds horizontal degree $a\ge0$ and lowers internal
degree by $a$.  Thus $\nabla_U$ raises relative degree
by one and does not increase internal degree.
Conjugation also gives
\[
 \mathbb D_{t,0}^2=\nabla_U^2=0,\qquad
 \mathbb D_{t,0}\nabla_U+\nabla_U\mathbb D_{t,0}=0.
\]
These identities involve the full operators: the higher
structure components participate in the curvature cancellations.

\subsection{The actual horizontal SDR}

Let $A^p=C^\infty(P_p)$, with $A$ acting through the
terminal point $m_*$.  The relative frame is
\[
 \alpha=d_\pi g\,g^{-1},\qquad dz_1,\ldots,dz_p.
\]
For the root-changing face,
\begin{equation}\label{eq:open-relative-frame}
 d_0^*\alpha=\alpha-Z_{m_*}dz_1,\qquad
 d_0^*dz_j=dz_{j+1}-dz_1.
\end{equation}
Indeed,
\[
 g'=g\exp(-z_1Z_{m_0})a_1^{-1},\qquad
 z'_j=H_{1,j+1}+z_{j+1}-z_1,
\]
and the prefix quantities have zero relative derivative.
Other faces delete the appropriate $z$, while the
initial degeneracy inserts zero and the remaining
degeneracies repeat an entry.

Let $\mathcal V^p$ be free over $A$ on
$\mathbf e_1,\mathbf e_2,\mathbf v_1,\ldots,\mathbf v_p$,
with the operators prescribed by
\eqref{eq:open-relative-frame}.  Put
\[
 \mathcal H^q=\Gamma(O,\underline{\Hom}^q(V_R,V_S)),
 \qquad
 W^{p,r,q}=\Lambda_A^r\mathcal V^p\otimes_A\mathcal H^q.
\]
The actual terminal relative forms are
$A^p\otimes_AW^{p,r,q}$.  This is a finite-frame and
finite-projective bundle identity; arbitrary jointly
smooth dependence remains in $A^p$.

Let $j$ pull coefficient functions back from $m_*$,
and let $\epsilon$ evaluate coefficient functions on
the constant simplex, retaining the formal $W$-tensor.
These are cosimplicial $A$-linear maps with
$\epsilon j=1$.

For $p\ge1$, define
\[
 f_j^p:[p+1]\longrightarrow[p],\qquad
 f_j^p(s)=
 \begin{cases}
 s,&s\le j,\\
 p,&s>j,
 \end{cases}
 \quad 0\le j\le p-1,
\]
and $h_j^p=G(f_j^p):P_{p-1}\to P_p$.  Set
\begin{equation}\label{eq:open-horizontal-L}
 L_p=
 \sum_{j=0}^{p-1}(-1)^{j+1}
 (h_j^p)^*\otimes\sigma_W^j,\qquad L_0=0.
\end{equation}
The first factor acts on coefficient functions, and
the second is the formal-frame codegeneracy.

\begin{lem}[Horizontal coefficient contraction]
\label{lem:open-horizontal-contraction}
The maps above satisfy
\[
 \delta_{\mathrm{full}}L+L\delta_{\mathrm{full}}
 =1-j\epsilon.
\]
After genuine cosimplicial normalization, there are
maps $i_0,p_0,h_0$ with
\[
 p_0i_0=1,\quad
 \delta_{\mathrm{full}}h_0+h_0\delta_{\mathrm{full}}
 =1-i_0p_0,\quad
 h_0^2=h_0i_0=p_0h_0=0.
\]
They preserve relative and internal degree, and $h_0$
lowers horizontal degree by one.
\end{lem}

\begin{proof}
For the first identity, write $\delta=\delta_{\mathrm{full}}$
within this proof and expand $\delta L+L\delta$.
If $\mathbf A_{ij}$ denotes an unsigned term from $\delta L$
and $\mathbf B_{ij}$ one from $L\delta$, the order-map identities give
\[
 \mathbf B_{ij}=
 \begin{cases}
 \mathbf A_{i,j-1},&i<j,\\
 \mathbf A_{i-1,j},&i>j+1.
 \end{cases}
\]
Their signs are opposite.  These pairs exhaust the
off-diagonal terms.  For $i=j,j+1$, the $W$-map is
the identity, and the coefficient terms telescope:
successively more vertices are retained instead of
being sent to the terminal vertex.  The sum is the
identity minus constant-terminal evaluation.  This
also proves the identity in horizontal degree zero.

Use
\[
 N^pT=\bigcap_{i=0}^{p-1}\ker s^i,\qquad
 \Pi_N^p=(1-d^ps^{p-1})\cdots(1-d^1s^0).
\]
The cosimplicial identities split the full complex
into its normalized part and the positive-coface
degenerate subcomplex.  Thus $\Pi_N$ and the
normalized inclusion are chain maps.  Put
$L_N=\Pi_NL\iota_N$.

The codegeneracy kernels on $\mathcal V^p$ are generated
by
\[
 \mathbf v_1,\quad
 \mathbf v_2-\mathbf v_1,\quad
 \ldots,\quad
 \mathbf v_p-\mathbf v_{p-1}.
\]
Together with $\mathbf e_1,\mathbf e_2$, these form a
free basis.  On the exterior algebra the corresponding
kernels are the ideals generated by these vectors.
Their intersection consists exactly of basis monomials
containing all $p$ vectors.  Consequently
\begin{equation}\label{eq:open-exterior-normalization}
 N^p\Lambda_A^r\mathcal V^\bullet
 =
 \Lambda_A^{r-p}(A\mathbf e_1\oplus A\mathbf e_2)
 \wedge\mathbf v_1\wedge\cdots\wedge\mathbf v_p.
\end{equation}
It vanishes unless $0\le r-p\le2$.  Tensoring this
split calculation with $\mathcal H^q$ is legitimate.

Identify the normalized small module with
\[
 \mathbb B_{R,S}
 =
 B_O\otimes_A\mathcal H,
\]
by fixing
\begin{equation}\label{eq:open-iota-phase}
 i_0(f\eta^p\theta_I X)
 =
 c_{p,j}f(m_*)\alpha_I\wedge dz_1\wedge\cdots\wedge dz_p\,X,
 \quad
 c_{p,j}=(-1)^{pj+p(p+1)/2},\quad j=|I|.
\end{equation}
The map $p_0$ is coefficient evaluation followed by
the inverse of this identification.  Set
\[
 P_0=1-i_0p_0,\qquad a=P_0L_NP_0,\qquad
 \boxed{h_0=a\delta_{\mathrm{full}}a}.
\]
Then
\[
 \delta_{\mathrm{full}}a+a\delta_{\mathrm{full}}=P_0,
 \qquad ai_0=p_0a=0,\qquad P_0a=aP_0=a.
\]
It follows that $(\delta_{\mathrm{full}}a)^2
=\delta_{\mathrm{full}}a$,
$(a\delta_{\mathrm{full}})^2=a\delta_{\mathrm{full}}$,
and
$\delta_{\mathrm{full}}a^2\delta_{\mathrm{full}}=0$.
These identities give the homotopy and
all three side conditions for $h_0$.
Every factor preserves $r,q$, so $h_0$ lowers $p$
by one.
\end{proof}

The homotopy is normalized using only the horizontal
differential.  It acts as the identity on the body-Hom
tensor factor and is $A$-linear through $m_*$.  Thus it
commutes with the plain body-Hom differential, and the
two total body signs cancel because it lowers $p$.
The same maps therefore give an SDR for
$\mathbb D_{t,0}$, with small differential
\begin{equation}\label{eq:open-small-unperturbed}
 D_{\mathbb B,0}(\omega X)
 =
 Q_h\omega\,X+(-1)^{|\omega|}\omega\,d_HX,
 \qquad Q_h=\eta\iota_{(-y,x)}.
\end{equation}
The small complex retains the full body and Koszul differentials.

Apply Lemma~\ref{lem:open-finite-perturbation} to
$\nabla_U$:
\begin{equation}\label{eq:open-horizontal-transfer}
 \begin{split}
 \Sigma_U&=\nabla_U(1+h_0\nabla_U)^{-1},\\
 D_{\mathrm{tr}}&=D_{\mathbb B,0}+p_0\Sigma_Ui_0,\\
 i_h&=i_0-h_0\Sigma_Ui_0,\qquad
 p_h=p_0-p_0\Sigma_Uh_0,\\
 H_h&=h_0-h_0\Sigma_Uh_0.
 \end{split}
\end{equation}
For an input of total degree $m$ and relative degree $r$,
both $h_0\nabla_U$ and $\nabla_Uh_0$ have total degree
zero and raise relative degree by one.  A nonzero
output after $k$ iterations satisfies
$r+k\le m-\ell$.  Hence
\begin{equation}\label{eq:open-horizontal-finite-bound}
 (h_0\nabla_U)^k=(\nabla_Uh_0)^k=0
 \quad\text{when }k>m-\ell-r.
\end{equation}
Thus all displayed maps are finite in each total degree,
and
\[
 p_hi_h=1,\qquad
 i_hp_h=1-\mathbb D_tH_h-H_h\mathbb D_t,
 \qquad H_h^2=H_hi_h=p_hH_h=0.
\]

\begin{prop}[The fixed comparison maps]
\label{prop:open-fixed-comparison-maps}
Define
\[
 V_{\mathrm{mod}}=r_\pi\mathfrak A^{-1}i_h,\qquad
 S_{\mathrm{mod}}=p_h\mathfrak A J_\pi,\qquad
 T=F V_{\mathrm{mod}}.
\]
Then $V_{\mathrm{mod}}$ and $S_{\mathrm{mod}}$ are
chain-homotopy inverse maps between the actual raw Hom
complex and $(\mathbb B_{R,S},D_{\mathrm{tr}})$.
The map $T$ is a chain map to the actual cohesive Hom
complex, and
\begin{equation}\label{eq:open-fixed-TS-homotopy}
 TS_{\mathrm{mod}}-F
 =
 -d_{\Coh}\mathcal K_F-\mathcal K_F D_{\mathrm{raw}},
 \qquad
 \mathcal K_F=
 Fr_\pi\mathfrak A^{-1}H_h\mathfrak A J_\pi.
\end{equation}
Every displayed comparison map is internally
nonincreasing.
\end{prop}

\begin{proof}
The two model homotopies are
\[
 K_C=r_\pi\mathfrak A^{-1}H_h\mathfrak A J_\pi,
 \qquad
 K_{\mathbb B}
 =p_h\mathfrak A H_\pi\mathfrak A^{-1}i_h.
\]
The two SDR identities give
\[
 V_{\mathrm{mod}}S_{\mathrm{mod}}
 =1-D_{\mathrm{raw}}K_C-K_CD_{\mathrm{raw}},
\]
\[
 S_{\mathrm{mod}}V_{\mathrm{mod}}
 =1-D_{\mathrm{tr}}K_{\mathbb B}
   -K_{\mathbb B}D_{\mathrm{tr}}.
\]
Applying $F$ to the first identity proves
\eqref{eq:open-fixed-TS-homotopy}, including its minus
signs.

The elementary column maps preserve internal degree;
$\Delta_\pi$ does not increase it.  Positive components
of $U,U^{-1}$ lower it.  The relative derivative
preserves it, and the horizontal maps preserve it.
Therefore all finite transfer expressions are internally
nonincreasing.  The map $F$ preserves the internal matrix
block by its entrywise quotient formula.
Thus internal degree orders the comparison operators, while
the body differential itself can raise internal degree.
\end{proof}

\subsection{The scalar inclusion and its actual differential}

For scalar unit coefficients, $U=1$ and
$\nabla_U=d_v$.  Denote \eqref{eq:open-iota-phase} by
$\iota$.  Split $Q=Q_h+Q_v$, with
\[
 Q_vf=f_x(x\theta_1+y\theta_2),\quad
 Q_v\theta_1=0,\quad
 Q_v\theta_2=\theta_1\theta_2,\quad Q_v\eta=0.
\]
The exact frame formulas give
\[
 \delta_{\mathrm{full}}
 (\alpha_I\wedge dz_1\wedge\cdots\wedge dz_p)
 =
 (-1)^j(\iota_Z\alpha_I)
 \wedge dz_1\wedge\cdots\wedge dz_{p+1}.
\]
The terms without a root-face change cancel by expanding
$\bigwedge_{i=2}^{p+1}(dz_i-dz_1)$.
Since $c_{p+1,j-1}=(-1)^jc_{p,j}$,
\[
 \delta_{\mathrm{full}}\iota=\iota Q_h.
\]
Also
\[
 d_\pi f(m_*)
 =f_x(m_*)(x_*\alpha_1+y_*\alpha_2),\quad
 d_\pi\alpha_1=0,\quad
 d_\pi\alpha_2=\alpha_1\wedge\alpha_2,\quad d_\pi dz_i=0.
\]
The identity $c_{p,j+1}=(-1)^pc_{p,j}$ gives
\begin{equation}\label{eq:open-scalar-iota-chain}
 d_v\iota=\iota Q_v.
\end{equation}
Since $h_0\iota=0$, also $h_0d_v\iota=0$.  Thus
\begin{equation}\label{eq:open-scalar-transfer-differential}
 i_h^{\mathrm{sc}}=\iota,\qquad
 D_{\mathrm{tr}}^{\mathrm{sc}}=Q.
\end{equation}
It remains to compute the composite
$\operatorname{ve}r_\pi^{\mathrm{sc}}\iota$.

\subsection{Weighted smooth functions at the unit section}

Use root coordinates
\[
 u_i=(\log\mathfrak a_i,b_i),\quad u_0=0,\qquad
 H_{ij}=h_{0ij},\quad H_{0j}=H_{ii}=0.
\]
Set
\begin{equation}\label{eq:open-adjacent-coordinates}
 e_i=u_i-u_{i-1},\qquad
 c_{ij}=H_{ij}-H_{i-1,j}-H_{i,j-1}+H_{i-1,j-1}.
\end{equation}
The inverse is
\[
 u_i=\sum_{a\le i}e_a,\qquad
 H_{ij}=\sum_{a\le i<b\le j}c_{ab}.
\]
The degeneracy repeating vertices $k-1,k$ is exactly
\begin{equation}\label{eq:open-degeneracy-planes}
 S_k=\{e_k=0,\ c_{ik}=0\ (i<k),\
                    c_{kj}=0\ (k<j)\}.
\end{equation}
Indeed the equations recover repeated cumulative
$u$-entries, repeated rows and columns of $H$, and
$H_{k-1,k}=0$; deleting that slot gives the lower simplex.

Assign each component of $e_i$ weight one and label
$i$, and assign $c_{ij}$ weight two and labels $i,j$.
Let $\mathcal W^{\ge d}$ be the smooth unit-normal
germs whose Taylor coefficients of weighted degree
less than $d$ vanish.  Base functions have weight zero.

\begin{lem}[The weighted function bound]
\label{lem:open-weighted-functions}
For every $n\ge1$,
\[
 C_N^n(G)\cap\mathcal W^{\ge n+1}\subseteq J_n
\]
at unit-germ level and, after germ localization, globally.
In particular,
\[
 \operatorname{ve}_n
 \bigl(C_N^n(G)\cap\mathcal W^{\ge n+1}\bigr)=0.
\]
Every normalized degree-$n$ function has
normal weight at least $n$.
\end{lem}

\begin{proof}
Choose a small normal-coordinate polydisc, with base
parameters, inside $G_n$, stable under setting any
subset of the coordinates in
\eqref{eq:open-adjacent-coordinates} to zero.
Let $P_k$ zero all coordinates carrying label $k$.
The maps commute and retract onto
\eqref{eq:open-degeneracy-planes}.  Each $P_j^*$
sends a generator of $I(S_k)$ to itself or to zero,
so preserves $I(S_k)^2$.  The operator
\[
 \mathcal N_n=\prod_{k=1}^n(1-P_k^*)
\]
normalizes functions and fixes normalized functions.

For a normalized $f$ of weight at least $n+1$,
ordinary Taylor expansion gives
\[
 f(\xi)=T_nf(\xi)+
 \sum_{|\gamma|=n+1}\xi^\gamma R_\gamma(\xi),
\]
\[
 R_\gamma(\xi)=
 \frac{n+1}{\gamma!}
 \int_0^1(1-s)^n\partial^\gamma f(s\xi)\,ds.
\]
The remainders are smooth with all base and passive
parameter derivatives.  Restriction to a coordinate
plane commutes with Taylor truncation.  Hence $T_nf$
is normalized, and every surviving coordinate monomial
covers all $n$ labels.  Its weight is its number of
label occurrences.  Weight greater than $n$ forces
one label to occur twice, placing the monomial in a
selected $I(S_k)^2$.

Each remainder monomial has $n+1$ coordinate factors,
so it also repeats a label.  Its smooth multiple
need not be normalized separately.  Applying
$\mathcal N_n$ gives
\[
 \mathcal N_n(\xi^\gamma R_\gamma)
 \in C_N^n\cap I(S_{k(\gamma)})^2,
\]
because every factor of $\mathcal N_n$ preserves the
selected squared ideal.  Applying the normalizer to
the Taylor formula therefore gives exactly a sum of
separately normalized squared-ideal summands.

Patch these local expressions near the unit section
using a locally finite family of cutoffs supported
inside their coordinate domains.  The coordinate
generators are finite in each arity.  Grouping the
locally finite coefficients by their finitely many
repeated-label generator pairs gives global smooth
squared-ideal summands, still normalized.  The remaining
normalized difference vanishes near the whole unit
section and belongs to $J_n$ by
\cite[Lemma~3.1.5]{CabreraDelHoyo2026}.

Finally, a coordinate monomial of weight less than
$n$ cannot cover every label.  Applying the same
restriction argument to the low Taylor polynomial
shows that such terms do not occur in a normalized
function.
\end{proof}

We apply the lemma entrywise after the column retraction,
when the outputs are normalized functions.  The full
differentiating ideal remains
$J+\delta_{\mathrm{full}}J$, with its
allowed-transposition relations.

\subsection{Filteredness of the prescribed operators}

Give relative coordinate differentials the weights of
their coordinates.  Equivalently, $\alpha_1,\alpha_2$
have weight one and $dz_i$ weight two.  The change
between these frames is filtered:
\[
 \alpha_1=d_\pi u^1,\qquad
 \alpha_2=d_\pi u^2-b\,d_\pi u^1.
\]
The term $b\,d_\pi u^1$ has weight two.

\begin{lem}[Filtered operators]
\label{lem:open-filtered-operators}
All maps used in
\eqref{eq:open-column-transfer} and
\eqref{eq:open-horizontal-transfer}, as well as
$\mathfrak A,\mathfrak A^{-1}$, preserve normal weight.
The assertion uses the fixed charts $\psi_m$ and the
fixed horizontal homotopy $h_0$.
\end{lem}

\begin{proof}
Group multiplication and inversion have additive
linear part in $(\log\mathfrak a,b)$ and higher terms
of weight at least two.  The $H$-coordinates transform
linearly through triangle identities.  Root-changing
group formulas contain $\exp(-HZ)$, whose nonconstant
terms have weight at least two.  Base movement satisfies
$f(a_im_0)=f(m_0)+\mathcal W^{\ge1}$.
Thus every genuine face and degeneracy is filtered.
Equation~\eqref{eq:open-relative-frame} is compatible:
$Z_{m_*}dz_1$ has weight two, not one.

Write $\chi_m(u)=\psi_m(g(u))$.  The actual radial
group-coordinate map is
$\chi_m^{-1}(s\chi_m(u))$.  Its value at $u=0$ is zero
and its derivative there is $s\,1$.  Smooth Taylor
expansion, uniformly for $s\in[0,1]$ locally in $m$,
gives
\[
 \chi_m^{-1}(s\chi_m(u))
 =su+\mathcal W^{\ge2}.
\]
The $ds$-coefficient and relative-differential
coefficients of the pulled-back group coordinates
have weight at least one.  For a central coordinate,
\[
 d(sz_i)=s\,dz_i+z_i\,ds
\]
has weight two in both terms.  Contraction with
$\partial_s$ and compact integration therefore
preserve weight.  This proves filteredness of the
$k_p$ determined by the fixed charts.

The root-appending evaluation is also filtered.
It retains prefix data: for example,
$e_{p+1}=-u_p$ there, and
$c_{i,p+1}=-H_{i,p}+H_{i-1,p}$.
These terms have their original weights.

For \eqref{eq:open-horizontal-L}, the order map sends
a rooted group coordinate to a rooted group coordinate
or zero, and an $H$-coordinate to an $H$-coordinate
or zero.  On the formal frame, codegeneracies preserve
the weight-one $\mathbf e$-vectors and the weight-two
$\mathbf v$-vectors or send them to zero.  Thus $L$
is filtered.  Genuine normalization gives the same
property for $L_N$.

The map $p_0$ evaluates coefficient functions at the
constant simplex while retaining their formal frame
tensor.  Positive-weight coefficients vanish; surviving
frame weights remain unchanged.  The filtered frame
isomorphism ensures that this conclusion does not
depend on a chosen expression of the form.
The map $i_0$ reevaluates coefficients at $m_*$ and
multiplies by the prescribed frame monomial of exactly
the CE weight.  Hence $i_0,p_0,P_0,a,h_0$ are filtered.

Normalization and the preceding function lemma give
$r_k\in\mathcal W^{\ge k}$ for $k\ge2$, and strict
units give $r_1-1\in\mathcal W^{\ge1}$.
The same estimates apply to $U_0-1$ and its inverse,
using the exact identity
$U_0^{-1}-1=-U_0^{-1}(U_0-1)$.
Positive components of $U^{-1}$ are finite products.
The relative derivative is filtered because a
coordinate and its differential have the same weight.
It follows that $\Delta_\pi$, both conjugations,
and $\nabla_U$ are filtered.  Every transfer expression
is finite in its established relative-degree bound,
so every resulting operator is filtered.
\end{proof}

For clarity, the associated symbol complexes consist
of fixed-weight polynomial normal germs, with smooth
coefficients in $A$.  The relative symbol uses the
genuine $P$-normalization.  Finite normalization supplies
polynomial representatives.  This is a weighted
associated graded, not an ordinary tangent
linearization or a completed power series construction.
The scalar CE component $B_O^n$ is concentrated in
weight $n$, so its weight-preserving differential is
zero.  The raw and relative symbol differentials
generally remain nonzero.

\subsection{The prescribed scalar retraction is the identity}

\begin{thm}[Scalar normal-symbol calibration]
\label{thm:open-scalar-transfer-identity}
The actual map
\[
 \tau_O=\operatorname{ve}r_\pi^{\mathrm{sc}}\iota:
 (B_O,Q)\longrightarrow(B_O,Q)
\]
is literally the identity in every degree.
\end{thm}

\begin{proof}
The column SDR identities of
Proposition~\ref{prop:open-column-sdr} descend to the
weighted polynomial symbol:
\[
 r_{\pi,\mathrm{lin}}\pi_{\mathrm{lin}}^*=1,\qquad
 \pi_{\mathrm{lin}}^*r_{\pi,\mathrm{lin}}
 =
 1-D_{\mathrm{lin}}H_{\pi,\mathrm{lin}}
  -H_{\pi,\mathrm{lin}}D_{\mathrm{lin}}.
\]
The map $\pi_{\mathrm{lin}}^*$ is a DGA map.

We check that the leading inclusion is a DGA map
with the phase in \eqref{eq:open-iota-phase}.
Take $f\eta^p\theta_I$ and $g\eta^{p'}\theta_J$,
with $j=|I|$, $j'=|J|$.
For $p>0$, the back central differentials pull to
$dz_{p+a}-dz_p$.  Their terms containing $dz_p$
vanish against the front central block.
The leading back $\alpha$-symbol is the same
$\alpha$; its $Z\,dz_p$ correction has higher weight.
The relative product contributes
$(-1)^{(p+j)p'}$, and moving the back $\alpha_J$
past the front $p$ central differentials contributes
$(-1)^{pj'}$.  The exact identity
\[
 c_{p,j}c_{p',j'}(-1)^{(p+j)p'+pj'}
 =c_{p+p',j+j'}
\]
follows from
\[
 \frac{p(p+1)}2+\frac{p'(p'+1)}2+pp'
 =\frac{(p+p')(p+p'+1)}2.
\]
The remaining signs are exactly those of exterior
multiplication.  Overlapping odd generators give zero.
If $p=0$, the back map has no preceding root shift
and no $dz_0$ occurs.  If $p'=0$, the back central
block is empty.  If both are zero, this is ordinary
exterior multiplication.  Thus $\iota_{\mathrm{lin}}$
is a unital DGA map in every degree.

The induced $\operatorname{ve}_{\mathrm{lin}}$ is
the DGA quotient furnished by the actual
differentiation, including $\delta J$ and its
transposition relations.  On cohomology,
$H(r_{\pi,\mathrm{lin}})
=H(\pi_{\mathrm{lin}}^*)^{-1}$ is an algebra map.
Both scalar CE symbol endpoints have zero differential.
Consequently their actual composite
\[
 \tau_{\mathrm{lin}}
 =
 \operatorname{ve}_{\mathrm{lin}}
 r_{\pi,\mathrm{lin}}\iota_{\mathrm{lin}}
\]
is an algebra map.  Multiplicativity here follows from the
induced cohomology map; the column retraction $r_\pi$ is used
as a chain map.

In the following calculation, $\delta_W$ is the
weighted associated-graded map induced by the full
horizontal differential on normalized scalar relative
forms on $P_p/G_p$:
\[
 \delta_{\mathrm{full}}
   =\sum_{i=0}^{p+1}(-1)^i d_i^*.
\]
The original $G$-face $d_{p+2}$ is omitted.
It sends bidegree $(p,r)$ to $(p+1,r)$, with no
additional relative-form sign.  Unlike the scalar CE
symbol differential, this horizontal differential
need not vanish.
The index $a\in\{1,2\}$ in $u_a$ is a group-coordinate
index, and $u_{1,a}$ is that coordinate of the first
prefix root arrow.

It remains to calibrate generators.  Functions are
sent to themselves by root-appending evaluation.
For a group coordinate, at horizontal degree zero,
\[
 H_{\pi,0}(du_a)=u_a,\qquad
 \delta_Wu_a=-u_{1,a}
\]
in the symbol model.  The next homotopy vanishes
on a relative zero-form, so
\[
 r_{\pi,\mathrm{lin}}(du_a)=u_{1,a}.
\]
Equation~\eqref{eq:open-ve-generators} gives
$\tau_{\mathrm{lin}}(\theta_i)=\theta_i$.

For $\eta$, the inclusion is $-dz_1$.  Since
$H_{\pi,0}=-k_1$ in horizontal degree one,
\[
 H_{\pi,0}(-dz_1)=z_1.
\]
On $P_2=G_3$, put
$H_{12}=h_{012}$, $z_1=h_{013}$, $z_2=h_{023}$.
The tetrahedral equation gives
$h_{123}=H_{12}+z_2-z_1$, hence
\[
 \delta_Wz_1=h_{123}-h_{023}+h_{013}=H_{12}.
\]
The column projection correction has a minus sign,
so $r_{\pi,\mathrm{lin}}(-dz_1)=-H_{12}$.
Its differentiated image is $+\eta$.
Thus $\tau_{\mathrm{lin}}$ fixes
$A,\theta_1,\theta_2,\eta$, and is the identity on
the entire polynomial CE algebra, including every
power of $\eta$.

Degree zero was already settled by the identity on $A$.
For $n\ge1$ and $x\in B_O^n$, the actual
$r_\pi^{\mathrm{sc}}\iota(x)$ is a normalized raw
degree-$n$ cochain of weight at least $n$.
Its weight-$n$ raw symbol maps to $x$ under
$\operatorname{ve}_{\mathrm{lin}}$.
Lemma~\ref{lem:open-weighted-functions} annihilates
the smooth remainder of weight $n+1$.
Therefore $\operatorname{ve}r_\pi^{\mathrm{sc}}\iota(x)=x$.
Only a finite normal symbol is used.
\end{proof}

\subsection{The genuine coefficient diagonal}

\begin{thm}[Identity internal diagonal]
\label{thm:open-transfer-identity-diagonal}
For every invertible-edge pair $R,S$, the fixed chain
map
\[
 T=Fr_\pi\mathfrak A^{-1}i_h
\]
has identity internal-degree-preserving diagonal under
the graded identification
\[
 \mathbb B_{R,S}=B_O\otimes_A\mathcal H.
\]
\end{thm}

\begin{proof}
Every factor is internally nonincreasing by
Proposition~\ref{prop:open-fixed-comparison-maps}.
Thus the diagonal of a composite is the composite
of its diagonals.  This is a statement about
operators, not about an associated-graded differential.

Work in local homogeneous frames and take
$\omega X$, where $\omega\in B_O^n$ and $X$ is an
arbitrary constant homogeneous rectangular matrix of
internal degree $q$.  All smooth coefficients are
included in $\omega$.

On $\sigma\in P_p$, the two diagonal endpoint
transports are
\[
 L_S=r_{S,1}(\sigma[0,p+1]),\qquad
 L_R=r_{R,1}(\sigma[p,p+1]).
\]
They are distinct pullbacks.  Diagonal conjugation is
$\mathfrak A_0^{-1}(\chi)=L_S\chi L_R^{-1}$.
Put
\[
 \kappa_S=L_S^{-1}d_\pi L_S,\qquad
 \kappa_R=L_R^{-1}d_\pi L_R.
\]
For tridegree $(p,r,q)$, differentiating this expression
gives
\begin{equation}\label{eq:open-diagonal-connection}
 \nabla_{U,0}\chi
 =
 d_v\chi+\kappa_S\star_\pi\chi
 -(-1)^{p+r+q}\chi\star_\pi\kappa_R.
\end{equation}
The left product has phase $(-1)^p$ and the right
product phase $(-1)^q$.  Thus in bare matrix-valued
wedge notation the same formula is
\[
 (-1)^p\left(
 d_\pi\chi+\kappa_S\wedge_{\mathrm{bare}}\chi
 -(-1)^r\chi\wedge_{\mathrm{bare}}\kappa_R
 \right).
\]
This verifies both endpoints for every integer $q$.
Strict units give $L_R,L_S=1+\mathcal W^{\ge1}$,
so $\kappa_R,\kappa_S$ have weight at least one.
Hence
$\mathcal C=\nabla_{U,0}-d_v$ raises weight.

Using
$i_h=(1+h_0\nabla_U)^{-1}i_0$, consider its diagonal
finite expansion.  An all-$d_v$ positive correction
vanishes at its rightmost factor because
\[
 h_0d_vi_0(\omega X)
 =h_0i_0(Q_v\omega\,X)=0.
\]
Every remaining term contains a $\mathcal C$, so
\[
 \operatorname{diag}_q i_h(\omega X)
 =\iota(\omega)X+E_1,\qquad
 E_1\in\mathcal W^{\ge n+1}.
\]
The body differential remains in
\eqref{eq:open-small-unperturbed}.  It is absent from
these perturbation words because $h_0$ was normalized
using only $\delta_{\mathrm{full}}$.

Next, the exact identity
\[
 L_S\chi L_R^{-1}-\chi
 =(L_S-1)\chi L_R^{-1}
   +\chi(L_R^{-1}-1)
\]
shows that diagonal conjugation differs from the
scalar identity by one additional normal weight.

Finally, write $r_1=1+\rho$ in these local frames.
The diagonal of the column perturbation is
\[
 \Delta_{\pi,0}=\delta_{\mathrm{full}}+\mathcal E,
\]
\[
 \mathcal E\chi
 =
 \rho_S^\pi\star_\pi\chi
 -(-1)^{p+r+q}\chi\star_\pi\rho_R^\pi.
\]
The right unit has coefficient
$-(-1)^{p+r+q}(-1)^{r+q}=(-1)^{p+1}$, exactly
the scalar last-face sign.  Since $\rho$ has weight
at least one, $\mathcal E$ raises weight.
All $r_k$, $k\ge2$, lower internal degree.

Compare the finite column formulas with
$\Delta_{\pi,0}=\delta_{\mathrm{full}}+\mathcal E$
and with $\delta_{\mathrm{full}}$ alone.
Every difference word contains an $\mathcal E$;
all other factors preserve weight.  Consequently
\[
 \operatorname{diag}_q
 (r_\pi\mathfrak A^{-1}i_h)(\omega X)
 =
 r_\pi^{\mathrm{sc}}\iota(\omega)X+E,
 \qquad E\in\mathcal W^{\ge n+1}.
\]
This is an entrywise comparison of Hom operators.

The final difference is a normalized function-valued
raw degree-$n$ cochain: both complete outputs are
normalized.  Individual intermediate errors need
not be normalized, nor preserve every squared ideal.
For $n\ge1$, Lemma~\ref{lem:open-weighted-functions}
annihilates $E$. For $n=0$, the whole space $G_0=O$
is the unit section, so the weight-one function $E$
is identically zero. In either case,
Theorem~\ref{thm:open-scalar-transfer-identity} gives
\[
 \operatorname{diag}_qT(\omega X)
 =\tau_O(\omega)X=\omega X.
\]

To globalize, fix a point of the unit section.
There are finitely many operations in each total
degree.  All their geometric maps take units to units,
uniformly over the compact radial parameter.
Shrink the normal neighborhood so that every
intermediate vertex and radial image lies in the
chosen source and target frame patches.
The radial chart is the restriction of the fixed
$\psi_O$, not a new chart for the patch.
Testing all smooth scalar coefficients times the
constant homogeneous matrix units tests all local
sections.  Such sections can be extended near the
point using cutoffs, and the differentiated output
depends only on their germs.

Frame agreement at a unit is per bundle, not an
identification of source and target transition
matrices.  Frame and projector variation can enter
the connection, conjugation, edge, and coefficient
terms; all are included in the estimates above.
Equation~\eqref{eq:open-diff-corner} then gives the
identity on the original global projective summands.
\end{proof}

\begin{rem}\label{rem:open-transfer-scope}
The transferred and cohesive differentials act on the same
graded module.  The comparison between them is constructed
in real vector spaces, and its identity diagonal will give
a real-linear chain inverse in
Proposition~\ref{prop:open-transfer-chain-inverse}.  Thus the
differentials are related by conjugation with this particular
$T$, rather than identified in advance.

All operators preserve jointly smooth passive finite-dimensional
parameters under fixed global body bounds. The later choice of
kernel primitives and counits is made objectwise.
\end{rem}

\section{The all-open affine Riemann--Hilbert theorem}
\label{sec:affine-nonregular-rh}

We now pass from the prescribed transfer to the actual
raw Hom complexes and then to the actual full-MC
integration functor.  Throughout, the notation
$\Rep_{\mathrm{sm}}^{\mathrm{perf}}$ means the globally
bounded finite-projective raw category defined in
Section~\ref{sec:representations}.  The body cohomology
ranks may jump.

\begin{prop}[Finite chain inverse of the fixed transfer]
\label{prop:open-transfer-chain-inverse}
For an invertible-edge pair $R,S$, put
\[
 [\ell,u]=[a_S-b_R,b_S-a_R].
\]
Under the fixed graded Hom identification,
\[
 T=1+N,\qquad N^{u-\ell+1}=0,
\]
and
\begin{equation}\label{eq:open-transfer-inverse}
 T^{-1}=\sum_{j=0}^{u-\ell}(-N)^j
\end{equation}
is a two-sided real-linear chain inverse.
\end{prop}

\begin{proof}
Proposition~\ref{prop:open-fixed-comparison-maps}
makes $T$ internally nonincreasing, and
Theorem~\ref{thm:open-transfer-identity-diagonal}
makes its diagonal the identity.  Thus $N=T-1$
strictly lowers internal degree.  After $u-\ell+1$
iterations its value lies below the available
interval, so it is zero.  The finite geometric sum
is therefore a two-sided graded inverse.

The already proved chain-map identity is
$d_{\Coh}T=TD_{\mathrm{tr}}$.
Multiplication by the graded inverse gives
$T^{-1}d_{\Coh}=D_{\mathrm{tr}}T^{-1}$.
Thus \eqref{eq:open-transfer-inverse} is a chain
inverse.  The argument does not require either
differential to preserve the internal-degree
filtration.
\end{proof}

\begin{thm}[All-pair actual Hom comparison]
\label{thm:affine-open-hom}
For every ordinary open $O\subseteq\R^2$ and every
ordered pair of actual raw coefficients $R,S$,
\[
 F_{R,S}:
 \Hom^\bullet_{\mathrm{raw}}(R,S)
 \longrightarrow
 \Hom^\bullet_{\Coh(B_O)}(FR,FS)
\]
is a quasi-isomorphism in every integer degree.
Its actual global smooth kernel is acyclic.
\end{thm}

\begin{proof}
First suppose the edges are actual chain
isomorphisms.  The map $S_{\mathrm{mod}}$ of
Proposition~\ref{prop:open-fixed-comparison-maps}
is a chain-homotopy equivalence, and $T$ is a
chain isomorphism by
Proposition~\ref{prop:open-transfer-chain-inverse}.
Therefore $TS_{\mathrm{mod}}$ is a quasi-isomorphism.
Equation~\eqref{eq:open-fixed-TS-homotopy} identifies
its map on cohomology with that of the actual
$F_{R,S}$.

For arbitrary original edges, take the closed
rectification isomorphisms
\[
 E_R:R^{\mathrm e}\to R,\qquad
 E_S:S^{\mathrm e}\to S
\]
of Theorem~\ref{thm:open-edge-rectification}.
They give chain isomorphisms
\[
 \phi\longmapsto E_S^{-1}\star\phi\star E_R
\]
on raw Homs, and
\[
 \Phi\longmapsto F(E_S)^{-1}\Phi F(E_R)
\]
on cohesive Homs.  Strict dg differentiation gives
the commuting square
\[
 \begin{tikzcd}[column sep=large]
 \Hom_{\mathrm{raw}}(R,S)
 \ar[r,"F_{R,S}"] \ar[d,"\cong"']
 &
 \Hom_{\Coh}(FR,FS)\ar[d,"\cong"]
 \\
 \Hom_{\mathrm{raw}}(R^{\mathrm e},S^{\mathrm e})
 \ar[r,"F_{R^{\mathrm e},S^{\mathrm e}}"']
 &
 \Hom_{\Coh}(FR^{\mathrm e},FS^{\mathrm e}).
 \end{tikzcd}
\]
The differentiated gauges are retained; they need
not be identities.  The bottom map is a
quasi-isomorphism by the invertible-edge case,
and hence so is the top map.

Lemma~\ref{lem:open-graded-hom-lifting} supplies a
short exact sequence of actual global complexes
\[
 0\longrightarrow K_{R,S}
 \longrightarrow\Hom_{\mathrm{raw}}(R,S)
 \xrightarrow{F_{R,S}}\Hom_{\Coh}(FR,FS)
 \longrightarrow0.
\]
The long exact cohomology sequence now gives
$H^m(K_{R,S})=0$ for every integer $m$.
The acyclicity therefore holds in the global smooth
kernel itself.
\end{proof}

\begin{cor}[Calibrated identity-body counits]
\label{cor:open-calibrated-counits}
Every actual raw object $R$ admits a global smooth
normalized closed isomorphism
\[
 q_R:I(FR)\longrightarrow R,\qquad
 (q_R)_0=1,\qquad Fq_R=1.
\]
The map $q_R$ is chosen separately for each $R$.
\end{cor}

\begin{proof}
Put $R'=I(FR)$.  Theorem~\ref{thm:open-literal-retraction}
gives $FR'=FR$, and the two objects have the same
underlying graded bundle and body differential.
On that common body the identity zero-cochain has
\[
 D_{R',R}1=r_R-r_{R'}\in K_{R',R}^1,
 \qquad D(D1)=0.
\]
By Theorem~\ref{thm:affine-open-hom}, choose
$c\in K_{R',R}^0$ with $Dc=D1$.  Then
$q_R=1-c$ is closed and $Fq_R=1$.
Since $F$ is the identity on body components,
$c_0=0$ and $(q_R)_0=1$.

Every positive component of $c$ lowers internal
degree.  If the body's width is $w_R=b_R-a_R$,
then $c^{\star(w_R+1)}=0$ in the common underlying
graded cochain algebra.  Form there
\[
 q_R^{-1}=\sum_{j=0}^{w_R}c^{\star j},
\]
and assign the resulting cochain the reverse
object type $R\to R'$.  The two inverse identities
hold in that graded algebra; the Hom Leibniz rule
then makes the inverse closed.
The powers of $c:R'\to R$ are taken in the common underlying
graded cochain algebra before assigning the reverse object type.
\end{proof}

\begin{thm}[All-open affine Riemann--Hilbert]
\label{thm:affine-open-rh}
For every ordinary open $O\subseteq\R^2$, the actual
functors
\[
 I=\Psi_O^*\mathsf{IN}_{E_O}:
 \Coh(B_O)\longrightarrow
 \Rep_{\mathrm{sm}}^{\mathrm{perf}}(K[O])
\]
and
\[
 \mathsf{IN}_{E_O}:
 \Coh(B_O)\longrightarrow
 \Rep_{\mathrm{sm}}^{\mathrm{perf}}(\Int(E_O))
\]
are smooth strictly unital $\ai$ quasi-equivalences.
The bodies are globally bounded and degreewise
finite projective, with no regularity assumption
on their fibre cohomology.

Every finite-model raw object has the calibrated
actual isomorphism of
Corollary~\ref{cor:open-calibrated-counits}.
Every full-MC raw object has a global smooth
normalized closed actual isomorphism from the
integration of a cohesive object, with identity
body.
\end{thm}

\begin{proof}
For cohesive objects $P,Q$, literal retraction gives
\[
 F_{IP,IQ}I_1=1.
\]
The first map is a quasi-isomorphism by
Theorem~\ref{thm:affine-open-hom}, so $I_1$ is a
quasi-isomorphism in every Hom degree.
Corollary~\ref{cor:open-calibrated-counits} supplies
essential surjectivity on $H^0$.  Strict unitality
and smoothness are properties of the actual
integration functor.  Here $I_1$ is used as a chain map;
the higher Taylor maps govern compatibility with composition.

For the full-MC target, write
\[
 i=\Phi_O^*,\qquad p=\Psi_O^*,\qquad pi=1.
\]
The range-preserving homotopy of
Theorem~\ref{thm:affine-open-comparison} gives,
by the pointed prism, closed invertible maps
\[
 U_X:X\longrightarrow ipX,\qquad (U_X)_0=1.
\]
For a homogeneous cochain $\alpha:X\to Y$, its
all-degree identity is
\[
 D\mathcal P_{H_O}\alpha+\mathcal P_{H_O}D\alpha
 =
 U_Y\star\alpha-ip(\alpha)\star U_X.
\]
For objects pulled back from $K[O]$, the homotopy
being over $K[O]$ makes $U_{iR}=1$.
The identity therefore gives the Hom retraction,
and arbitrary objects are reduced to these
pullbacks by $U_X$.  This is precisely the
all-degree dg quasi-equivalence of
Corollary~\ref{cor:affine-open-raw-comparison}.
Since
\[
 p(\mathsf{IN}_{E_O})_1=I_1,
\]
the actual full-MC unary maps are quasi-isomorphisms.

For essential surjectivity, let $X$ be a full-MC
raw object, put $R=pX$, and take $P=FR$.
Let $q_K:I(P)\to R$ be its finite-model counit.
The required map is
\begin{equation}\label{eq:open-full-mc-counit}
 U_X^{-1}\star i(q_K)\star U_{\mathsf{IN}_{E_O}(P)}
 :
 \mathsf{IN}_{E_O}(P)\longrightarrow X.
\end{equation}
Its types are
\[
 \mathsf{IN}_{E_O}(P)
 \longrightarrow ip\mathsf{IN}_{E_O}(P)=iI(P)
 \longrightarrow iR=ipX
 \longrightarrow X.
\]
All three factors are closed, normalized, invertible,
and have identity body.  This proves the stated
full-MC reconstruction.
\end{proof}

\begin{cor}\label{cor:open-two-term-kernels}
For literal bodies in degrees $0,1$, the actual
differentiation kernel vanishes in cohomology in
every degree.  In particular, this holds separately
for pairs whose body differentials are injective
on an open dense locus and for pairs whose body
differentials are surjective on an open dense locus.
The all-pair theorem applies to arbitrary body differentials.
\end{cor}

\begin{proof}
This is the specialization of
Theorem~\ref{thm:affine-open-hom}.
\end{proof}

\begin{rem}[Coefficient categories and choices]
\label{rem:open-rh-scope}
The theorem concerns the global smooth Hom complexes of
globally bounded, degreewise finite-projective real cohesive
and raw coefficients on each ordinary open. Full-MC simplices
retain their whole base image. These are the concrete coefficient
categories, rather than $\Perf(B_O)$ or a Hom-sheafified replacement.

The comparison $T$ and its inverse conjugate the two Hom
differentials as real-linear chain maps. Integration is an
$\ai$-functor, with compatibility with composition governed by
its higher Taylor maps. Counits are constructed objectwise,
and the sections are chosen separately on each open. Naturality
over coefficient moduli and strict compatibility with open
restriction are additional conditions on these choices.
The kernel argument chooses a primitive in the actual global
complex for each closed kernel cochain.

The classification and vertical-family results below assume
regular body cohomology and use the minimalization and finite
MC-lifting arguments of
Appendices~\ref{sec:affine-hom} and~\ref{sec:affine-reconstruction}.
Questions of descent and broader reconstruction are discussed
in Section~\ref{sec:programme}.
\end{rem}

\section{Nonregular coefficients and a nonrectangular open}
\label{sec:nonregular-examples}

Theorem~\ref{thm:affine-open-rh} permits body cohomology
to jump, including at an infinitely flat rank drop.
The following examples show why neither regular minimalization
nor a generic injectivity or surjectivity hypothesis suffices
for the coefficient theorem.

\subsection{An infinitely flat coefficient with a higher term}

On $\R^2$, define
\[
 s(y)=
 \begin{cases}
 e^{-1/y^2},&y\ne0,\\
 0,&y=0,
 \end{cases}
 \qquad
 a(y)=
 \begin{cases}
 e^{-1/y^2}/y,&y\ne0,\\
 0,&y=0.
 \end{cases}
\]
Then $s=ya$ everywhere.

\begin{example}[Infinitely flat higher-MC coefficient]
\label{ex:open-infinitely-flat}
Let $V^0,V^1$ be trivial real lines. In homogeneous frames
ordered by degrees zero and one, put
\[
 d=\begin{pmatrix}0&0\\1&0\end{pmatrix},
 \qquad
 h=\begin{pmatrix}0&1\\0&0\end{pmatrix}.
\]
In forms-first notation define
\begin{equation}\label{eq:open-flat-example}
 \mathcal D
 =
 Q+s(y)d+a(y)\theta_1\,1+\eta h.
\end{equation}
This is a flat cohesive coefficient.  Its body
cohomology is zero off $y=0$ and has dimension one
in each of degrees zero and one on $y=0$.
\end{example}

\begin{proof}
Every derivative of $s$ or $a$ away from zero is
a finite sum of terms $y^{-N}e^{-1/y^2}$.
These tend to zero faster than every power of $y$.
Successively extending each derivative by zero
therefore proves that $s$ and $a$ are smooth and
infinitely flat at zero.

The matrices have degrees $|d|=1$, $|h|=-1$ and satisfy
\[
 d^2=h^2=0,\qquad dh+hd=1.
\]
Thus each of the three displayed terms of $\mathcal D-Q$
has total degree one.
Since $s,a$ depend only on $y$, $Qs=Qa=0$, and
\[
 Q(a\theta_1\,1)=-ay\eta\,1.
\]
The three squares vanish:
\[
 (sd)^2=0,\qquad
 (a\theta_1\,1)^2=0,\qquad
 (\eta h)^2=0.
\]
The mixed products with the one-form term are
\[
 (sd)(a\theta_1\,1)=-sa\theta_1d,\qquad
 (a\theta_1\,1)(sd)=as\theta_1d,
\]
and
\[
 (a\theta_1\,1)(\eta h)=a\theta_1\eta h,\qquad
 (\eta h)(a\theta_1\,1)=-a\eta\theta_1h.
\]
They cancel in pairs.  The remaining two products are
\[
 (sd)(\eta h)+(\eta h)(sd)
 =s\eta(dh+hd)=s\eta\,1.
\]
Consequently the full curvature is
\[
 \mathcal D^2=\eta(s-ay)\,1=0.
\]

The body complex is
$\R\xrightarrow{s(y)}\R$ in degrees zero and one.
For $y\ne0$, this is an isomorphism.
For $y=0$, its differential is zero.
The stated fibre-cohomology dimensions follow. These concern
the complexes at individual points; cohomology of the
global-sections complex is a different calculation.
\end{proof}

The term $\eta h$ is nonzero and cannot be deleted
while leaving the other displayed components fixed:
the resulting curvature would be $-s\eta\,1$.
The actual integration of
\eqref{eq:open-flat-example} is covered by
Theorem~\ref{thm:affine-open-rh}, as are its restrictions
to arbitrary ordinary opens.

\begin{example}[A nonregular body which is nowhere injective
or surjective]
\label{ex:open-nongeneric-direct-sum}
Take the direct sum of
Example~\ref{ex:open-infinitely-flat} with the trivial
coefficient on one additional line in degree zero
and one additional line in degree one.
Its body differential, in degree-ordered bases, is
\[
 \begin{pmatrix}s(y)&0\\0&0\end{pmatrix}:\R^2\to\R^2.
\]
This body map is nowhere injective and nowhere
surjective.  Its kernel and cokernel each have
dimension one off $y=0$ and dimension two on $y=0$.
The coefficient is flat and satisfies the all-open
Riemann--Hilbert theorem.
\end{example}

\begin{proof}
The superconnection is block diagonal, so its
curvature is the direct sum of the zero curvature
computed above and the zero curvature of the
trivial coefficient.  The body matrix has rank one
for $y\ne0$ and rank zero for $y=0$.  The kernel and
cokernel dimensions follow immediately.
In particular, neither the injectivity locus nor
the surjectivity locus is dense; both are empty.
The graded body is nevertheless bounded and finite
free, so Theorem~\ref{thm:affine-open-rh} applies.
\end{proof}

These examples lie beyond the regular normal-form
classification of Section~\ref{sec:applications}.
Together they exhibit rank jumps and a body map with
neither an injective nor a surjective locus. The
all-pair theorem instead uses a comparison that
retains the body differential.

\subsection{A finite coefficient invisible at the origin}

\begin{example}[Nonconservativity of origin observation]
\label{ex:open-origin-observer}
On the whole plane let $V^0=Ae_0$, $V^1=Ae_1$, let
$d(e_0)=e_1$, $d(e_1)=0$, and set
\[
 P=(V\otimes_A B,D),\qquad D=Q+(y-1)d.
\]
This is a bounded finite-free cohesive coefficient. Its body
complex at the origin is contractible, but $P$ is not
contractible in $\Coh(B)$. Thus the functor which observes
only the body complex at the origin does not reflect
isomorphisms on the $H^0$ categories, even on these
finite coefficients.
\end{example}

\begin{proof}
The curvature is zero because $Q(y-1)=0$ and $d^2=0$.
At the origin the body differential is $-d$; the map
$e_1\mapsto-e_0$, $e_0\mapsto0$ contracts it.
At $(0,1)$ the body differential vanishes on a nonzero
two-term graded vector space. A cohesive contraction of $P$
would induce a contraction of this body complex, which is
impossible. Consequently the map $0\to P$ is not a homotopy
equivalence, although its origin-body image is one.
\end{proof}

Integration preserves the body, while
Theorem~\ref{thm:affine-open-rh} reflects the nonzero object
class of $P$. Thus origin observation still has a nontrivial
kernel after integration. A reconstruction through that
observer must account for this kernel.

\subsection{A scalar class on a nonrectangular open}

The scalar cohomology of an arbitrary open is not
obtained by restricting the full-plane answer.
In degree one the difference can already be seen
directly.

\begin{prop}[Degree-one scalar cohomology]
\label{prop:open-scalar-degree-one}
For every ordinary open $O\subseteq\R^2$, put
$\xi=x\theta_1+y\theta_2$.  Then
\[
 H^0(B_O,Q)=\ker\partial_x,\qquad
 H^1(B_O,Q)\cong
 C^\infty(O)/\partial_xC^\infty(O),
\]
where the degree-one class of $b$ is represented
by $b\xi$.  The latter is a quotient of real vector
spaces, not an ideal quotient.
\end{prop}

\begin{proof}
Direct calculation gives
\[
 Q(P\theta_1+R\theta_2)
 =
 (-yP+xR)\eta+
 (-yP_x+xR_x+R)\theta_1\theta_2.
\]
If $-yP+xR=0$, then away from the origin there is
a unique smooth
\[
 b=\frac{xP+yR}{x^2+y^2}
\]
with $P=xb$ and $R=yb$.
If the origin belongs to $O$, choose a rectangle
around it.  The relation implies $P(0,y)=0$,
first for $y\ne0$ and then by continuity.
Hadamard division on that rectangle gives $P=xb$,
and $x(R-yb)=0$ gives $R=yb$ everywhere there.
Uniqueness glues this local expression to the
one away from the origin.
For these values the second displayed curvature
coefficient vanishes as well.  Thus the degree-one
cycles are exactly $b\xi$.

Multiplication by $\xi$ is injective on smooth
functions: a function killed by both $x$ and $y$
vanishes away from the origin and then everywhere
by continuity.  Since $Qf=f_x\xi$, the boundaries
are exactly $(\partial_xf)\xi$.
The same injectivity gives the degree-zero formula.
\end{proof}

\begin{example}[A nonzero class on a punctured open]
\label{ex:open-punctured-scalar}
Let
\[
 O=\R^2\setminus\{(1,0)\},\qquad
 b(x,y)=\frac{1}{(x-1)^2+y^2}.
\]
Then $[b\xi]\ne0$ in $H^1(B_O,Q)$.
For an admissible edge from $(x,y)$ to $(x',y)$,
its actual scalar integration is
\begin{equation}\label{eq:open-punctured-edge-integral}
 I_1(b\xi)(g,m)
 =
 \int_x^{x'}b(u,y)\,du
 =
 (x'-x)\int_0^1b(x+t(x'-x),y)\,dt.
\end{equation}
\end{example}

\begin{proof}
The function $b$ is smooth on $O$.
If $b=\partial_x f$ with $f\in C^\infty(O)$, then
for $0<|y|<1/2$,
\[
 \begin{split}
 f(3/2,y)-f(1/2,y)
 &=
 \int_{1/2}^{3/2}
 \frac{dx}{(x-1)^2+y^2}\\
 &\ge
 \int_{1-|y|}^{1+|y|}
 \frac{dx}{(x-1)^2+y^2}
 =
 \frac{\pi}{2|y|}.
 \end{split}
\]
The two endpoint functions are smooth near $y=0$,
so their difference is locally bounded, a
contradiction.  Proposition~\ref{prop:open-scalar-degree-one}
therefore gives the nonzero class.

On a full-MC edge $\sigma$, the base equation is
$\sigma^*\xi=d_\Delta x_\sigma$.
The unary scalar kernel is integration, and the
selected base path lies in one horizontal interval
component of $O$.  A one-form $b(x,y)\,dx$ on that
interval has an ordinary primitive, so its integral
depends only on the two endpoints.  This gives
\eqref{eq:open-punctured-edge-integral}.
The hull condition ensures that the displayed
straight segment lies in $O$, making the last
formula jointly smooth on the admissible arrow
manifold.
\end{proof}

Thus scalar cohomology can change on a nonrectangular open:
degree-one cohomology vanishes on the full plane, whereas
this open has the displayed nonzero class. The all-open
theorem compares the actual cohomology on each open.

\section{Two-term normal forms and ordinary non-descent}
\label{sec:applications}

Throughout this section $M=\R^2$, $B=B_{E_{\mathrm{aff}}}$ and $K$
is the actual crossed-module nerve. All equivalences of coefficient
objects mean isomorphisms in the $H^0$ category of the indicated
dg category.

\subsection{Normal forms for cohesive coefficients}

Let $\Coh^{\mathrm{reg}}_{[0,1]}(B)$ be the full subcategory whose
internal cohomology is supported in degrees $0,1$. The original
coefficient complex may have acyclic terms outside these degrees.
For raw representations, the same amplitude condition concerns
vertex cohomology, not the support of the original graded bundle.
For finite-dimensional $V^0,V^1$ and a linear map $N:V^1\to V^0$,
let $P_N$ be the cohesive module on the constant graded bundle
$V=V^0\oplus V^1$, defined on constant vectors by
\begin{equation}\label{eq:two-term-constant}
 D_N(v_0)=0,\qquad D_N(v_1)=N(v_1)\eta.
\end{equation}
It is extended by the right-module Leibniz rule. Extend $N$ by zero
on $V^0$; then $N^2=0$, and $Q\eta=0$ proves flatness.
Write $\mathbf1_j$ for the trivial coefficient with one generator
of internal degree $j$ and zero differential on that generator.
The elementary eta-block $\mathsf E_\eta$ has one generator in
each of degrees $0,1$ and $D(v_1)=v_0\eta$.

Taking CE degree zero and then evaluating at $m\in M$ gives a dg
functor from cohesive modules to finite complexes. The analogous
vertex functor on raw representations is dg by
\eqref{eq:representation-differential}. Thus an $H^0$-equivalence
induces a chain homotopy equivalence of vertex complexes.
This uses actual homotopies, not exactness of point evaluation.

\begin{thm}[Two-term normal forms]\label{thm:two-term-normal-forms}
Every object of $\Coh^{\mathrm{reg}}_{[0,1]}(B)$ is $H^0$-isomorphic
to some $P_{N_0}$ with constant real matrix $N_0$.
Two such objects are isomorphic exactly when their triples
\[
 (d_0,d_1,r)
 =(\dim V^0,\dim V^1,\operatorname{rank}N_0)
\]
agree. The possible triples are the nonnegative integers with
$r\le\min(d_0,d_1)$, and the corresponding normal form is
\begin{equation}\label{eq:two-term-decomposition}
 \mathsf E_\eta^{\oplus r}
 \oplus\mathbf1_0^{\oplus(d_0-r)}
 \oplus\mathbf1_1^{\oplus(d_1-r)}.
\end{equation}
\end{thm}

\begin{proof}
Apply Lemmas~\ref{lem:affine-finite-minimalization} and
\ref{lem:affine-minimal-gauge}. The resulting cohomology bundle
is concentrated in degrees $0,1$ and is trivial on $\R^2$.
Its internal differential and CE-degree-one connection vanish.
A remaining connection component of CE degree $p$ has internal
degree $1-p$. Hence the only possible component has $p=2$ and
maps $V^1$ to $V^0$.

In graded frames, write its coefficient matrix as
$N(x,y)\eta+T(x,y)\theta_1\theta_2$.
Products of these negative internal endomorphisms vanish.
The scalar differential calculation
\eqref{eq:affine-ce-even} therefore gives
\[
 (N_x-T)(x\theta_1+y\theta_2)\eta=0.
\]
Entrywise, multiplication by both $x$ and $y$ kills $N_x-T$.
It is zero away from the origin and then at the origin by
continuity. Thus the remaining coefficient is
\[
 b_N=N\eta+N_x\theta_1\theta_2.
\]

Set $N_0=N(0,0)$ and define smooth matrices
\[
 P=-\int_0^1N_y(tx,ty)\,dt,\qquad
 R=\int_0^1N_x(tx,ty)\,dt.
\]
The fundamental theorem of calculus gives
$N-N_0=-yP+xR$. These compact-interval integrals are globally
smooth; no division by $x$ or $y$ is involved.
For $c=P\theta_1+R\theta_2$, equation~\eqref{eq:affine-ce-odd} gives
\begin{equation}\label{eq:two-term-gauge-primitive}
 Qc=(N-N_0)\eta+N_x\theta_1\theta_2=b_N-N_0\eta.
\end{equation}

Denote the reduced differential by $D_{b_N}$ and regard $v^j$
as the row of degree-$j$ basis vectors.
The right-linear degree-zero map from the $b_N$ module to $P_{N_0}$,
\[
 g(v^0)=v^0,\qquad g(v^1)=v^1+v^0c,
\]
has inverse with the negative correction, since its off-diagonal
part squares to zero. Its orientation is fixed by
\[
 D_{N_0}g(v^1)=v^0(N_0\eta+Qc)
             =v^0b_N=gD_{b_N}(v^1).
\]
The equality on $v^0$ is immediate, and the Leibniz rule extends
these identities to the whole module. This proves existence.

For uniqueness, a general closed degree-zero map
$\varphi:P_N\to P_{N'}$ between constant-residue models has the form
\[
 \varphi(v^0)=w^0S_0,\qquad
 \varphi(v^1)=w^1S_1+w^0C,
 \qquad C=P\theta_1+R\theta_2,
\]
where $S_0,S_1$ are smooth function matrices.
Nonnegative CE degree excludes any other component.
Closedness gives
\begin{equation}\label{eq:two-term-intertwining}
 QS_0=QS_1=0,\qquad N'\eta S_1+QC=S_0N\eta.
\end{equation}
The eta coefficient of $QC$ is $-yP+xR$, which vanishes at
the origin. Consequently
\[
 N'S_1(0)=S_0(0)N.
\]
If $\varphi$ is invertible in $H^0$, its vertex map is a homotopy
equivalence between zero-differential complexes. Thus $S_0(0)$
and $S_1(0)$ are actual isomorphisms. The dimensions and ranks agree.
Conversely, equal dimensions and ranks give constant invertible
$S_0,S_1$ intertwining $N,N'$. Setting $C=0$ in
\eqref{eq:two-term-intertwining} gives a strict cohesive isomorphism.
Putting $N_0$ in rank normal form proves
\eqref{eq:two-term-decomposition}.

Zero matrices cover $r=0$, and empty graded pieces cause no change
in the argument.
If both cohomology bundles vanish, the initial finite minimalization
contracts the object to zero. The same argument removes all
acyclic tails in the original bounded complex.
\end{proof}

The matrix $N_0$ depends on frames and reductions. Its
$\mathrm{GL}(V^0)\times\mathrm{GL}(V^1)$ orbit, equivalently its
rank together with the two dimensions, is the object invariant.
The classification concerns objects up to $H^0$-isomorphism;
their morphism complexes retain the dg structure.

\subsection{The actual integrated normal forms}

\begin{cor}\label{cor:two-term-representations}
The actual affine RH functors restrict to $\ai$ quasi-equivalences
on the regular full subcategories of internal cohomology amplitude
$[0,1]$. On $\y K$, their normal forms are
\begin{equation}\label{eq:two-term-integrated}
 R_{N_0}:\quad
 r_0=0,\quad r_1=1,\quad r_2=-h_{012}N_0,\quad r_p=0\quad(p\ge3).
\end{equation}
The same triples $(d_0,d_1,r)$ classify $H^0$-isomorphism classes
on $\y K$ and on full MC integration.
\end{cor}

\begin{proof}
The actual functors preserve vertex complexes. If a regular target
has internal cohomology in degrees $0,1$, the preimage and counit
of Theorem~\ref{thm:affine-rh} have homotopy equivalent vertex
complexes. That preimage therefore already belongs to
$\Coh^{\mathrm{reg}}_{[0,1]}(B)$.
Together with the all-degree Hom comparison, this proves the
quasi-equivalence on these full subcategories.

For $P_{N_0}$, the argument of Proposition~\ref{prop:in-affine-probe}
applies to the whole matrix. Every connection word with two
occurrences of $N_0$ vanishes because $N_0^2=0$.
The actual unary kernel is $\Gamma_1=\int$, so on full MC the
two-component is $(\int_{[012]}\beta)N_0$.
Pullback by $\Psi$ and $h=-\int\beta$ give
\eqref{eq:two-term-integrated}. Theorem~\ref{thm:two-term-normal-forms}
and the restricted RH equivalences give the classification.
\end{proof}

\subsection{An actual ordinary quotient}

\begin{prop}\label{prop:ordinary-quotient}
There is an internal ordinary groupoid $\overline{\CG}$ in smooth
sheaves with object sheaf $\y M$, and a strict simplicial sheaf map
\[
 q:\y K\longrightarrow N(\overline{\CG}),\qquad q_0=\id_{\y M},
\]
obtained by identifying parallel $2$-isomorphic $1$-arrows.
Its arrow object is the sheafification of the presheaf of
these classes.
\end{prop}

\begin{proof}
For each test $U$, evaluate the strict affine $2$-groupoid and take
$2$-isomorphism classes of its parallel $1$-arrows. Horizontal
composition and inverses preserve this relation, so the classes
form an ordinary groupoid. In source coordinates, $(g,m)$ is
related to $(g\exp(hZ_m),m)$; equivalently the multiplier on the
left is $\exp(hZ_{gm})$. Both arrows have target $gm$.

Restriction pulls back the actual $2$-isomorphisms, giving a
presheaf of groupoids with object presheaf $\y M$.
Let $P_1$ be its arrow presheaf and $a$ ordinary sheafification.
Left exactness gives
\[
 a(P_1\times_{\y M}P_1)
   \cong aP_1\times_{\y M}aP_1.
\]
Thus all groupoid structure maps sheafify, with the object
sheaf unchanged, to $\overline{\CG}$.

Send a Duskin simplex to the consecutive classes
$[g_{01}],\ldots,[g_{n-1,n}]$.
Its triangle $2$-isomorphisms imply
$[g_{ik}]=[g_{jk}][g_{ij}]$, so inner faces become composition,
outer faces omit end arrows, and degeneracies insert identities.
This is a simplicial presheaf map to the ordinary nerve.
Sheafification commutes with its finite nerve products and
therefore yields the stated $q$.
\end{proof}

This is an ordinary groupoid object in the same smooth sheaf topos
used for raw coefficients. In particular its object sheaf is represented,
as required by Definition~\ref{def:smooth-representations}.
It records ordinary $2$-isomorphism classes of $1$-arrows.
Comparison with AS holonomy, refined higher holonomy, a
geometric Morita quotient or an infinity-topos truncation
requires additional geometric input.

\subsection{The higher-isotropy obstruction}

There is a strict inclusion
\[
 \iota:T=B^2\Z\longrightarrow\y K
\]
at the origin: all root edges are identities and all rooted
triangle coordinates are integers.
Equation~\eqref{eq:affine-root-coordinates}, with $Z_0=0$,
shows that every edge stays an identity and every other
triangle coordinate stays integral under the simplex operators.
Here $T_n$ is the locally constant sheaf
$\Z^{\binom n2}$, not a constant presheaf on disconnected tests.
In particular $T_0=T_1=*$ and $T_2$ contains the simplex $h_{012}=1$.
The simplicial object has all degrees.

\begin{thm}[Ordinary-groupoid non-descent]
\label{thm:ordinary-nondescent}
Let $\CG$ be any internal ordinary groupoid in smooth sheaves
with represented object sheaf $\CG_0=\y M'$.
For any strict simplicial sheaf map
$f:\y K\to N(\CG)$ and any
$W\in\Rep_{\mathrm{sm}}^{\mathrm{perf}}(N(\CG))$,
the normal form $R_N$ with $N\ne0$ is not $H^0$-isomorphic to $f^*W$.
The same assertion holds for
$\mathsf{IN}_{E_{\mathrm{aff}}}(P_N)$ and any strict map
$\Int(E_{\mathrm{aff}})\to N(\CG)$.
The coefficient $W$ may be nonregular or nonminimal.
\end{thm}

\begin{proof}
An ordinary nerve is determined by its consecutive $1$-arrows.
Since the unique $1$-simplex of $T$ is degenerate, $f\iota$
sends every simplex to the totally degenerate simplex at
$p=f_0(0)$. Consequently $(f\iota)^*W$ is the constant
representation of the finite vertex complex $C=W_p$:
its internal differential is $d_C$, its edge transport is the
identity, and its higher components vanish by normalization
on their degenerate images.

Suppose an $H^0$-equivalence existed. After restriction, represent
it by a closed degree-zero morphism
$F:\iota^*R_N\to\underline C$.
The vertex map $F_0:(V,0)\to(C,d_C)$ is a chain homotopy
equivalence. Normalization gives $F_1=0$ because every
$1$-simplex of $T$ is degenerate.
At simplicial degree two, the actual differential
\eqref{eq:representation-differential} becomes
\begin{equation}\label{eq:ordinary-obstruction}
 d_CF_2+hF_0N=0.
\end{equation}
Indeed the sign of $r_{C,0}\star F_2$ is $(-1)^{1\cdot2}=1$,
and $-F_0\star r_{N,2}=hF_0N$.
All other terms contain $F_1$, the zero source internal
differential or a zero target higher component.

Evaluate at the specified simplex $h=1$ and apply to $v\in V^1$.
The term $F_2(v)$ can lie in $C^{-1}$, and its differential
contributes a boundary.
Equation~\eqref{eq:ordinary-obstruction} therefore says that
$F_0Nv$ is a boundary. Thus
$H^0(F_0)(Nv)=0$ in $H^0(C)$.
Injectivity of $H^0(F_0):V^0\to H^0(C)$ forces $Nv=0$
for every $v$, a contradiction. The argument evaluates the
cochain identity on the specified simplex of $B^2\Z$.

For a strict map $f:\Int(E_{\mathrm{aff}})\to N(\CG)$,
pull a hypothetical descent back along the actual global
section $\Psi$. The equality
$\Psi^*\mathsf{IN}_{E_{\mathrm{aff}}}(P_N)
=\mathsf{IN}^K(P_N)=R_N$
would give descent along $f\Psi$, contradicting the first part.
\end{proof}

\begin{cor}[Exact two-term descent criterion]
\label{cor:two-term-descent}
A regular two-term raw representation of $\y K$ lies in the
essential image of $H^0(q^*)$ exactly when its residue rank is zero.
For full MC the same criterion holds for
$q_{\mathrm{MC}}=q\Phi$.
\end{cor}

\begin{proof}
Corollary~\ref{cor:two-term-representations} puts every such object
in a unique normal-form class. Nonzero rank is excluded by
Theorem~\ref{thm:ordinary-nondescent}. At rank zero the normal form
is a sum of the trivial degree-zero and degree-one coefficients.
It is pulled back from the constant graded representation of
$\overline{\CG}$ with those dimensions.
The same constant representation works for $q\Phi$ on full MC;
also $q\Phi\Psi=q$.
\end{proof}

Here descent means membership in the specified essential image.
The classification is global and regular, and the obstruction
concerns strict maps to ordinary nerves in the raw coefficient
theory. Ordinary-open RH follows from the direct comparison in
Theorem~\ref{thm:affine-open-rh}; regular minimalization enters
here to classify the resulting objects. Effective descent and
comparison with refined higher holonomy are discussed in
Section~\ref{sec:programme}.

\section{Vertical families and coefficient bundles}
\label{sec:affine-families}

We extend the scalar parameter comparison of
Corollary~\ref{cor:affine-scalar-parameters} to regular
bundle-valued coefficients. A frame domain on the object
manifold need not contain every vertex of a simplex.
Normalized projective cochain modules give a global
description, and reconstruction then takes place over the
central algebra of parameter functions.

Throughout this section, $S$ is a finite-dimensional Hausdorff
second-countable smooth manifold. It need not be compact or
connected. Put
\[
 M_S=S\times\R^2,\qquad \pi:M_S\longrightarrow S,\qquad
 A_S=C^\infty(M_S),\qquad C_S=C^\infty(S).
\]
All coefficient complexes are globally bounded and have
finite-projective global section modules. Regularity continues
to mean locally constant internal cohomology ranks.
The empty parameter manifold gives zero coefficient modules;
the assertions below include that case.

\subsection{The vertical resolution and product integration}

Let $E_S$ be the graded-bundle extension of $E_{\mathrm{aff}}$
with zero anchor in the $S$ directions. Thus
\[
 \partial\varepsilon=-ye_1+xe_2,\qquad
 \rho(e_1)=x\partial_x,\qquad \rho(e_2)=y\partial_x.
\]
Its CE algebra is
\begin{equation}\label{eq:family-ce}
 \begin{gathered}
 B_S=A_S[\theta_1,\theta_2,\eta],\qquad
 |\theta_i|=1,\quad|\eta|=2,\\
 Qf=f_x(x\theta_1+y\theta_2),\qquad
 Q\theta_1=-y\eta,\qquad
 Q\theta_2=\theta_1\theta_2+x\eta,\qquad Q\eta=0.
 \end{gathered}
\end{equation}
This is not the inverse-image algebroid that adds $TS$.
Write $K_{S,n}=S\times K_n$, so
$\y K_S=c(\y S)\times\y K$, where $c$ denotes the constant
simplicial object, not the smooth singular complex.

\begin{lem}\label{lem:family-product}
The vertical extension is a global geometric resolution of
$\langle x\partial_x,y\partial_x\rangle$ on $M_S$.
There is a natural identification
\[
 \Int(E_S)=c(\y S)\times\Int(E_{\mathrm{aff}}).
\]
Taking the product of $\Phi,\Psi,H$ with $c(\y S)$ gives
maps $\Phi_S,\Psi_S,H_S$, respectively. The homotopy $H_S$
remains relative to the whole $c(\y M_S)$.
The resulting actual raw pullbacks are dg quasi-equivalences,
also on regular full subcategories.
\end{lem}

\begin{proof}
Suppose $xf+yg=0$ in $A_S$. At $x=0$, the function $g$ vanishes
for $y\ne0$, and then for $y=0$ by continuity. The globally smooth
parameterized Hadamard formula
\[
 h(s,x,y)=\int_0^1 g_x(s,tx,y)\,dt
\]
gives $g=xh$. The original relation now gives $f=-yh$ on
$x\ne0$, and hence everywhere. The map $h\mapsto(-yh,xh)$ is
injective by the same continuity argument. This proves exactness
over $A_S$.

For a full MC simplex, the zero $S$-anchor makes its $S$-coordinate
constant in the connected simplex variables. On a test manifold
$U$ that coordinate is precisely a smooth map $U\to S$; all
remaining equations are those of $E_{\mathrm{aff}}$. This proves
the displayed identification in every degree, naturally in
tests and simplex operators. Products of the actual affine
comparison maps and homotopy give the product retraction.
The product homotopy remains relative to the whole $c(\y M_S)$,
not merely to the endpoint vertices.
Theorem~\ref{thm:raw-representation-invariance} therefore applies.
The vertex maps are identities, so regularity is preserved.
\end{proof}

\subsection{A global normalized projective raw dictionary}

Put $C=C^\bullet_{\mathrm{norm,sm}}(\y K_S)$, with its full
cochain differential $d_{\mathrm{full}}$ and Alexander--Whitney
product. Distinguish this differential from the inner-face
operator $\delta_{\mathrm{in}}$ of
\eqref{eq:inner-differential}. For a finite bounded graded real
vector space $W$, the constant raw representation gives
\[
 \mathbb F_W=\Hom_{\mathrm{raw}}(\mathbf1,\underline W)
             =C\otimes_\R W,\qquad
 (c\otimes w)a=(-1)^{|w||a|}(c\smile a)\otimes w.
\]
Homogeneous module maps use the signless right-linearity
convention. A finite graded idempotent
$p\in A_S\otimes\End^0_\R(W)$ presents a vertex bundle
$V=\im p$. Define $\mathbb F_V=p\mathbb F_W$.
Its components are
\begin{equation}\label{eq:family-normalized-module}
 \mathbb F_V^{k,j}
   =\Gamma_{\mathrm{norm}}(K_{S,k},v_0^*V^j).
\end{equation}
For $k>0$, the subscript means vanishing on all degeneracies.
These are the graded modules induced from $V$ over $A_S$
through the displayed finite projectors.

\begin{lem}[Global projective raw dictionary]
\label{lem:family-projective-dictionary}
Raw representations on globally bounded finite-projective
vertex bundles are dg equivalent to flat right superconnections
on the modules $\mathbb F_V$. The equivalence identifies their
actual normalized Hom complexes and preserves the regular
full subcategories.
\end{lem}

\begin{proof}
Multiplication by $p(v_0)$ commutes with degeneracies, giving
\eqref{eq:family-normalized-module}. For presentations
$V=\im p\subset W$ and $V'=\im p'\subset W'$, finite free-module
linearity gives
\begin{equation}\label{eq:family-hom-corner}
 \Hom_C(\mathbb F_V,\mathbb F_{V'})
 =p'\Hom_C(\mathbb F_W,\mathbb F_{W'})p.
\end{equation}
In cochain degree $k$, the corner condition is exactly
\[
 \alpha(\sigma)
   =p'(v_0\sigma)\alpha(\sigma)p(v_k\sigma).
\]
Indeed, postcomposition multiplies at the first vertex, whereas
precomposition uses the last vertex in the cochain product.
This supplies the bundle type
$\Hom(v_k^*V,v_0^*V')$ without requiring a common object-frame
domain for the endpoints. The free-module calculation and the
braided action give precisely the sign in
\eqref{eq:cochain-product}. In particular, a degree-$j$
zero-cochain bundle map $T$ extends as
\begin{equation}\label{eq:family-zero-cochain-sign}
 \bigl(\widehat T(c\otimes w)\bigr)(\sigma)
   =(-1)^{jk}c(\sigma)T(v_0\sigma)w,\qquad c\in C^k.
\end{equation}

The reference connection on $\mathbb F_V$ is
$\nabla_p=p\,d_{\mathrm{full}}$. The Leibniz rule follows from
$p\alpha=\alpha$, but this connection need not be flat.
For an arbitrary right superconnection $D$ on $\mathbb F_V$,
its difference $A=D-\nabla_p$ is a normalized degree-one
right-linear corner. Let $e_W$ be the ambient identity
one-cochain, which is not normalized, and set
\begin{equation}\label{eq:family-projective-structure}
 r=p\star e_W\star p+A.
\end{equation}
Here $\delta_{\mathrm{in}}p=0$, whereas
$d_{\mathrm{full}}p$ need not vanish. Applying
\eqref{eq:in-full-inner} to $\alpha\in\mathbb F_V$ therefore
gives the exact formula
\begin{equation}\label{eq:family-module-differential}
 D\alpha=\delta_{\mathrm{in}}\alpha+r\star\alpha
              -(-1)^{|\alpha|}\alpha\star e_{\mathbf1}.
\end{equation}
Projector variation contributes to curvature. For example the
reference structure $p\star e_W\star p$ has two-curvature
\[
 \kappa^p_2(012)=p_0(p_1-1)p_2,\qquad p_i=p(v_i),
\]
which can be nonzero.

The derivation identity of Lemma~\ref{lem:cochain-identities}
and flatness of the source unit yield
\[
 D^2\alpha=(\delta_{\mathrm{in}}r+r\star r)\star\alpha.
\]
The curvature $\kappa=\delta_{\mathrm{in}}r+r\star r$ is itself
an endpoint corner: $\kappa=p\star\kappa\star p$.
If $D^2=0$, test this identity on the finitely many zero-cochain
generators $pe_j$ coming from a homogeneous basis of $W$.
On a $k$-simplex it says
$\kappa_k(\sigma)p(v_k\sigma)e_j=0$ for every $j$.
These vectors span the last-vertex fiber, so $\kappa_k=0$.
Thus flatness implies the raw MC equation globally.
The reverse implication follows from the same identity.

Since $A$ is normalized, \eqref{eq:family-projective-structure}
has $r_1(s_0m)=p(m)|_{V_m}=1_{V_m}$ and normalized higher
components. The degree-two MC identities on the $h=0$
inverse-arrow triangles in $K_S$ give chain-homotopy inverses
for its edge maps. Hence it is a raw representation, including
the edge quasi-isomorphism requirement.
Conversely, \eqref{eq:family-module-differential} is a flat
right superconnection for a raw representation.
The constructions are inverse in projector presentations.

Finally, substitute \eqref{eq:family-module-differential} in
the superconnection commutator on
\eqref{eq:family-hom-corner}. The source-unit terms cancel,
leaving
\[
 D_{\Hom}\phi
 =\delta_{\mathrm{in}}\phi+r_{V'}\star\phi
                         -(-1)^{|\phi|}\phi\star r_V.
\]
This is the actual raw Hom differential, with its normalization,
and composition is the same star product. Changing a projector
presentation gives the corresponding module isomorphism and
does not change this dg equivalence.
\end{proof}

The lemma gives a concrete dg description on induced projective
modules equipped with flat superconnections, before any passage
to a derived-perfect localization.

\subsection{Regular reduction to bundles over the parameter base}

\begin{lem}[Basic minimal coefficients]
\label{lem:family-basic-minimalization}
Every regular cohesive $B_S$-module is homotopy equivalent to
one on $\pi^*H$ with zero internal differential and zero
CE-degree-one connection. Every regular raw representation
of $\y K_S$ is homotopy equivalent to one on $\pi^*H$ with
$r_0=0$ and $r_1=1_H$.
In both statements $H$ is globally bounded and
$\Gamma(S,H)$ is finite projective over $C_S$.
\end{lem}

\begin{proof}
Regularity and boundedness give a smooth contraction
$(i,p_H,h)$ of the vertex complex onto its cohomology bundle
$H_M$, as in Lemma~\ref{lem:affine-finite-minimalization}.
For the raw category, first use
Lemma~\ref{lem:family-projective-dictionary}.
Extend the contraction by
\eqref{eq:family-zero-cochain-sign}: on $k$-cochains both
$d_0$ and $h$ carry $(-1)^k$. These factors cancel in
$d_0h+hd_0=1-ip_H$, and all side conditions are preserved.
The finite transfer formulas
\eqref{eq:affine-finite-transfer} now apply to the actual
normalized projective modules. For $\Delta=D-d_0$, every
component has positive cochain degree, so $h\Delta$ and
$\Delta h$ strictly lower internal degree. The global
amplitude bound makes their inverse series uniformly finite.
The right-linearity proof uses the contraction side conditions
and the Leibniz rule, not centrality of $A_S$.
The perturbation starts with the flat superconnection $D$;
the dictionary then returns raw maps and homotopies.

In particular the transferred edge is the cohomology transport.
The degree-one part of the closed inclusion $i'$ is
\[
 -r_0(v_0)i'_1+r_1i(v_1)-i(v_0)r_{H,1}=0.
\]
The minus sign is $(-1)^{1\cdot1}$ in the first star product.
Multiplication by $p_H(v_0)$ gives
$r_{H,1}=p_H(v_0)r_1i(v_1)$.
It is an isomorphism on the zero-differential fibers.

The bundle $H_M$ is a smooth direct summand of the original
finite-projective bundle. Restrict its finite idempotent to
$S\times\{(0,0)\}$ to obtain a bundle $H$ with finite-projective
$C_S$-module of sections. Choose a grading-preserving ordinary
connection on $H_M$. Its transport along
$(s,tx,ty)$, $0\le t\le1$, gives a smooth graded isomorphism
\[
 \tau:\pi^*H\longrightarrow H_M.
\]
Local smooth dependence for this finite-interval transport
requires no compactness of $S$. This is an ordinary bundle
identification, not an algebroid path from the origin.
Pull the minimal structure back by $\tau$.

For a minimal cohesive object write its one-form part as
$\theta_1L_1+\theta_2L_2$, now with values in
$\pi^*\End^0(H)$. The $\eta$-coefficient of flatness is
$-yL_1+xL_2=0$. The global Hadamard argument of
Lemma~\ref{lem:affine-minimal-gauge}, in this fixed bundle
over each $s$, gives
\[
 L_1=xL,\qquad L_2=yL,\qquad
 L(s,x,y)=\int_0^1\partial_xL_1(s,tx,y)\,dt.
\]
Solve the intrinsic endomorphism ODE
\begin{equation}\label{eq:family-ode}
 u_x=-Lu,\qquad u(s,0,y)=1_{H_s}.
\end{equation}
It has an invertible solution for every finite $x$, jointly
smooth in $(s,x,y)$. The local parameter estimates are on
compact charts, not on all of $S$; uniqueness makes solutions
in parameter frames agree on overlaps. Since $Qu=(x\theta_1+
y\theta_2)u_x$, the transformed one-form part vanishes.

For a minimal raw object on $\pi^*H$, its triangle equation
is \eqref{eq:affine-minimal-edge}, parameterwise. At $h=0$ it
is a contravariant $\Aff_+(1)$-action, and its boundary
identity is $R(\exp(hZ_m),m)=1$. Differentiation gives the
same relation $-yL_1+xL_2=0$. The ODE
\eqref{eq:family-ode} kills these infinitesimal matrices.
The path-and-cocycle computation in
Lemma~\ref{lem:affine-target-minimalization} then kills the
whole edge representation, since $\Aff_+(1)$ is connected.
Thus $r_1=1_H$ after this gauge, not merely to first order.

For clarity, let $R_{H_M}$ denote the minimal raw object on $H_M$.
Let $R_{\mathrm{min}}$ be its transport to $\pi^*H$ through
$\tau$, and let $R_{\mathrm{triv}}$ be the result of the ODE gauge.
Their structures satisfy
\[
 r_{\mathrm{triv}}=u^{-1}\star r_{\mathrm{min}}\star u.
\]
Since $\delta_{\mathrm{in}}u=0$ for a zero-cochain,
$u:R_{\mathrm{triv}}\to R_{\mathrm{min}}$ is closed.
The map $\tau:R_{\mathrm{min}}\to R_{H_M}$ is the preceding
bundle-transport isomorphism. Thus
\begin{equation}\label{eq:family-ordinary-gauge}
 g=\tau u:R_{\mathrm{triv}}\longrightarrow R_{H_M}.
\end{equation}
Composing with the perturbed inclusion $i'$ gives the
required homotopy equivalence. If $H=0$, the finite transfer
instead contracts the object to zero. Empty graded pieces,
acyclic tails and disconnected bases cause no change.
\end{proof}

\subsection{The actual family Riemann--Hilbert theorem}

\begin{thm}[Global regular vertical-family RH]
\label{thm:family-rh}
For every $S$ as above, the actual functors
\[
 \begin{aligned}
 \mathsf{IN}^{K_S}
   &=\Psi_S^*\mathsf{IN}_{E_S}:
     \Coh^{\mathrm{reg}}(B_S)\longrightarrow
     \Rep_{\mathrm{sm}}^{\mathrm{perf,reg}}(\y K_S),\\
 \mathsf{IN}_{E_S}&:
     \Coh^{\mathrm{reg}}(B_S)\longrightarrow
     \Rep_{\mathrm{sm}}^{\mathrm{perf,reg}}(\Int(E_S))
 \end{aligned}
\]
are smooth strictly unital $\ai$ quasi-equivalences.
Their unary Hom maps are quasi-isomorphisms in every integer
degree. Every regular target has a global regular cohesive
preimage and an actual closed degree-zero homotopy equivalence
from its integration.
\end{thm}

\begin{proof}
The actual functors, their smoothness and strict units are
those of Theorem~\ref{thm:coefficient-integration}.
We prove the additional Hom and reconstruction assertions.

\emph{Finite-projective scalar extensions.}
By Lemma~\ref{lem:family-basic-minimalization}, source objects
reduce to basic minimal coefficients on $\pi^*H,\pi^*H'$.
Put
\[
 P^q=\Gamma(S,\Hom^q(H,H')).
\]
This is a finite-projective $C_S$-module. If $P^q=eC_S^N$,
extension to forms of CE degree $r$ gives $e(B_S^r)^N$,
and extension to normalized $k$-cochains gives $e(C^k)^N$.
Consequently the actual source and target coefficient spaces
are respectively
\begin{equation}\label{eq:family-projective-tensors}
 B_S^r\otimes_{C_S}P^q,\qquad C^k\otimes_{C_S}P^q.
\end{equation}
The idempotent depends only on $s$ and commutes with
normalization. These identities retain all jointly smooth
functions; they do not identify smooth functions on a product
with an algebraic tensor product of two function spaces.

We use $C_S$, rather than $A_S$, as the central scalar algebra
annihilated by the differentials. Every vertex of a simplex has
the same $s$, so $C_S$ is central in $C$ and
$d_{\mathrm{full}}(C_S)=Q(C_S)=0$. Each actual integration kernel is
separately $C_S$-multilinear: an input factor $c(s)$ is
constant throughout that simplex and its integrals.
Thus the graded matrix extension descends to
\eqref{eq:family-projective-tensors}. Its suspended and
unsuspended signs remain the distinct formulas
\eqref{eq:in-matrix-sign} and
\eqref{eq:in-unsuspended-matrix-sign}. Parameter frames are
legitimate here because an entire simplex lies over one
parameter, unlike an arbitrary object-frame domain.

\emph{All Hom degrees.}
Filter the Hom complexes increasingly by internal coefficient
degree. The filtration is finite. The source internal
differential and one-form connection vanish, and the
integrated target has $r_0=0$, $r_1=1_H$.
The associated-graded differentials are therefore $Q$ and
$d_{\mathrm{full}}$, tensored with $1_{P^q}$.
Every non-bare word in the twisted unary map contains a
strictly negative internal-degree connection component.
Its only filtration-preserving term is consequently
\[
 \operatorname{gr}^q(\mathsf{IN}^{K_S}_1)
       = I_S\otimes_{C_S}1_{P^q},
\]
where $I_S$ is the actual scalar comparison of
Corollary~\ref{cor:affine-scalar-parameters}.
Since $P^q$ is finite projective, this is a direct summand
of a finite sum of quasi-isomorphisms $I_S$.
Induction through the finite filtration proves every integer
Hom degree. As in Lemma~\ref{lem:affine-hom-reduction}, the
actual second Taylor identity and strict unit transfer the
conclusion back through the source homotopy equivalences.
This transfer uses homotopy equivalences within the regular
subcategories.

\emph{Global reconstruction.}
For a regular raw target $R$, apply
Lemma~\ref{lem:family-basic-minimalization} and write its
edge-trivialized replacement as $e_H+b$ on $\pi^*H$.
Let
\[
 \mathfrak n_H=\bigoplus_{q<0}\Gamma(S,\End^qH),\qquad
 F^\ell\mathfrak n_H
    =\bigoplus_{q\le-\ell}\Gamma(S,\End^qH).
\]
This is a finite multiplicative nilpotent filtration.
Its layers are finite projective over $C_S$, regardless
of the ranks of individual endomorphism sections.
Restrict the actual functor to the endomorphisms of the
scalar unit and use its just-established parameter-linear
matrix extension. It gives an $\ai$ map on global DGAs
\[
 F^H:B_S\otimes_{C_S}\mathfrak n_H
          \rightsquigarrow C\otimes_{C_S}\mathfrak n_H .
\]
Matrix composition gives the arity-additive filtration
bound of \eqref{eq:mc-arity-filtration}. In layer $\ell$
the unary map is
$I_S\otimes1_{\Gamma(S,\End^{-\ell}H)}$, a quasi-isomorphism.
Lemma~\ref{lem:finite-mc-lifting} applies to these global
complexes over $\R$; it does not require them to be
finite-dimensional real vector spaces.
It gives a global MC element $\mu$ and an actual finite
gauge $q$ with $b^q=F^H_*(\mu)$.
In particular the primitives are global smooth sections,
not unrelated choices on parameter fibers.

The source object
\[
 P=\bigl(\Gamma(S,H)\otimes_{C_S}B_S,Q+\mu\bigr)
\]
uses the basic scalar connection $Q$ and the forms-first
matrix action of $\mu$. It is bounded finite projective
and regular. Its internal
differential and one-form part vanish, because a total
degree-one element with negative internal degree has CE
degree at least two. The actual object formula is literally
\begin{equation}\label{eq:family-literal-object}
 \mathsf{IN}^{K_S}(P)=e_H+F^H_*(\mu).
\end{equation}
Indeed every MC input has total degree one, so the outer
unsuspension tensor sign is zero, as in the proof of
Theorem~\ref{thm:affine-rh}. The matrix sign remains
\eqref{eq:in-unsuspended-matrix-sign}, and
\eqref{eq:in-object-formula} needs no second target conversion.
All these object words are finite by negative internal
degree.

The gauge $q$ is a normalized closed degree-zero raw
morphism from \eqref{eq:family-literal-object} to
$e_H+b$, with direction fixed by
\eqref{eq:mc-gauge-orientation}. It is strictly invertible
by a finite series and has $q_0=1_H$; it can have higher
cochain components. With the preceding $\tau,u,i'$,
the actual counit is
\begin{equation}\label{eq:family-counit}
 \epsilon_R=i'\tau u q:
       \mathsf{IN}^{K_S}(P)\longrightarrow R.
\end{equation}
Here $\tau u$ is a vertex-only invertible map, whereas
$i'$ is a homotopy equivalence, not necessarily invertible.
This proves essential surjectivity on $H^0$.

\emph{Full MC.}
Lemma~\ref{lem:family-product} gives the actual raw
dg quasi-equivalence $\Psi_S^*$ and its unit-relative
invertible gauges
$U_Z:Z\to\Phi_S^*\Psi_S^*Z$.
For a reconstructed map
$v:\Psi_S^*\mathsf{IN}_{E_S}(P)\to\Psi_S^*R$, the map
\begin{equation}\label{eq:family-full-mc-counit}
 U_R^{-1}\Phi_S^*(v)U_{\mathsf{IN}_{E_S}(P)}
     :\mathsf{IN}_{E_S}(P)\longrightarrow R
\end{equation}
is an actual closed degree-zero homotopy equivalence.
The Hom comparison transfers through the same $\Psi_S^*$.
This proves the theorem for the intrinsic full-MC functor.
\end{proof}

\begin{cor}[Parameter pullback]\label{cor:family-parameter-naturality}
For every smooth map $f:S'\to S$, the actual family
integration functors commute with parameter pullback in
every Taylor arity, with the canonical identifications
compatible with identities and composition.
They give an equivalence of contravariant
$\Cat_\infty$-valued diagrams on the ordinary category
of smooth parameter manifolds and maps.
\end{cor}

\begin{proof}
Parameter pullback uses the scalar map
$C^\infty(S)\to C^\infty(S')$, $c\mapsto c\circ f$.
Finite-projective coefficients pull back through finite
idempotents. Regularity is preserved by pulling back the
actual smooth splitting of the internal complex; no
flatness of this ring map is required.
The product descriptions of $E_S,K_S,\Phi_S,\Psi_S$ and
the absence of parameter differentiation give, in every
Taylor arity, the actual identities
\[
 (f\times K)^*\mathsf{IN}^{K_S}
   \cong\mathsf{IN}^{K_{S'}}(f\times1_{\R^2})^*,
\]
and the corresponding identities for full MC integration.
For bundles, these use the canonical composition
identifications of pullback. Equivalently, in idempotent
models pullback is entrywise composition; the identities
then respect composition directly.
Passing to the dg/$\ai$ nerves gives the map of
contravariant diagrams. Theorem~\ref{thm:family-rh} makes
it pointwise an equivalence, hence an equivalence in the
functor infinity-category.
\end{proof}

The coherent inverse exists in the functor infinity-category.
Metrics, contractions, ODE gauges and MC primitives in the proof
are chosen objectwise; parameter naturality comes from the
forward integration functors. The parameter is passive: the
coefficient structure uses functions on $S$ and differentiates
only in the affine directions. A flat connection in the $S$
directions would be additional coefficient data. The indexing
category is the ordinary category of smooth parameter manifolds
and maps, before localization at parameter homotopies.

\subsection{Two-term bundle arrows and their examples}

\begin{thm}[Two-term family normal forms]
\label{thm:family-two-term}
Regular objects of internal cohomological amplitude $[0,1]$
in the three categories of Theorem~\ref{thm:family-rh} are
classified up to $H^0$-isomorphism by isomorphism classes of
smooth bundle maps
\[
 N:H^1\longrightarrow H^0
\]
over $S$, with $\Gamma(S,H^0)$ and $\Gamma(S,H^1)$ finite
projective. The rank of $N$ may vary.
The cohesive representative on $\pi^*(H^0\oplus H^1)$ is
\begin{equation}\label{eq:family-two-term-object}
 D_N(v_0)=0,\qquad D_N(v_1)=N(s)(v_1)\eta
\end{equation}
on basic sections, extended by the right Leibniz rule.
Denote this object by $P_N$. Its actual integrated
representative $R_N$ on $K_S$ is
\begin{equation}\label{eq:family-two-term-integrated}
 r_0=0,\qquad r_1=1_H,\qquad
 r_2=-h_{012}N(s),\qquad r_p=0\quad(p\ge3).
\end{equation}
\end{thm}

\begin{proof}
Finite minimalization removes acyclic tails and leaves
$\pi^*(H^0\oplus H^1)$. After the one-form gauge, degree
counting permits only a two-form component from $H^1$
to $H^0$. Write it as
$N(s,x,y)\eta+T(s,x,y)\theta_1\theta_2$.
Products of these negative-degree endomorphisms vanish.
Equation~\eqref{eq:affine-ce-even}, with parameters, gives
$T=N_x$: the difference is killed by both $x$ and $y$,
and continuity covers the origin.
Thus the remaining coefficient is $b_N=N\eta+
N_x\theta_1\theta_2$.

All its values at fixed $s$ lie in
$\Hom(H_s^1,H_s^0)$. The following integrals therefore
define global bundle-valued functions:
\[
 \begin{gathered}
 N_0(s)=N(s,0,0),\qquad
 P=-\int_0^1N_y(s,tx,ty)\,dt,\qquad
 R=\int_0^1N_x(s,tx,ty)\,dt,\\
 c=P\theta_1+R\theta_2,\qquad
 Qc=(N-N_0)\eta+N_x\theta_1\theta_2=b_N-N_0\eta .
 \end{gathered}
\]
The last equality follows from the fundamental theorem
of calculus and \eqref{eq:affine-ce-odd}.
The off-diagonal right-linear map $1+c$, with inverse
$1-c$, points from the unreduced $b_N$ object to $P_{N_0}$.
Indeed on a degree-one basic section,
$D_{N_0}(1+c)=(1+c)D_{b_N}$ is precisely the displayed
primitive identity. This proves existence globally.

For uniqueness, a closed degree-zero morphism between
normal forms has diagonal function maps $T_0,T_1$ and
an off-diagonal one-form $C=P\theta_1+R\theta_2$.
The closedness calculation
\eqref{eq:two-term-intertwining} is valid in these bundles.
At the zero section, the $\eta$-coefficient of $QC$ is
$-yP+xR=0$, so
\begin{equation}\label{eq:family-arrow-isomorphism}
 N'T_1(s,0,0)=T_0(s,0,0)N .
\end{equation}
An actual inverse modulo Hom differentials restricts,
together with its homotopies, to a vertex homotopy
inverse there. The internal differentials of the normal
forms are zero. Consequently the restricted $T_i$ are
actual smooth bundle isomorphisms.
Here evaluation acts on the inverse maps and their specified
homotopies.
Conversely any bundle isomorphisms intertwining $N,N'$
are basic closed degree-zero maps and give a strict
cohesive isomorphism. This proves the bundle-arrow
classification.

The actual functors and counits preserve vertex complexes
up to the specified homotopies, so
Theorem~\ref{thm:family-rh} restricts to the indicated
amplitude full subcategories.
For \eqref{eq:family-two-term-object}, every object word
with two occurrences of $N$ vanishes. The actual unary
kernel integrates $\eta$ to $-h$ on $K_S$, giving
\eqref{eq:family-two-term-integrated} and transferring
the same classification to both raw targets.
Regularity concerns the zero internal differential,
not the possibly varying rank of $N$.
\end{proof}

The invariant is a global bundle-arrow isomorphism class,
not a pointwise rank profile. This classifies objects up to
$H^0$-isomorphism; their Hom complexes remain those of the
respective dg categories.

\begin{example}[Nontrivial coefficient topology]
\label{ex:family-mobius}
On $S^1$, the periodic idempotent
\[
 p(\theta)=\frac12
 \begin{pmatrix}
  1+\cos\theta&\sin\theta\\
  \sin\theta&1-\cos\theta
 \end{pmatrix}
\]
has image the real M\"obius line $L$. It explicitly makes
$\Gamma(S^1,L)$ finite projective. A unit spanning vector
on $[0,2\pi]$ is
$v(\theta)=(\cos(\theta/2),\sin(\theta/2))$, with
$v(2\pi)=-v(0)$. A nowhere-zero section would require
a continuous nowhere-zero real coefficient
$a(\theta)$ with $a(2\pi)=-a(0)$, which is impossible.
Thus $L$ is not trivial.

The arrows $\id_L:L\to L$ and
$\id_\R:\underline\R\to\underline\R$ have the same
fiberwise dimension/rank triple $(1,1,1)$, but are not
isomorphic bundle arrows. Their cohesive and integrated
family coefficients are therefore inequivalent.
The nontriviality here is in the coefficient bundles;
the vertical affine resolution itself has free terms.
\end{example}

\begin{example}[Smooth vanishing order beyond rank]
\label{ex:family-vanishing-order}
For $S=\R$ and trivial coefficient lines, the residues
$N(t)=t$ and $N'(t)=t^3$ have the same rank at every
parameter. An isomorphism of bundle arrows would give
nowhere-zero smooth functions $a,b$ with $t^3b=at$.
For $t\ne0$ this forces $a=t^2b$, and continuity gives
$a(0)=0$, a contradiction.
The compact example $\sin\theta$ versus $\sin^3\theta$
on $S^1$ has the same local obstruction at a simple zero.
Both are regular coefficient families: changing residue
rank is allowed.
\end{example}

\subsection{The family ordinary-descent obstruction}

Let $\overline{\CG}_S$ be the product of the ordinary
quotient of Proposition~\ref{prop:ordinary-quotient} with
the unit groupoid on $S$. It is an ordinary groupoid in
smooth sheaves, with represented objects $\y M_S$.
The strict product map is
\[
 q_S:\y K_S\longrightarrow N(\overline{\CG}_S),
\]
and on full MC use $q_S\Phi_S$.

\begin{cor}[Exact family descent criterion]
\label{cor:family-ordinary-descent}
A regular two-term family lies in the essential image
on $H^0$ of $q_S^*$ exactly when its bundle-map residue
$N$ is identically zero. The same criterion holds for
$(q_S\Phi_S)^*$ on full MC.
If $N\ne0$, descent is impossible through every strict
map to an ordinary groupoid nerve with represented
objects, even when the bounded finite-projective target
coefficient is not regular or minimal.
\end{cor}

\begin{proof}
If $N=0$, the basic graded bundle over $S$ has identity
transport along every arrow of $\overline{\CG}_S$.
Its pullback is \eqref{eq:family-two-term-integrated}
with zero higher components.
The normal-form equivalences give the assertion for
any representative of that class.

If $N\ne0$, choose $s\in S$ with $N_s\ne0$.
Strict restriction to $\{s\}\times\y K$ preserves the
actual inverse maps and homotopies in any proposed
$H^0$-equivalence. A hypothetical pullback through any
strict ordinary-nerve map would therefore give the
forbidden point-parameter descent of
Theorem~\ref{thm:ordinary-nondescent}.
That theorem permits these ordinary targets, including
nonregular target coefficients.
For full MC, first pull back along the actual $\Psi_S$,
then restrict the parameter. The literal identity
$\Psi_S^*\mathsf{IN}_{E_S}(P_N)=R_N$ gives the same
contradiction.
\end{proof}

The quotient is used as a groupoid in smooth sheaves,
without a representability hypothesis on its arrow sheaf.
Comparison with AS or refined higher holonomy is a further
geometric question. These regular families vary the
parameter manifold and coefficient bundles while keeping
the affine algebroid fixed.

\section{Further directions}\label{sec:programme}

Theorem~\ref{thm:affine-open-rh} covers every ordinary open
$O\subseteq\R^2$, including nonrectangular, disconnected,
noncompact and origin-containing opens, on globally bounded
degreewise finite-projective real coefficient bodies whose
cohomology ranks may jump.
The whole base simplex, not merely its vertices, must remain
in $O$. The even generator $\eta$ retains every polynomial
power. The theorem compares the concrete cohesive and raw
categories, rather than the ordinary category $\Perf(B_O)$.

Actual differentiation and global smooth kernel acyclicity give
invertible identity-body counits, chosen objectwise. Internal
triangularity supplies a real-linear chain inverse that conjugates
the transferred and cohesive differentials. Choosing counits
coherently over coefficient moduli, and making the inverse sections
compatible with restriction between opens, are further
compatibility problems.

\paragraph{Parameters and effective descent.}
The regular vertical-family theorem of
Section~\ref{sec:affine-families} assumes regular body cohomology.
Its globally bounded finite-projective bundles need not be
trivial, and the parameter manifold need not be compact or
connected. The anchor is vertical, rather than that of an
inverse-image algebroid adding $TS$.
Corollary~\ref{cor:family-parameter-naturality} gives natural
forward functors and an equivalence of parameter diagrams.
The proof chooses metrics, contractions and MC primitives
objectwise, and the parameter $S$ remains passive.

Beyond fixed-product families, one can ask for a coefficient
comparison when the affine domain $O_s\subseteq\R^2$ varies
with $s\in S$, or when the parameter directions carry
connection or anchor data. Such extensions require a jointly
smooth global comparison compatible with the new geometry.
Effective descent asks a different
question: when do compatible local coefficients assemble into
globally bounded finite-projective objects? Such an argument
must control the amplitude and projectivity bounds together
with the higher gluing data.

\paragraph{Holonomy and reconstruction.}
Monodromy and derived-transverse refined
higher holonomy \cite{ZengHigherHolonomy} are distinct
constructions. The ordinary quotient used here is a groupoid in
smooth sheaves, obtained from parallel $2$-isomorphism classes.
Theorem~\ref{thm:ordinary-nondescent} concerns strict
maps to ordinary nerves in the actual raw coefficient theory.
Corollary~\ref{cor:family-ordinary-descent} retains its exact
zero-residue criterion for regular vertical families.
An intrinsic marked affine image and base change along surjective
submersions with connected fibres for closed saturated marks
are constructed in
\cite[Theorems~C--D]{ZengHigherHolonomy}.
Comparing our coefficient functor with these images requires
compatibility of the entire source maps, not only agreement on
the central degree-two class. An equivalence of coefficient
categories would additionally need a suitable descent theorem.

Example~\ref{ex:open-origin-observer} shows that coefficient RH
does not make origin observation reflect coefficient homotopy
equivalences, even on finite coefficients.
A reconstruction statement through a chosen observer must specify
whether it concerns the full coefficient category or a quotient
detected by that observer. A monoidal reconstruction would
additionally require compatible monoidal structures, and a
geometric reconstruction would require the corresponding geometric
interpretation. For derived-complete enlargements, the corresponding
question must also specify the completion and its coefficient
hypotheses.

A complementary infinitesimal perspective comes from CE filtrations.
For finite-length universal resolutions, the fibre-weight spectral
sequence is intrinsic from $E_1$
\cite[Theorem~3.2]{ZengIsotropySpectral2026}.
Finite-type minimal isotropy carries a separate pointwise sequence
\cite[Theorem~4.1]{ZengIsotropySpectral2026}; in the holomorphic cubic
Koszul example, a pointwise weight-one class has no unary-closed
stalk lift \cite[Proposition~6.2]{ZengIsotropySpectral2026}.
These filtrations suggest finer tests for what origin observation
retains from coefficient transport.

\paragraph{Convolution and the affine relation.}
The same affine action appears in the convolution-algebra
construction of \cite[Theorem~6.3]{FoliationFactorization2026}.
Its origin-germ AS relation module is the image of an injective
right-module operator
$D=-ay\,\partial_a+ax\,\partial_b-y$
in the total-density frame. The vector-field part differentiates
along the stabilizer relation, while the zeroth-order term comes
from its density action. A comparison with the coefficient residue
studied here would have to relate these distinct module categories
and account for this density twist.

\paragraph{Beyond the affine model.}
For general Lie- and $\li$-algebroids, one can ask when the
integration functor of \cite{ZengGlobalResolutions2026} is a
quasi-equivalence. Related questions concern local-resolution
effectivity and coherence at the level of coefficient categories.
A source-valued cohesive factorization would require an enriched
construction compatible with the integration and invariance
results used here. Further comparison problems concern
varying-base Morita maps and extensions to derived, Banach
or noncommutative geometries.

\appendix
\section{Global affine lifting and deformation retraction}
\label{app:affine-global-lifting}

Throughout this appendix, $M=\R^2$, $G=\Aff_+(1)$, and
$K$ is the crossed-module nerve of Section~\ref{sec:affine-integration}.
All differentials and forms are in the simplex direction,
with jointly smooth external parameters. We prove the lifting
and retraction results stated in that section; the extension
operators are also used in Section~\ref{sec:affine-open-geometry}.

\subsection{Holonomy and flat correction}

For an edge matrix
$\mathsf A(\partial_t)=
\left(\begin{smallmatrix}q_1(t)&q_2(t)\\0&0\end{smallmatrix}\right)$,
the solution has entries
\[
 r(t)=\exp\left(\int_0^tq_1(s)\,ds\right),\qquad
 s(t)=r(t)\int_0^tr(v)^{-1}q_2(v)\,dv.
\]
Thus it exists on the whole edge and stays in $G$.
The formulas preserve arbitrary joint smooth $U$-dependence:
the integrations have compact fixed fibers, and their
parameter derivatives are computed locally in $U$.
The base equation gives $\sigma(v_j)=g_{ij}\sigma(v_i)$.

Let $H_v$ be the radial Poincare operator on a simplex
based at a vertex $v$. It integrates the $ds$ component
of pullback by $(s,t)\mapsto v+s(t-v)$, so
\begin{equation}\label{eq:radial-homotopy}
 dH_v+H_vd=\id-v^*.
\end{equation}
For a closed 2-form, $u=H_v\beta$ satisfies $du=\beta$
and restricts to zero on all edges issuing from $v$.

\begin{lem}[Rank-uniform flat correction]\label{lem:affine-flat-correction}
Suppose $(\sigma,\mathsf A,\beta)$ satisfies
\eqref{eq:affine-full-mc}, and let $u$ be a one-form with $du=\beta$.
Then
\[
 \mathsf A'=\mathsf A-Z_\sigma u
\]
is flat and satisfies $d\sigma=\mathsf A'\sigma$.
A flat matrix form
$\mathsf A'=\left(\begin{smallmatrix}\xi&\zeta\\0&0\end{smallmatrix}\right)$
on a simplex has a unique potential $g$ with
$dg\,g^{-1}=\mathsf A'$ and any prescribed initial value
in $G$ at one vertex.
\end{lem}

\begin{proof}
The base equation is unchanged because $Z_\sigma\sigma=0$.
Expanding curvature, the terms $Z_\sigma\beta-Z_\sigma du$
cancel. The mixed terms cancel by
$dZ_\sigma=\mathsf AZ_\sigma-Z_\sigma\mathsf A$
and $u\wedge\mathsf A=-\mathsf A\wedge u$.
The remaining quadratic term vanishes because $u\wedge u=0$.

For the second assertion, flatness means
$d\xi=0$ and $d\zeta=\xi\wedge\zeta$.
Put
\[
 p=H_v\xi,\qquad q=H_v(e^{-p}\zeta),\qquad
 g_{\mathrm{norm}}=
 \begin{pmatrix}e^p&e^p q\\0&1\end{pmatrix}.
\]
Equation~\eqref{eq:radial-homotopy} gives
$dp=\xi$ and $dq=e^{-p}\zeta$. Hence
$dg_{\mathrm{norm}}g_{\mathrm{norm}}^{-1}=\mathsf A'$
and $g_{\mathrm{norm}}(v)=1$.
Right multiplication by the prescribed initial value
gives $g$. If $g_1,g_2$ are two such potentials with the same
initial value, $g_1^{-1}g_2$ has zero differential on the
connected simplex and is the identity.
All these operations preserve joint smooth parameters.
\end{proof}

For a potential based at $v_0$ with $g(v_0)=1$, the
base equation gives $\sigma=g m_0$:
$d(g^{-1}\sigma)=0$ and $m_0=\sigma(v_0)$.
Equivariance gives $Z_\sigma=gZ_{m_0}g^{-1}$.
Solving the edge equation in this gauge yields
\begin{equation}\label{eq:affine-gauged-holonomy}
 g_{ij}(\mathsf A)
 =g(v_j)\exp\left(Z_{m_0}\int_{[ij]}u\right)g(v_i)^{-1}.
\end{equation}
This can be verified directly by differentiating
$g(t)\exp(Z_{m_0}\int_0^tu)\,g(v_i)^{-1}$.

On a triangle based at $v_i$, take $u=H_{v_i}\beta$.
Its two root-edge restrictions vanish, so Stokes gives
$\int_{[jk]}u=\int_{[ijk]}\beta=-h_{ijk}$.
Equation~\eqref{eq:affine-gauged-holonomy} then gives
\[
 g_{jk}=g_{ik}\exp(-h_{ijk}Z_{\sigma(v_i)})g_{ij}^{-1},
\]
which is exactly \eqref{eq:affine-triangle-convention}.
Stokes on a tetrahedron gives $\delta h=0$.
Thus the actual edges and triangles form a $K$-simplex.
They commute with faces and degeneracies: affine
simplex maps carry an edge to the corresponding edge
or a constant edge, and a triangle to the corresponding
triangle or a map of rank at most one.
Degenerate coordinates are consequently the identity
and zero. This proves that $\Phi$ is a natural
simplicial sheaf map, independently of choices in
the subsequent lifting argument.

\subsection{Compatible boundary forms and prescribed periods}

Write $\Omega^r(\partial\Delta^n)$ for compatible
tangential forms on the codimension-one faces.
Only their pullbacks to intersections are identified;
normal derivatives are not specified.

\begin{lem}[Smooth face extension]\label{lem:smooth-face-extension}
For $n\ge1$ there is a linear extension operator
\[
 \operatorname{Ext}_n:\Omega^r(\partial\Delta^n)
 \longrightarrow\Omega^r(\Delta^n)
\]
right inverse to restriction, including for functions
when $r=0$. It preserves joint
smooth parameter families and commutes with
pullback in the external parameter.
\end{lem}

\begin{proof}
For $F_i=\{t_i=0\}$ use the retraction
\[
 (p_i(t))_i=0,\qquad
 (p_i(t))_j=\frac{t_j}{1-t_i}\quad(j\ne i)
\]
on $t_i<1$. Fix a smooth cutoff $\chi$ equal to
one near zero and zero on a neighborhood of
$[1/2,1]$. Extend a face form $\omega$ by
$\chi(t_i)p_i^*\omega$, with zero continuation
away from that collar. This is smooth through
all corners; the denominator does not degenerate
on the support.

Process the faces in order. On $F_i$ subtract the
restriction of the extension already made, and
extend this residual by the displayed operator.
Compatibility makes the residual vanish after
pullback to earlier face intersections.
The map $p_i$ sends $F_j$ into $F_i\cap F_j$
where defined, so the new term restricts to
zero on every earlier $F_j$.
After finitely many steps all face restrictions
are correct. The construction uses only fixed
smooth pullbacks, multiplication and addition,
so it has the stated parameter properties.
\end{proof}

For an ordered $p$-face $I=(i_0,\ldots,i_p)$ put
\[
 \omega_I=p!\sum_{r=0}^p(-1)^r t_{i_r}
   dt_{i_0}\wedge\cdots\wedge
   \widehat{dt_{i_r}}\wedge\cdots\wedge dt_{i_p}.
\]
Extend this notation alternatingly, with repeated
indices giving zero. Barycentric identities give
\begin{equation}\label{eq:whitney-identities}
 d\omega_I=\sum_j\omega_{jI},\qquad
 \int_J\omega_I=\delta_{IJ}.
\end{equation}
Indeed the first identity follows from
$\sum t_i=1$ and $\sum dt_i=0$; on its own
face $\omega_I$ is $p!$ times the oriented
coordinate volume form, and on another
$p$-face it vanishes. For a simplicial
cochain $c$, write $W(c)=\sum_I\omega_I c_I$.
Equation~\eqref{eq:whitney-identities} gives
$dW(c)=W(\delta c)$.

\begin{lem}[A smooth degree-two radial identity]
\label{lem:degree-two-radial}
For a closed smooth 2-form $\gamma$, put
\[
 \alpha_i=H_{v_i}\gamma,\qquad
 f_{ij}=H_{v_j}H_{v_i}\gamma,\qquad
 c_{ijk}=\int_{[ijk]}\gamma.
\]
The operator
\[
 D_n\gamma=\sum_i t_i\alpha_i
              -\sum_{i<j}\omega_{ij}f_{ij}
\]
satisfies
\[
 dD_n\gamma=\gamma-W(c).
\]
It commutes with face restrictions, vanishes
after restriction to any edge, and preserves
arbitrary joint smooth parameter dependence.
\end{lem}

\begin{proof}
The radial identity gives $d\alpha_i=\gamma$
and $df_{ij}=\alpha_i-\alpha_j$.
The function $f_{ij}$ vanishes at $v_i$ and
$v_j$: one endpoint is its radial base, and
$\alpha_i$ vanishes on the radial edge $[ij]$.
Moreover $d(f_{ij}+f_{ji})=0$, and their common
zero endpoint makes this constant zero.
Thus $f_{ji}=-f_{ij}$.
The function $f_{ij}+f_{jk}-f_{ik}$ is constant.
At $v_k$ it equals $\int_{[jk]}\alpha_i$,
which is $c_{ijk}$ by Stokes on $[ijk]$.

Differentiate the formula for $D_n$. Since
\[
 \sum_{i<j}\omega_{ij}\wedge(\alpha_i-\alpha_j)
 =\sum_{i,j}t_i\,dt_j\wedge(\alpha_i-\alpha_j)
 =-\sum_i dt_i\wedge\alpha_i,
\]
the mixed terms cancel, leaving
\[
 dD_n\gamma
 =\gamma-2\sum_{i<j}dt_i\wedge dt_j\,f_{ij}.
\]
Pointwise in the coefficient functions,
\eqref{eq:whitney-identities} gives
\[
 2\sum_{i<j}dt_i\wedge dt_j\,f_{ij}
 =\sum_{i<j<k}\omega_{ijk}
                (f_{ij}+f_{jk}-f_{ik})
 =W(c).
\]
This proves the identity for smooth forms
directly, not by polynomial approximation.

On a face, terms involving an absent vertex
vanish because its barycentric coordinate
and differential restrict to zero.
All remaining radial operators stay within
the face. This proves naturality.
A 2-form restricts to zero on an edge,
so its $D$-primitive does too.
Finite compact-interval integrations preserve
the stated joint smoothness.
\end{proof}

\begin{lem}[Closed extension with fixed triangle integrals]
\label{lem:closed-period-extension}
Let $n\ge2$, let $c$ be a smooth parameter family of
simplicial 2-cocycles on $\Delta^n$, and let
$\beta_\partial$ be compatible closed boundary
2-forms with integrals $c$ on the boundary
triangles. There is a jointly smooth closed
2-form $\beta$ on the simplex with this
boundary and all triangle integrals equal to $c$.
\end{lem}

\begin{proof}
On each face apply Lemma~\ref{lem:degree-two-radial}
to $\gamma_\partial=\beta_\partial-W(c)|_\partial$.
Its triangle integrals vanish, so
\[
 \eta_\partial=D\gamma_\partial,\qquad
 d\eta_\partial=\gamma_\partial.
\]
The resulting 1-forms agree on intersections
and vanish on every edge. Set
\[
 \beta=W(c)+d\operatorname{Ext}_n(\eta_\partial).
\]
It is closed because $\delta c=0$, and its
restriction to each face is $\beta_\partial$.
Every edge lies in the boundary when $n\ge2$,
so Stokes makes the exact correction have zero
integral on every triangle.
When $n=2$, the boundary 2-forms and
$\eta_\partial$ are zero, and $\beta=W(c)$
has the prescribed unrestricted top-triangle
integral. For $n=3$, $\delta c=0$ is exactly
the tetrahedral compatibility. All operations
have the parameter properties already proved.
\end{proof}

\subsection{All relative boundary lifts}

\begin{proof}[Proof of Theorem~\ref{thm:affine-integration}]
Fix $U$ and a relative boundary problem over
$k\in C^\infty(U,K_n)$.
Thus compatible full MC fields are prescribed
on all faces, and their images are the faces
of $k$. In degree zero $\Phi_0$ is the identity.
In degree one the target specifies $m_0,a_1$
and $m_1=a_1m_0$. Using
$\psi(a,b)=(\log a,b)$, take
\[
 g(t)=\psi^{-1}(t\psi(a_1)),\qquad
 \sigma=g m_0,\qquad
 \mathsf A=dg\,g^{-1},\qquad\beta=0.
\]
This has the required endpoints and holonomy.

Suppose $n\ge2$. Set $c_{ijk}=-h_{ijk}$ from
the target. Its tetrahedral equations say
$\delta c=0$, and boundary compatibility gives
the required integrals of $\beta_\partial$.
Lemma~\ref{lem:closed-period-extension} supplies
$\beta$. Put $u=H_{v_0}\beta$; then
$du=\beta$ and $u|_{[0i]}=0$.

On a face $F_i$, the prescribed fields give
the flat connection
\[
 \mathsf A_i'=\mathsf A_i-Z_{\sigma_i}u|_{F_i}
\]
with the same base equation. These connections
agree on intersections.
For $i\ne0$, use the potential $g_i$ based
at $v_0$ with value $1$, from
Lemma~\ref{lem:affine-flat-correction}.
Then $\sigma_i=g_i m_0$, and
$g_i(v_r)=a_r$ whenever that vertex lies
in $F_i$, since the root-edge correction is zero.
The potentials agree on intersections containing
$v_0$.

On the opposite face $F_0$, take the potential
$g_0$ with initial value $a_1$ at $v_1$.
It satisfies $\sigma_0=g_0m_0$.
It agrees with every $g_i$, $i\ne0,1$, on
their intersection, by their common value at
$v_1$ and uniqueness of a flat potential.
For $F_0\cap F_1$, which contains $v_2$
even when $n=2$, the holonomy formula gives
\[
 g_{12}=g_0(v_2)
       \exp\left(Z_{m_0}\int_{[12]}u\right)a_1^{-1}.
\]
Stokes and the zero root-edge restrictions give
$\int_{[12]}u=-h_{012}$. Hence
\[
 g_0(v_2)=g_{12}a_1\exp(h_{012}Z_{m_0})=a_2
\]
by the target triangle convention.
Uniqueness on the connected intersection
now gives $g_0=g_1$ there.
This constructs a compatible jointly smooth
boundary potential $g_\partial$.

Extend $\psi(g_\partial)$ componentwise by
Lemma~\ref{lem:smooth-face-extension}, and
apply $\psi^{-1}$ to obtain $g$ on the full
simplex. This step uses the global chart
$G\cong\R^2$, not an extension assertion for
an arbitrary Lie group.
Define
\[
 \sigma=g m_0,\qquad
 \mathsf A'=dg\,g^{-1},\qquad
 \mathsf A=\mathsf A'+Z_\sigma u.
\]
The form $\mathsf A'$ is flat,
$d\sigma=\mathsf A'\sigma$, and equivariance
gives $dZ_\sigma=\mathsf A'Z_\sigma-Z_\sigma\mathsf A'$.
Since $Z_\sigma\sigma=0$, the reverse
flat-correction calculation yields
\[
 d\sigma=\mathsf A\sigma,\qquad
 d\mathsf A-\mathsf A\wedge\mathsf A=Z_\sigma\beta.
\]
Together with $d\beta=0$, these are all of
\eqref{eq:affine-full-mc}.
On each face the potential, base map,
matrix form and 2-form recover the prescribed
data exactly.

The root holonomies are $a_i$, since $u$
vanishes on the root edges and $g(v_0)=1$.
All triangle integrals are $-h_{ijk}$.
The root coordinates
\eqref{eq:affine-root-coordinates} therefore
show that the image is exactly $k$.
We have solved every relative boundary
problem, naturally in the external parameter.
The right lifting property against
$\partial\Delta[n]\hookrightarrow\Delta[n]$
for all $n\ge0$ is precisely the trivial
Kan fibration criterion in simplicial sets.
The same smooth formulas apply across the rank-drop locus.
\end{proof}

\subsection{Monodromy and the fixed-object source}

\begin{cor}\label{cor:affine-monodromy}
The comparison gives an equivalence under $\y M$
\[
 \operatorname{Mon}^{\mathrm{sh}}(E_{\mathrm{aff}})
 \simeq|\y K|
\]
and the corresponding hypercomplete equivalence.
These affine monodromy objects are 2-truncated
smooth infinity-sheaves and hence hypercomplete.
Their associated derived fixed-object Cech
groupoid objects are equivalent.
\end{cor}

\begin{proof}
Objectwise weak equivalence gives equivalence
of presheaves of spaces and therefore of their
open-cover sheafifications. Vertex preservation
makes this an equivalence under $\y M$.
For every $U$, $K(U)$ is the nerve of a strict
2-groupoid. Its realization is a homotopy
2-type, so the presheaf of spaces is
2-truncated. Left-exact sheafification
preserves that property, and truncated
infinity-sheaves are hypercomplete.
Finally, the Cech levels are finite derived
fiber products over these equivalent
realizations with the common object $\y M$.
Their comparison is compatible with every
simplicial operator.
\end{proof}

The derived Cech groupoid in the corollary is not the raw
Duskin nerve $K$. The result compares monodromy sources. A refined
transverse-holonomy comparison would also have to specify
its target, slice and transport.
Likewise, 2-truncation here is a statement
in the smooth infinity-topos, not about
the ordinary topological shape of a
smooth stack.

\begin{cor}\label{cor:affine-all-global-resolutions}
Every other global geometric resolution of
this same affine foliation satisfying the
hypotheses of Section~\ref{sec:polynomial}
has equivalent full smooth integration and
the same finite monodromy model $K$.
In either smooth infinity-target these
comparisons form a coherent diagram,
with a contractible space of choices
after framing at the specified affine
comparison.
\end{cor}

\begin{proof}
Apply Theorem~\ref{thm:coherent-resolutions}
to $\mathfrak f=\CF(M)$ and use
Corollary~\ref{cor:coherent-monodromy},
with $E_{\mathrm{aff}}$ as the selected
resolution. Compose its constant value
with the equivalence induced by $\Phi$.
At the presheaf level, the contractible
Kan comparison spaces supply inverse
vertices and polynomial Hom paths for
each actual comparison.
Theorem~\ref{thm:enriched-integration}
turns them into simplicial homotopies.
Thus these maps are also objectwise
homotopy equivalences, independently
of the general derived-source
identification.
\end{proof}

\subsection{A global simplicial deformation retraction}

An objectwise trivial Kan fibration does not by
itself supply smooth simplicial inverse data.
We now construct such data for the actual map
$\Phi$. Local smooth simplicial deformation
retractions are already available in
\cite[Theorems~8.1--8.2]{SeveraSiran2020}.
The assertion here concerns the global affine
comparison and its specified crossed-module target.

For fixed $n$, jointly smooth simplex-direction
forms constitute a sheaf of real vector spaces
\[
 U\longmapsto\Omega^q(U\times\Delta^n/U).
\]
Data parametrized by an arbitrary smooth sheaf
$P$ mean morphisms from $P$ to this vector sheaf,
defined on plots. The sheaf $P$ may be nonrepresentable;
the data are sheaf morphisms rather than elements of a
function algebra on $P$.

\begin{lem}[Ordered split-locus extension]
\label{lem:split-locus-extension}
Let $P_1,\ldots,P_N\subseteq P$ be split
subsheaves with retractions $r_j:P\to P_j$.
Suppose
\[
 r_j(P_i)\subseteq P_i\cap P_j\qquad(i<j).
\]
Compatible morphisms $v_j:P_j\to\mathcal V$
to a sheaf of real vector spaces extend to
$v:P\to\mathcal V$. The same holds inside
any specified linear subsheaf of $\mathcal V$.
\end{lem}

\begin{proof}
Starting with $\widetilde v_0=0$, set
\[
 \widetilde v_j=\widetilde v_{j-1}
  +\bigl(v_j-\widetilde v_{j-1}|_{P_j}\bigr)r_j.
\]
The residual in parentheses vanishes on each
earlier intersection $P_i\cap P_j$. Its pullback
therefore vanishes on $P_i$, while $r_j|_{P_j}=\id$
gives the desired new restriction. All operations
preserve a linear subsheaf. This also covers
$N=0$, and does not require commuting retractions.
\end{proof}

\begin{lem}[Affine lifting with split parameter loci]
\label{lem:affine-split-lifting}
Put $S=\Int(E_{\mathrm{aff}})$ and $T=\y K$.
Let $k:P\to T_n$ have compatible prescribed
faces $b_i:P\to S_{n-1}$ over its faces, when
$n>0$. Suppose lifts $\ell_j:P_j\to S_n$ are
given on finitely many split loci as in
Lemma~\ref{lem:split-locus-extension}, agreeing
with these faces and with one another on
intersections. Then there is a lift
$\ell:P\to S_n$ with
\[
 \Phi_n\ell=k,\qquad d_i\ell=b_i,\qquad
 \ell|_{P_j}=\ell_j.
\]
The retractions need not preserve $k$.
\end{lem}

\begin{proof}
The face condition means
$d_a b_i=d_{i-1}b_a$ for $a<i$ and
$\Phi_{n-1}b_i=d_i k$. Thus all data are
morphisms of sheaves of the indicated degrees.
We extend only vector-valued residuals, using
baselines constructed from the current target $k$.

For $n\ge2$, put $c=-h(k)$ and let
$\beta_\partial$ be the prescribed boundary
2-forms. Lemma~\ref{lem:closed-period-extension}
gives the baseline
\[
 \beta^0=W(c)+d\operatorname{Ext}_n
       \bigl(D_\partial(\beta_\partial-W(c)|_\partial)\bigr).
\]
It has the prescribed boundary and triangle
periods. If $\beta_j$ is the 2-form of $\ell_j$,
then
\[
 \gamma_j=\beta_j-\beta^0|_{P_j},\qquad
 \eta_j=D_n\gamma_j
\]
satisfy $d\eta_j=\gamma_j$: the closed form
$\gamma_j$ has zero triangle periods.
Face naturality of $D_n$ makes $\eta_j$
vanish on the whole boundary, since
$\gamma_j|_\partial=0$. The $\eta_j$ agree
on parameter-locus intersections.
Apply Lemma~\ref{lem:split-locus-extension}
in the linear sheaf of 1-forms vanishing
on $\partial\Delta^n$, and set
\[
 \beta=\beta^0+d\eta.
\]
This recovers every $\beta_j$ and the full
boundary. Every edge lies in the boundary
for $n\ge2$, so Stokes shows that the
correction has zero triangle periods.
In particular the current periods $c(k)$
have not been pulled back along any $r_j$.

Put $u=H_{v_0}\beta$. Build the compatible
boundary potential $g_\partial$ exactly as
in the preceding all-boundary proof, using
this $u$ and the current root values $a_i(k)$.
The opposite face uses the potential with
value $a_1(k)$ at $v_1$; its compatibility
is the same plus-$h$ triangle calculation.
Let $g^0$ be the extension of $g_\partial$
through the global chart
$\psi(a,b)=(\log a,b)$.

On $P_j$, flatten the prescribed simplex
with $u_j=H_{v_0}\beta_j=u|_{P_j}$, and
let $g_j$ be its potential normalized by
$g_j(v_0)=1$. It agrees with $g_\partial$
on all faces: on a face containing $v_0$
the flat equation and initial value agree;
on the opposite face its value at $v_1$
is $a_1(k)$, because $u_j$ vanishes on
the root edge. Uniqueness of the flat
potential gives the required agreement.
Consequently
\[
 w_j=\psi(g_j)-\psi(g^0)|_{P_j}
\]
is a compatible family of $\R^2$-valued
functions vanishing on the simplex boundary.
Extend it within that linear subsheaf and put
\[
 g=\psi^{-1}\bigl(\psi(g^0)+w\bigr),\qquad
 \sigma=g m_0(k),\qquad
 \mathsf A=dg\,g^{-1}+Z_\sigma u.
\]
There is no range constraint on the chart-valued
sum because $\psi$ is global.

For completeness, the modified lift still
solves the entire problem. The form
$\mathsf A'=dg\,g^{-1}$ is flat,
$d\sigma=\mathsf A'\sigma$, and
$dZ_\sigma=[\mathsf A',Z_\sigma]$.
With $du=\beta$, $d\beta=0$ and
$Z_\sigma\sigma=0$, reverse flat correction
gives all of \eqref{eq:affine-full-mc}.
The forms and potential recover every
prescribed face. Since $w|_\partial=0$,
$g(v_0)=1$ and $g(v_i)=a_i(k)$.
The root-edge restrictions of $u$ vanish,
so \eqref{eq:affine-gauged-holonomy} gives
the current root holonomies $a_i(k)$.
The triangle periods are $-h_{ijk}(k)$.
The root coordinates
\eqref{eq:affine-root-coordinates} therefore
give $\Phi_n\ell=k$ exactly.
On $P_j$ we have $\beta=\beta_j$,
$u=u_j$, $g=g_j$, and
$\sigma=g_jm_0(k|_{P_j})=\sigma_j$.
The matrix form is then $\mathsf A_j$,
so the whole prescribed simplex is recovered.

For $n=1$, $\beta=u=0$. Take
\[
 g^0(t)=\psi^{-1}\bigl(t\psi(a_1(k))\bigr).
\]
Here $\psi(1)=0$, so the endpoints are
$1$ and $a_1(k)$. Each prescribed simplex
has a normalized flat potential $g_j$.
Extend $\psi(g_j)-\psi(g^0)|_{P_j}$
in the vector sheaf of functions vanishing
at both endpoints, and reconstruct as above.
For $n=0$ the lift is forced by $\Phi_0=\id$.

All operations are natural on external plots:
fixed pullbacks, finite linear operations,
compact-fiber integrations, and the displayed
smooth group formulas. These formulas are defined
uniformly across the rank-drop locus.
\end{proof}

\begin{lem}[Ordered degeneracy loci]
\label{lem:ordered-degeneracy-loci}
For a simplicial sheaf $Z$, put
$B_i=\im(s_i:Z_{n-1}\to Z_n)$ and
$r_i=s_i d_{i+1}$. These are split loci
satisfying $r_i(B_j)\subseteq B_i\cap B_j$
for $j<i$. If $i<j$ and $\sigma\in B_i\cap B_j$,
then
\[
 \sigma=s_i s_{j-1}z=s_j s_i z,\qquad
 z=d_i d_j\sigma.
\]
\end{lem}

\begin{proof}
The identity $d_{i+1}s_i=\id$ gives the splitting.
For $j<i$,
\[
 r_i s_j=s_i d_{i+1}s_j
        =s_i s_jd_i=s_j s_{i-1}d_i.
\]
This proves preservation of the earlier locus.
For the intersection formula write
$\sigma=s_j b$ with $b=d_j\sigma$.
Then $d_i\sigma=s_{j-1}d_i b$, and
$\sigma=s_i d_i\sigma$ because $\sigma\in B_i$.
The two displayed double degeneracies agree
by the simplicial identity.
\end{proof}

\begin{thm}[Global affine simplicial retraction]
\label{thm:affine-sdr}
The comparison $\Phi:S\to T$ admits a smooth
simplicial section $\Psi:T\to S$ and a
simplicial homotopy
\[
 H:S\times\Delta[1]\longrightarrow S,\qquad
 H(-,0)=\Psi\Phi,\quad H(-,1)=\id_S.
\]
Moreover $\Phi H=\Phi\pr_S$ and $H$ is fixed
on $\Psi(T)$. In particular it is relative
to $c(\y M)$. Thus the specified affine map
is a global smooth simplicial strong
deformation retraction.
\end{thm}

\begin{proof}
First construct the section by induction.
In degree zero set $\Psi_0=\id_M$.
At degree $n>0$, take parameter sheaf $T_n$,
current target $k=\id_{T_n}$, and boundary
$b_a=\Psi_{n-1}d_a$. On
$B_i=\im(s_i)\subset T_n$ prescribe
\[
 \ell_i(k)=s_i\Psi_{n-1}(d_{i+1}k).
\]
The face identities
$d_a s_i=s_{i-1}d_a$ for $a<i$,
$d_i s_i=d_{i+1}s_i=\id$, and
$d_a s_i=s_i d_{a-1}$ for $a>i+1$
show that these prescriptions have boundary
$b_a|_{B_i}$. Their target is $k|_{B_i}$
by the lower-degree section identities.
On $B_i\cap B_j$, $i<j$, use
Lemma~\ref{lem:ordered-degeneracy-loci}:
the two values are
\[
 s_i s_{j-1}\Psi_{n-2}(z)
       =s_j s_i\Psi_{n-2}(z).
\]
The same lemma orders the split loci.
Lemma~\ref{lem:affine-split-lifting} gives
$\Psi_n$, with all faces and degeneracies
prescribed. This constructs a simplicial
section, not just levelwise sections.

Now construct $H_n(-,\tau)$ inductively for
each label $\tau\in\Delta[1]_n$.
Constant labels give the two stated endpoint
maps. For a nonconstant label take $P=S_n$,
$k=\Phi_n$, and boundary
\[
 b_a(\sigma)=H_{n-1}(d_a\sigma,d_a\tau).
\]
The lower homotopy identities make this
boundary compatible and over $k$.
First impose the fixed locus
$B_*=\Psi(T_n)$, with retraction $\Psi\Phi$
and prescribed value $\sigma$ there.
For every $i$ with $\tau(i)=\tau(i+1)$,
also impose $B_i=\im(s_i)$ and prescribe
\[
 \ell_i(\sigma)
   =s_i H_{n-1}(d_{i+1}\sigma,d_{i+1}\tau).
\]
The three face--degeneracy identities above,
applied to both $\sigma$ and its label,
give the stated boundary. The target is
$s_i d_{i+1}k=k$ on $B_i$.

Here is the intersection check, including
the labels. Write $\varsigma_i$ for the
order-preserving codegeneracy and, on
$B_i\cap B_j$ with $i<j$, write
$\sigma=s_i s_{j-1}z=s_j s_i z$.
The equal adjacent values of $\tau$ let
it factor through the common surjection
\[
 \varsigma_{j-1}\varsigma_i
        =\varsigma_i\varsigma_j.
\]
The two reduced labels coincide: their
pullbacks along this surjection are $\tau$,
and the surjection is onto. Calling this
label $\widehat\tau$, the two prescribed
values are
\[
 s_i s_{j-1}H_{n-2}(z,\widehat\tau)
       =s_j s_iH_{n-2}(z,\widehat\tau).
\]
On $B_*\cap B_i$, if $\Psi(k)=s_i\sigma'$,
then $k=s_i\Phi(\sigma')$ and simpliciality
of $\Psi$, followed by the splitting of
$s_i$, gives $\sigma'=\Psi\Phi(\sigma')$.
The lower homotopy fixes this lower
simplex, so both prescriptions return
$\Psi(k)$. Boundary compatibility on
$B_*$ follows from the same lower
fixed-locus condition.

Every $r_i=s_i d_{i+1}$ preserves $B_*$
because $\Psi$ and $\Phi$ are simplicial.
It also preserves all earlier degeneracy
loci. Thus order $B_*$ first, followed by
the relevant $B_i$ increasingly, and apply
Lemma~\ref{lem:affine-split-lifting}.
The resulting $H_n$ has all required
faces and degeneracies of the product.
Labels on a test manifold are locally
constant; the constructions on constant
labels glue, including on disconnected
test manifolds.

In degree one the sole nonconstant label
is $(0,1)$, so only the fixed locus is imposed.
In degree two the nonconstant labels
$(0,0,1)$ and $(0,1,1)$ impose $B_0$ and
$B_1$, respectively. Their face conditions
are precisely the degree-one conditions
already constructed; the opposite-face
calculation is included in the lifting lemma.
All formulas include the origin and $y=0$.
Finally, $\Psi$ preserves degeneracies and
is the identity in degree zero, so
$c(\y M)\subseteq\Psi(T)$ and $H$ fixes it.
\end{proof}

Section~\ref{sec:representation-invariance}
extracts the raw coefficient consequence
from these actual inverse data, not from
weak equivalence alone. This particular
whole-plane section and homotopy can move
the interior base simplex and cannot simply
be restricted to arbitrary opens.
The construction of Section~\ref{sec:affine-open-geometry} instead
uses range-preserving charts on admissible
arrow domains, chosen separately for each open.
The global chart $G\cong\R^2$ is essential to this proof.
For other simply connected Lie groups, an analogous construction
would require a replacement for this chart-based extension.

\section{Regular reductions and scalar parameters}
\label{sec:affine-hom}

This appendix proves the regular-reduction and scalar-parameter
results used by the normal forms and the vertical-family theorem.
It also gives the regular Hom comparison by that route. The
nonregular proof in Sections~\ref{sec:affine-relative-comparison}
and~\ref{sec:affine-nonregular-rh} does not use the cohomology
splitting made here. Throughout this appendix,
$E=E_{\mathrm{aff}}$, $B=B_E$, $G=\Aff_+(1)$, and $K$ is the
crossed-module nerve of Section~\ref{sec:affine-integration}. Set
\[
 \mathcal R_K=\Rep_{\mathrm{sm}}^{\mathrm{perf}}(\y K),
 \qquad
 \mathsf{IN}^K=\Psi^*\mathsf{IN}_E.
\]
The base is the whole $\R^2$, including the origin. The section
$\Psi$ is the whole-plane choice; ordinary opens use the
range-preserving sections of Section~\ref{sec:affine-open-geometry}.

\begin{thm}[Global affine Hom comparison]
\label{thm:affine-hom}
\label{proposal:affine-ve-hom}
For every $P,Q\in\Coh^{\mathrm{reg}}(B)$, the actual first Taylor map
\[
 \mathsf{IN}^K_1:
 \Hom^\bullet_{\Coh(B)}(P,Q)
 \longrightarrow
 \Hom^\bullet_{\mathcal R_K}
       (\mathsf{IN}^K P,\mathsf{IN}^K Q)
\]
is a quasi-isomorphism in every integer degree. The same Hom
assertion holds for $\mathsf{IN}_E$ on the full MC target.
\end{thm}

The proof has two parts. A finite regular
coefficient reduction makes a scalar comparison sufficient.
An explicit relative de Rham calculation then proves that scalar
comparison for the actual nerve. The coefficient reduction
is what makes the scalar comparison determine all Hom degrees.

\subsection{Finite regular minimalization}

The following elementary perturbation statement also applies to an
associative coefficient DGA whose degree-zero algebra is not central.
This distinction will be useful for normalized simplicial cochains.

\begin{lem}[Finite superconnection perturbation]
\label{lem:affine-finite-minimalization}
Let $\mathcal D$ be a nonnegative cochain DGA with
$\mathcal D^0=A=C^\infty(M)$, and use the right-module convention of
Section~\ref{sec:cohesive}. Let
$(V\otimes_A\mathcal D,D)$ have a flat right superconnection $D$,
bounded finite-projective coefficient module $V$, and regular
internal differential $d_0=D^{[0]}$.
It is cohesively homotopy equivalent to a module
$(H\otimes_A\mathcal D,D_H)$ with $D_H^{[0]}=0$, where
$H=H^\bullet(V,d_0)$ is a bounded finite-projective graded module.
The displayed formulas are finite and define smooth maps and
homotopies on the resulting globally bounded finite-projective
bundles. If the differential and the chosen contraction data
vary jointly smoothly with auxiliary finite-dimensional parameters
and with a fixed global degree bound, the same is true of the
transferred formulas. For a smaller coefficient class, this
conclusion requires preservation by the chosen contraction maps
and by all operations used in the formulas.
\end{lem}

\begin{proof}
Boundedness and locally constant cohomology ranks imply locally
constant ranks of the differentials, by induction from either end
of the finite complex. Thus kernels and images are smooth
subbundles. A smooth fiber metric gives a decomposition
$V^k=B^k\oplus H^k\oplus C^k$, where $B^k=\operatorname{im}d_0$
and $d_0:C^k\to B^{k+1}$ is an isomorphism.
Take $i,p$ to be inclusion and projection of $H$, and let $h$
be this inverse on $B$, zero on $H\oplus C$. This gives
\[
 (H,0)\mathrel{\mathop{\rightleftarrows}^{i}_{p}}(V,d_0),
 \qquad
 pi=1,\quad d_0h+hd_0=1-ip,\quad h^2=hi=ph=0.
\]
In particular, $H$ is a direct summand of $V$ and is finite
projective. Extend $i,p,h$ right $\mathcal D$-linearly.

Put $\Delta=D-d_0$. A component of positive $\mathcal D$-degree $r$
changes internal coefficient degree by $1-r$. Consequently both
$h\Delta$ and $\Delta h$ strictly lower internal degree. This also
holds for the scalar-differential term in the Leibniz rule.
Their inverse series therefore terminate. Define
\begin{align}
 \Sigma&=\Delta(1+h\Delta)^{-1}
             =(1+\Delta h)^{-1}\Delta,
 \nonumber\\
 i'&=i-h\Sigma i=(1+h\Delta)^{-1}i,
 &p'&=p-p\Sigma h=p(1+\Delta h)^{-1},
 \label{eq:affine-finite-transfer}\\
 h'&=h-h\Sigma h=h(1+\Delta h)^{-1},
 &D_H&=p\Sigma i=p\Delta i'.
 \nonumber
\end{align}
The perturbation $\Delta$ is not assumed $\mathcal D$-linear.

To check right-linearity, observe that the maps
$(h\Delta)^n i$ are right-linear and are killed by $h$:
at the induction step the extra Leibniz term is
$(-1)^{|z|}hz\,d_{\mathcal D}b=0$.
For the other two series, put
\[
 f_n=p(\Delta h)^n,\qquad k_n=h(\Delta h)^n.
\]
Inductively each is right-linear and
$f_nh=k_nh=0$. Indeed, if $f_n$ has these properties, then
\[
 f_n\Delta h(zb)
 =f_n\Delta h(z)b+
       (-1)^{|z|-1}f_nh(z)d_{\mathcal D}b
 =f_n\Delta h(z)b,
\]
and $f_n\Delta h\,h=0$; the argument for $k_n$ is identical.
Thus $i',p',h'$ are right-linear. This proof uses neither
commutativity nor centrality of $\mathcal D^0$.

Since $pi'=pi=1$,
\[
 D_H(zb)=D_H(z)b+(-1)^{|z|}z\,d_{\mathcal D}b.
\]
For completeness, the finite operator identities give
\begin{gather}
 Di'=i'D_H,\qquad p'D=D_Hp',\qquad p'i'=1,
 \label{eq:affine-transfer-identities}\\
 Dh'+h'D=1-i'p',\qquad D_H^2=0.
 \nonumber
\end{gather}
For example, $i'=i-h\Delta i'$ and
$d_0\Delta=-\Delta D$ imply
$(1+h\Delta)Di'=iD_H$.
Similarly $p'=p-p'\Delta h$ and $\Delta d_0=-D\Delta$ give
$p'D(1+\Delta h)=D_Hp$; here
$D_H=p\Delta i'=p'\Delta i$.
These identities prove the two chain-map assertions.
The side conditions give $p'i'=1$ directly.
To check the homotopy identity, multiply it on the left by
$1+h\Delta$ and on the right by $1+\Delta h$.
Both sides become
\[
 1-ip+\Delta h+h\Delta+h\Delta^2h,
\]
using $D^2=d_0^2=0$. Finally $Di'=i'D_H$ and $p'i'=1$
give $D_H^2=0$. Every term of $D_H$ has positive
$\mathcal D$-degree. All expressions are finite in the global
coefficient amplitude. The chosen contraction data and the
differential are smooth, so the transferred maps, homotopies
and superconnection are smooth on the stated finite-projective
bundles.
\end{proof}

\subsection{The special affine one-form gauge}

Recall the actual CE differential:
\begin{align}
 B&=C^\infty(\R^2)[\theta_1,\theta_2,\eta],
 &|\theta_i|&=1,\quad |\eta|=2,
 \nonumber\\
 Qf&=f_x(x\theta_1+y\theta_2),
 &Q\theta_1&=-y\eta,
 \label{eq:affine-hom-ce}\\
 Q\theta_2&=\theta_1\theta_2+x\eta,
 &Q\eta&=0.
 \nonumber
\end{align}
Write $\xi=x\theta_1+y\theta_2=Qx$.

\begin{lem}\label{lem:affine-minimal-gauge}
Every regular cohesive $B$-module is cohesively homotopy equivalent
to one whose internal differential and CE-degree-one connection
term both vanish.
\end{lem}

\begin{proof}
Apply Lemma~\ref{lem:affine-finite-minimalization}. The resulting
graded cohomology bundle is globally trivial on $\R^2$.
In a graded frame, write its minimal superconnection as $Q+A$,
with
\[
 A_1=\theta_1L_1+\theta_2L_2.
\]
The $\eta$ coefficient of the CE-degree-two flatness equation is
$-yL_1+xL_2=0$. Higher connection terms do not contribute to
this equation, since the internal differential is zero.

The relation gives $L_1(0,y)=0$ for $y\ne0$, and then at $y=0$
by continuity. The global Hadamard formula
\[
 L(x,y)=\int_0^1\partial_xL_1(tx,y)\,dt
\]
gives $L_1=xL$. The same relation and continuity give $L_2=yL$.
Thus $A_1=\xi L$ on the whole plane.

In the forms-first convention, a degree-zero frame change gives
$A^g=g^{-1}Ag+g^{-1}Qg$. Solve
\[
 g_x=-Lg,\qquad g(0,y)=1.
\]
This linear matrix ODE has an invertible solution for every finite
$x$, smoothly in $y$. The usual local estimates on compact
parameter sets give joint smoothness; no compactness of $\R^2$
is assumed. Since $Qg=\xi g_x$, the transformed one-form part
vanishes. All remaining connection components strictly lower
internal degree.
\end{proof}

\begin{lem}[Reduction of the actual Hom map]
\label{lem:affine-hom-reduction}
Let
\[
 I:(B,Q)\longrightarrow
 C_K^\bullet:=C^\bullet_{\mathrm{norm,sm}}(K)
\]
be the scalar first Taylor map, and write $\delta$ for the full
simplicial coboundary on $C_K^\bullet$ throughout this appendix.
If $I$ is a quasi-isomorphism, then the first assertion of
Theorem~\ref{thm:affine-hom} holds.
\end{lem}

\begin{proof}
First replace the two objects using
Lemma~\ref{lem:affine-minimal-gauge}. An $\ai$-functor induces
a functor on graded cohomology categories: its second Taylor
identity identifies the two compositions on cohomology.
For the actual functor these are the identities of
Theorem~\ref{thm:coefficient-integration}, and strict unitality
gives $\mathsf{IN}_1(1)=1$.
Thus source homotopy equivalences become target homotopy
equivalences, and the Hom assertion can be transferred back
after it is proved for these replacements. This transfer uses
homotopy equivalences within the regular category.

Choose constant graded frames $W_P,W_Q$ for the replacements and put
\[
 T^q=\bigoplus_i\Hom_{\R}(W_P^i,W_Q^{i+q}).
\]
This is a finite-dimensional constant matrix space, placed in
internal degree $q$; it is not the full cohesive Hom module.
Filter both Hom complexes increasingly by matrix degree $q\le s$.
The filtration has finite length, independent of total degree.
Minimalization removes the degree-raising internal differential;
the subsequent one-form gauge makes the degree-preserving part
scalar. These are different roles for the two reductions.
In forms-first order its associated-graded complexes are
\[
 (B,Q)\otimes_\R T^q,\qquad
 (C_K^\bullet,\delta)\otimes_\R T^q.
\]
The leading differentials are $Q\otimes1$ and $\delta\otimes1$,
with $T^q$ in the indicated degree. For the target, $r_0=0$ and
$r_1=1$, so the two endpoint terms complete the inner differential
to the full simplicial coboundary. In particular the final
endpoint sign is
$-(-1)^{p+q}(-1)^q=(-1)^{p+1}$.

Every higher connection component has negative internal degree.
Every non-bare word in the twisted first Taylor map contains
such a component. Hence the only filtration-preserving term
is the ordinary unary kernel $\Gamma_1$, and
\[
 \operatorname{gr}^q(\mathsf{IN}^K_1)=I\otimes1_{T^q}.
\]
The unary matrix-sign sum in
\eqref{eq:in-unsuspended-matrix-sign} is empty.
This uses the actual kernel, with no second target-sign conversion.
Finite-filtration induction, using the successive short exact
sequences, proves the Hom quasi-isomorphism in every integer
degree.
\end{proof}

\subsection{Scalar source cohomology and detection}

\begin{prop}\label{prop:affine-scalar-ce-cohomology}
As a graded ring,
\[
 H^\bullet(B,Q)=C^\infty(\R_y)[u]/(yu),
 \qquad |u|=2,\qquad u=[\eta].
\]
\end{prop}

\begin{proof}
Directly from \eqref{eq:affine-hom-ce},
\begin{align}
 Q(a\theta_1+b\theta_2)
 &=(-ya+xb)\eta+
      (-ya_x+xb_x+b)\theta_1\theta_2,
 \label{eq:affine-ce-odd}\\
 Q(v\eta^n+w\theta_1\theta_2\eta^{n-1})
 &=(v_x-w)\xi\eta^n,\qquad n\ge1.
 \label{eq:affine-ce-even}
\end{align}
If a smooth scalar satisfies $g\xi=0$, it vanishes away from the
origin and then everywhere by density. In degree zero, $Qf=0$
therefore means $f_x=0$, hence $f=f(y)$.

The same smooth syzygy division as above shows that a closed
odd element is $c\xi\eta^n$. For $n=0$ it has primitive
$\int_0^x c(s,y)\,ds$; for $n\ge1$ it has primitive
$-c\theta_1\theta_2\eta^{n-1}$.
Thus all odd cohomology vanishes.

In degree $2n>0$ the cycles are exactly
\[
 z_v^{(n)}=v\eta^n+v_x\theta_1\theta_2\eta^{n-1}.
\]
Equation~\eqref{eq:affine-ce-odd} shows that boundaries correspond
exactly to $v\in(x,y)$. Global radial Hadamard division identifies
this ideal with the functions vanishing at the origin.
Consequently $H^{2n}(B)=A/(x,y)=\R$.
Finally
$z_v^{(n)}z_w^{(m)}=z_{vw}^{(n+m)}$.
The degree-zero action on $[\eta]$ is evaluation at $y=0$,
which proves the ring formula.
\end{proof}

\begin{lem}\label{lem:affine-scalar-injective}
The actual scalar integration map $I$ is injective on cohomology.
\end{lem}

\begin{proof}
At the origin, restrict $K$ to the simplices with identity
$G$-edges and integer triangle labels. Since $Z_0=0$, these
form the exact simplicial subset $B^2\Z=K(\Z,2)$.
The restriction of $h$ is its universal integral two-cocycle.
Thus its powers are nonzero in
$H^\bullet(K(\Z,2);\R)=\R[c]$.

The actual kernel gives $I(\eta)=-h$ and is multiplicative on
cohomology by the second Taylor identity in
Theorem~\ref{thm:coefficient-integration}, applied to the unit
coefficient object and pulled back by $\Psi$.
Every positive even source
generator is therefore detected. In degree zero $I$ is literally
$f(y)\mapsto f(y)$ on vertices. Together with
Proposition~\ref{prop:affine-scalar-ce-cohomology}, this proves
injectivity. The argument uses the entire exact simplicial subset,
not sampling at finitely many integer labels.
\end{proof}

\subsection{The upper-decalage relative de Rham complex}

Let
\[
 P_p=(\operatorname{Dec}^{\mathrm{up}}K)_p=K_{p+1},
 \qquad \pi_p=d_{p+1}:P_p\longrightarrow K_p.
\]
The faces and degeneracies of $P_p$ have indices $0,\ldots,p$;
the last degeneracy of $K$ is an extra degeneracy, not one of
these operators. The maps $\pi_p$ form a simplicial map.
In the root coordinates of Section~\ref{sec:affine-integration},
\begin{equation}\label{eq:affine-decalage-fibers}
 P_p=K_p\times G\times\R^p,\qquad
 g=a_{p+1},\quad z_i=h_{0,i,p+1}.
\end{equation}
In particular $\pi_p$ is globally a product submersion, even
at the singular base point.

Define the first-quadrant bicomplex
\begin{equation}\label{eq:affine-relative-bicomplex}
 \mathcal C^{p,q}=\Omega^q(P_p/K_p),\qquad
 \delta=\sum_{i=0}^{p+1}(-1)^i d_i^*,\qquad
 D_{\Tot}=\delta+(-1)^p d_{\mathrm{rel}}.
\end{equation}
These are actual jointly smooth relative forms. The simplicial
square for $\pi$ makes every face pullback well-defined on
vertical forms and commuting with $d_{\mathrm{rel}}$; hence
$D_{\Tot}^2=0$.

\begin{lem}[The relative frame]\label{lem:affine-relative-frame}
Let $m_*=gm_0$ be the last vertex. The forms
\[
 \alpha=d_{\mathrm{rel}}g\,g^{-1}
   =\begin{pmatrix}\alpha_1&\alpha_2\\0&0\end{pmatrix},
 \qquad dz_1,\ldots,dz_p
\]
give a global frame of the relative cotangent bundle.
For the root-changing face, in degrees $p\ge1$,
\begin{equation}\label{eq:affine-relative-face}
 d_0^*\alpha=\alpha-Z_{m_*}\,dz_1,\qquad
 d_0^*dz_j=dz_{j+1}-dz_1\quad(1\le j\le p-1).
\end{equation}
Other faces delete the corresponding $z_i$ and leave $\alpha$
unchanged. Degeneracy zero inserts $z_1=0$; other degeneracies
repeat a $z$-coordinate and leave $\alpha$ unchanged.
\end{lem}

\begin{proof}
The frame assertion follows from \eqref{eq:affine-decalage-fibers}
and the global right Maurer--Cartan frame on $G$.
When the root vertex is deleted, the actual nerve formulas give
\[
 g'=g\exp(-z_1Z_{m_0})a_1^{-1},\qquad
 z'_j=h_{0,1,j+1}+z_{j+1}-z_1.
\]
The additive tetrahedron identity is part of the crossed-module
nerve, including at $Z_{m_0}=0$; it is not inferred by cancelling
an exponential there.
The quantities $a_1,m_0,h_{0,1,j+1}$ are base coordinates for
$\pi_p$, so their relative differentials vanish.
Since $gZ_{m_0}g^{-1}=Z_{m_*}$, differentiation gives
\eqref{eq:affine-relative-face}. All formulas also hold when
$m_*=0$. The remaining assertions follow by deleting or repeating
the indicated vertex in the rooted triangle labels.
\end{proof}

Put $A^p=C^\infty(P_p)$ and $A=C^\infty(M)$, acting through
the last vertex. Let $V^p$ be free over $A$ on
$e_1,e_2,v_1,\ldots,v_p$, with cofaces and codegeneracies given
by Lemma~\ref{lem:affine-relative-frame}. All these operators
preserve the last vertex and are $A$-linear. The frame gives
an isomorphism of horizontal cosimplicial modules
\begin{equation}\label{eq:affine-relative-coefficients}
 \mathcal C^{p,q}=A^p\otimes_A\Lambda_A^q V^p.
\end{equation}
The second factor is finite free. Arbitrary smooth dependence on
all fiber coordinates remains in $A^p$; no algebraic tensor
description of a smooth function space on a product is assumed.

\subsection{An explicit horizontal coefficient contraction}

\begin{lem}\label{lem:affine-coefficient-contraction}
For any cosimplicial $A$-module $W^\bullet$, last-vertex pullback
and constant-simplex evaluation induce inverse horizontal
chain-homotopy equivalences
\[
 W^\bullet
 \mathrel{\mathop{\rightleftarrows}^{j}_{r}}
 A^\bullet\otimes_A W^\bullet.
\]
The assertion holds after normalization and with arbitrary
external jointly smooth parameters. In particular, $W^\bullet$
may be nonflat over $A$.
\end{lem}

\begin{proof}
Let $i^p:A\to A^p$ be last-vertex pullback, and let
$\epsilon^p:A^p\to A$ evaluate on the completely constant
$(p+1)$-simplex. They are cosimplicial $A$-algebra maps.
Set $j=i\otimes1$, $r=\epsilon\otimes1$; then $rj=1$.

There is a geometric contraction behind the formulas. For
$\lambda:[p]\to[1]$, define an order map
$a_\lambda:[p+1]\to[p+1]$ by retaining $s$ if
$s\le p$ and $\lambda(s)=0$, and sending it to $p+1$
otherwise. It fixes the last vertex.
If $\beta:[t]\to[p]$ is extended to $\beta^+$ by adding a
last vertex, then
\[
 a_\lambda\beta^+=\beta^+a_{\lambda\beta}.
\]
Both sides send a vertex either to $\beta(s)$ or to $p+1$.
Thus $K(a_\lambda)$ is a simplicial homotopy, over the last
vertex, from the identity of $P$ to the constant-last-vertex
map. It fixes constant simplices. Pullback gives a cosimplicial
right homotopy of $A^\bullet$.
Tensor its degree-$p$ operators with $1_{W^p}$.
The displayed order-map identity and functoriality of $W$
give their cosimplicial naturality and the endpoints $1,jr$.
Preservation of the last vertex makes all these operators
$A$-linear. This is tensoring an actual cosimplicial homotopy,
not tensoring a chain homotopy degreewise.

We give the resulting cochain homotopy explicitly. For
$p\ge1$ and $0\le j\le p-1$, let
\[
 f_j^p:[p+1]\longrightarrow[p],\qquad
 f_j^p(s)=
 \begin{cases}s,&s\le j,\\p,&s>j,\end{cases}
\]
and put $h_j^p=K(f_j^p):P_{p-1}\to P_p$.
Let $\sigma_W^j:W^p\to W^{p-1}$ be the codegeneracy
identifying $j$ and $j+1$. Define
\begin{equation}\label{eq:affine-coefficient-homotopy}
 L_p=\sum_{j=0}^{p-1}(-1)^{j+1}
           (h_j^p)^*\otimes\sigma_W^j,\qquad L_0=0.
\end{equation}
Then
\begin{equation}\label{eq:affine-coefficient-homotopy-identity}
 \delta L+L\delta=1-jr.
\end{equation}
To verify the identity in every degree,
write $\mathbf A_{ij}$ for the unsigned tensor operator in
$\delta_{p-1}L_p$, and $\mathbf B_{ij}$ for that in
$L_{p+1}\delta_p$. Their signs are $(-1)^{i+j+1}$.
As pairs of order maps, the off-diagonal terms satisfy
\[
 \mathbf B_{ij}=
 \begin{cases}
   \mathbf A_{i,j-1},&i<j,\\
   \mathbf A_{i-1,j},&i>j+1.
 \end{cases}
\]
The matched signs are opposite, and these cases exhaust all
$\mathbf A$ terms. In the remaining cases $i=j,j+1$, the
$W$-map is the identity. On the coefficient factor their
signed sum is $c_j-c_{j-1}$, where $c_j$ retains vertices
$0,\ldots,j$ and sends all others to $p+1$.
Summing over $j=0,\ldots,p$ gives $c_p-c_{-1}$, namely
the identity minus constant-last-vertex evaluation.
This proves \eqref{eq:affine-coefficient-homotopy-identity},
including $p=0$.

Use the normalization
$N^pT=\bigcap_{i=0}^{p-1}\ker(s^i:T^p\to T^{p-1})$.
A compatible finite projection is
\[
 \Pi_N^p=(1-d^p s^{p-1})\cdots(1-d^1s^0).
\]
Its complement is the sum of the positive-coface images;
the cosimplicial identities show that both summands are
subcomplexes. Thus $\Pi_N$ and the normalized inclusion
$\iota_N$ are cochain maps. Since $j,r$ are cosimplicial,
$L_N=\Pi_N L\iota_N$ satisfies the same homotopy identity.
For $T=\mathcal C^{\bullet,q}$ the projection also commutes
with $d_{\mathrm{rel}}$: it uses only genuine $P$-operators.
All formulas are finite smooth pullbacks and remain valid
with jointly smooth external parameters.
\end{proof}

\begin{rem}\label{rem:affine-horizontal-only}
The retraction and homotopy in the lemma are horizontal only.
For example, at $p=0$ the coefficient $F(g,m)=\log a$ satisfies
$r(F)=0$, whereas $r(d_{\mathrm{rel}}F)=e_1$.
Thus this horizontal contraction does not commute with the
relative differential. The total comparison instead uses the
total chain map constructed below and a first-quadrant filtration.
\end{rem}

\begin{lem}[Exact exterior normalization]
\label{lem:affine-exterior-normalization}
Over $A$ itself,
\[
 N^p\Lambda_A^q V^\bullet
 =\Lambda_A^{q-p}(Ae_1\oplus Ae_2)
       \wedge v_1\wedge\cdots\wedge v_p.
\]
It vanishes unless $0\le q-p\le2$.
\end{lem}

\begin{proof}
The codegeneracy kernels on $V^p$ are the independent lines
generated by
$w_1=v_1,w_2=v_2-v_1,\ldots,w_p=v_p-v_{p-1}$.
Together with $e_1,e_2$ these form a free $A$-basis.
On the exterior algebra, each kernel is the ideal generated
by the corresponding $w_i$. A basis monomial belongs to all
these ideals precisely when it contains all $p$ symbols.
Their wedge is $v_1\wedge\cdots\wedge v_p$.
For $p=0$ the intersection is empty and the formula is
$\Lambda^q(Ae_1\oplus Ae_2)$.
The argument uses a free basis over $A$, not a fiberwise rank
calculation at regular points.
\end{proof}

\subsection{The actual scalar comparison}

\begin{thm}[Scalar affine comparison]\label{thm:affine-scalar-comparison}
The actual map
$I:(B,Q)\to(C_K^\bullet,\delta)$ is a quasi-isomorphism.
Consequently
\[
 H^\bullet(C_K^\bullet)
   =C^\infty(\R_y)[u]/(yu),\qquad u=[-h],\quad |u|=2,
\]
with this ring identification induced by $I$.
\end{thm}

\begin{proof}
There are two different filtrations to use.
First, $\pi_p^*$ defines a cochain map from unnormalized smooth
cochains of $K$ to $\Tot\mathcal C$.
The fiber in \eqref{eq:affine-decalage-fibers} has global
coordinates $(\log a,b,z)$ and contracts radially to $(0,0,0)$.
The relative Poincare operator integrates over a compact
interval and is smooth in all base parameters. Each column
therefore has cohomology $C^\infty(K_p)$ in relative degree
zero and zero otherwise. Filtering by the horizontal degree
$p$ gives a quasi-isomorphism
\[
 J:C^\bullet_{\mathrm{sm}}(K)\longrightarrow\Tot\mathcal C.
\]
No compatibility of the chosen Poincare operators with faces
is needed, since $J$ itself is a cochain map.
Ordinary cosimplicial normalization is a natural
chain-homotopy equivalence, constructed from the genuine
face and degeneracy maps. These commute with
$d_{\mathrm{rel}}$, as does $J$.
Thus the same conclusion holds for $C_K^\bullet$ and
$\Tot N\mathcal C$.
This is full cosimplicial normalization, not normalization of
the inner-chain complex used for higher prisms.

For the second comparison give $B$ the bigrading
\[
 B^{p,q}=A\eta^p\Lambda^{q-p}(\theta_1,\theta_2),
 \qquad q-p\in\{0,1,2\}.
\]
The total degree is $p+q$.
Split $Q=Q_h+Q_v$, where $Q_h=\eta\,\iota_Z$ for
$Z=(-y,x)$, and
\[
 Q_v f=f_x\xi,\quad Q_v\theta_1=0,\quad
 Q_v\theta_2=\theta_1\theta_2,\quad Q_v\eta=0.
\]
These components square to zero and anticommute.
For the mixed identity on $\theta_2$ the two terms are
$-\xi\eta$ and $+\xi\eta$; on functions use $Q_h\xi=0$.
The other generators give zero directly.
Thus this is a genuine first-quadrant double complex.
For a theta monomial $\theta_I$ of exterior degree $r$, put
\begin{equation}\label{eq:affine-iota}
 \iota(f\eta^p\theta_I)
   =c_{p,r}f(m_*)\alpha_I\wedge dz_1\wedge\cdots\wedge dz_p,
 \qquad
 c_{p,r}=(-1)^{pr+p(p+1)/2}.
\end{equation}
Lemma~\ref{lem:affine-exterior-normalization} identifies its
finite free part exactly with $N^p\Lambda^q V$.

We verify the total chain-map identity directly. In the horizontal
differential the terms without the root-face change of $\alpha$ cancel:
\[
 \delta(dz_1\wedge\cdots\wedge dz_p)=0.
\]
Expanding
$\prod_{j=2}^{p+1}(dz_j-dz_1)$ gives the signed sum of the
wedges omitting one $dz_i$, which is exactly this cancellation.
The remaining root-face term is
\[
 \delta(\alpha_I\wedge dz_1\wedge\cdots\wedge dz_p)
 =(-1)^r(\iota_Z\alpha_I)\wedge
                    dz_1\wedge\cdots\wedge dz_{p+1}.
\]
Since $c_{p+1,r-1}=(-1)^r c_{p,r}$, it gives
$\delta\iota=\iota Q_h$.
For the vertical differential the right Maurer--Cartan
equation gives
\[
 d_{\mathrm{rel}}f(m_*)=f_x(m_*)(x_*\alpha_1+y_*\alpha_2),
 \quad d_{\mathrm{rel}}\alpha_1=0,\quad
 d_{\mathrm{rel}}\alpha_2=\alpha_1\wedge\alpha_2.
\]
All $dz_i$ are closed.
The identity $c_{p,r+1}=(-1)^pc_{p,r}$ therefore gives
$(-1)^p d_{\mathrm{rel}}\iota=\iota Q_v$.
Thus \eqref{eq:affine-iota} is an actual total chain map.

For every fixed $q$, its horizontal part is the normalized
coaugmentation of Lemma~\ref{lem:affine-coefficient-contraction},
after the exact exterior identification. It is a rowwise
quasi-isomorphism. Filtering now by the relative degree $q$
shows that $\iota:B\to\Tot N\mathcal C$ is a
quasi-isomorphism. Both comparisons are first quadrant with
only finitely many summands in each total degree. We have
therefore constructed a chain zigzag
\begin{equation}\label{eq:affine-scalar-zigzag}
 C_K^\bullet\xrightarrow[\sim]{J}\Tot N\mathcal C
                         \xleftarrow[\sim]{\iota}(B,Q).
\end{equation}
The zigzag computes the graded vector spaces. To identify
the product for the actual map $I$, we use its cohomological
multiplicativity and the injectivity proved above.
Positive even target cohomology has dimension one and odd
target cohomology vanishes. Lemma~\ref{lem:affine-scalar-injective}
therefore makes $H(I)$ an isomorphism in positive degrees.
For degree zero argue directly: a smooth $G$-invariant
function on $M$ is independent of $x$ for $y\ne0$, since the
translation parameter moves $x$ arbitrarily. Smoothness
extends this independence to $y=0$.
Thus target $H^0=C^\infty(\R_y)$, and $I$ is literally the
identity on it. Negative scalar cochain degrees vanish.
This proves quasi-isomorphism of the actual $I$.
Its $\ai$ multiplicativity on cohomology, together
with $I(\eta)=-h$, gives the stated ring identification.
\end{proof}

\begin{cor}[External scalar parameters]
\label{cor:affine-scalar-parameters}
The same actual scalar comparison is a quasi-isomorphism for
jointly smooth coefficients on $U\times\R^2$ and cochains on
$U\times K_p$, with the differential only in the affine
variables, for every external smooth parameter manifold $U$.
\end{cor}

\begin{proof}
Put $A_U=C^\infty(U\times M)$ and
$A_U^p=C^\infty(U\times P_p)$, with the latter an $A_U$-algebra
through $(u,m_*)$. The actual relative complex is
\[
 \mathcal C_U^{p,q}
   =\Omega^q((U\times P_p)/(U\times K_p))
   =A_U^p\otimes_{A_U}\Lambda^q V_U^p,
 \qquad V_U^p=A_U^{p+2},
\]
in the same frame $\alpha_1,\alpha_2,dz_1,\ldots,dz_p$.
The tensor description uses this finite free frame and retains
all joint smooth dependence in $A_U^p$.
All parameterized operators are obtained by taking the product
with $1_U$, so they preserve the external parameter $u$.
The genuine upper-decalage simplicial operators and the horizontal
order maps used in $L$ preserve the terminal vertex $m_*$.
Their coefficient pullback and evaluation maps are therefore
$A_U$-linear through $(u,m_*)$, and the horizontal identity
$\delta L+L\delta=1-jr$ and the exterior normalization hold
over $A_U$.

The column radial homotopy is instead a homotopy over
$U\times K_p$. It generally moves $m_*=gm_0$ and need not
be $A_U$-linear for the terminal action. Its relative Poincare
operator is $C^\infty(U\times K_p)$-linear, which is the
property required for the column cohomology calculation.
Both constructions are $C^\infty(U)$-linear and commute
with smooth pullback of the external parameter.

The column Poincare operator still integrates only the radial
fiber variables $(\log a,b,z)$ over $[0,1]$.
It is smooth in $(u,K_p)$, and commutes with every smooth
change of external parameters. The displayed formula for
$\iota$ also remains a total chain map, since no differential
acts in the $U$ directions. The same two first-quadrant
filtrations therefore give actual quasi-isomorphisms
\[
 C^\bullet_{\mathrm{norm,sm}}(U\times K)
   \xrightarrow[\sim]{J_U}\Tot N\mathcal C_U
   \xleftarrow[\sim]{\iota_U}(B_U,Q),
 \qquad B_U=A_U[\theta_1,\theta_2,\eta].
\]
All three complexes and both arrows commute with pullback
along $U'\to U$ through the explicit parameterized
operators above.

The primitives and radial Hadamard formula in
Proposition~\ref{prop:affine-scalar-ce-cohomology} integrate
only the $x,y$ variables and remain jointly smooth in $u$.
The source degree-zero cohomology is
$C^\infty(U\times\R_y)$, the positive even cohomology is
$C^\infty(U)$, and the odd cohomology vanishes.
Let $S_U=H(\iota_U)^{-1}H(J_U)$ be the resulting natural
$C^\infty(U)$-linear cohomology isomorphism.
The actual $I_U$ is also parameter-natural: its simplex
integrals keep $u$ fixed, and hence commute with such pullbacks.
For $\pi:U\to *$ and $n\ge1$ this gives
\[
 S_UH(I_U)[\eta^n]
   =\pi^*\bigl(S_*H(I_*)[\eta^n]\bigr)
   =\lambda_n[\eta^n],\qquad \lambda_n\in\R^\times.
\]
Nonvanishing follows from the unparameterized theorem.
Multiplication by this fixed real constant is invertible.
The fixed nonzero scalar gives an inverse on the rank-one
$C^\infty(U)$-module in each positive even degree.
The degree-zero invariance argument works for each $u$
and gives literally $f(u,y)\mapsto f(u,y)$.
Joint smooth dependence is carried by the displayed
parameterized complexes on the whole affine plane.
\end{proof}

\begin{proof}[Proof of Theorem~\ref{thm:affine-hom}]
Apply Theorem~\ref{thm:affine-scalar-comparison} to
Lemma~\ref{lem:affine-hom-reduction}.
For the full MC target, the actual functor $\Psi^*$ is the
dg quasi-equivalence inverse to $\Phi^*$ from
Corollary~\ref{cor:affine-representation-invariance}.
Since $\mathsf{IN}^K=\Psi^*\mathsf{IN}_E$, its Hom
quasi-isomorphisms give the corresponding assertion for
$\mathsf{IN}_E$ through the constructed raw comparison.
\end{proof}

\section{Regular reconstruction and finite MC lifting}
\label{sec:affine-reconstruction}
\label{subsec:affine-reconstruction}

We use the whole-plane notation $E$, $B$, $K$, $\mathcal R_K$
and $\mathsf{IN}^K$ of Appendix~\ref{sec:affine-hom}.
We complement the regular Hom comparison with a finite
reconstruction proving essential surjectivity in the regular
raw category. The finite MC-lifting lemma below is also used
in the vertical-family argument of Section~\ref{sec:affine-families}.

\begin{lem}[The raw cochain-module dictionary]
\label{lem:affine-raw-module-dictionary}
In global graded frames, raw representations of $K$ are the finite
free graded right $C_K^\bullet$-modules with superconnections
encoded by
\[
 r=e+b,\qquad d_{\mathrm{full}}b+b\star b=0,
 \qquad e_1=1,\quad e_p=0\ (p\ne1).
\]
Their Hom differential is exactly the raw representation differential.
The degree-zero algebra $C_K^0=C^\infty(M)$ acts at the
two endpoints.
\end{lem}

\begin{proof}
Every bounded graded vertex bundle is trivial on $M=\R^2$;
choose its constant graded vector space $W$. In forms-first order,
the endomorphism algebra of the corresponding right module is
$C_K^\bullet\otimes\End_\R(W)$, with product
\[
 (\omega\otimes T)(\nu\otimes S)
    =(-1)^{|T||\nu|}(\omega\smile\nu)\otimes TS.
\]
Indeed its action on $C_K^\bullet\otimes W$ is
$\Phi_{\omega,T}(c\otimes v)=(-1)^{|T||c|}\omega c\otimes Tv$.
Here the right module structure is the braided one,
$(c\otimes v)a=(-1)^{|v||a|}(ca)\otimes v$.
This makes the displayed endomorphisms right-linear with
the signless convention for homogeneous morphisms.
The formula uses associativity, not commutativity of $C_K$.

The unit representation satisfies
$\delta_{\mathrm{in}}e+e\star e=0$, and
\[
 d_{\mathrm{full}}\alpha
 =\delta_{\mathrm{in}}\alpha+
       e\star\alpha-(-1)^{|\alpha|}\alpha\star e.
\]
Substituting $r=e+b$ in the raw MC equation therefore gives
the displayed ordinary DGA MC equation.
For a morphism from $b_P$ to $b_Q$, the differential becomes
\[
 d_{\mathrm{full}}\alpha+
 b_Q\star\alpha-(-1)^{|\alpha|}\alpha\star b_P,
\]
the superconnection Hom differential.
Degree-zero functions multiply a cochain at its first vertex
on the left and its last vertex on the right. Their
noncentrality is exactly what the cochain product records.
Normalization gives the strict units.
The degree-two MC identities on the $h=0$ inverse-arrow
triangles exhibit the reverse edge as a chain-homotopy
inverse, so the edge quasi-isomorphism condition also holds.
Changing a vertex frame by $g$ gives
$b^g=g^{-1}bg+g^{-1}d_{\mathrm{full}}g$,
so the dictionary is independent of the chosen frames.
\end{proof}

\begin{lem}[Minimal target and trivial edge transport]
\label{lem:affine-target-minimalization}
Every regular object of $\mathcal R_K$ is homotopy equivalent
inside that dg category to one with $r_0=0$ and $r_1=1$.
\end{lem}

\begin{proof}
Apply Lemma~\ref{lem:affine-finite-minimalization} to the
generally noncommutative DGA $C_K$ using
Lemma~\ref{lem:affine-raw-module-dictionary}.
In the forms-first frames, extend the contraction by left
$\star$-multiplication by its zero-cochains.
For a degree-$j$ matrix $T$ and a scalar $p$-cochain $c$,
the actual extension is
\[
 (\widehat T(c\otimes v))(\sigma)
    =(-1)^{jp}c(\sigma)T(v_0\sigma)v.
\]
Thus $h$ carries $(-1)^p$, not unsigned pointwise extension.
Associativity gives right $C_K$-linearity.
The signs from $d_0$ and $h$ cancel in $d_0h+hd_0$;
$i,p$ have degree zero, so all contraction side conditions
remain valid.
The proof of right-linearity there uses only the Leibniz rule
and the contraction side conditions. Its perturbation has
positive cochain degree, so the internal-degree series is finite.
All operations preserve actual smooth normalized cochains.
The resulting vertex cohomology bundle $H$ is again globally
trivial, and its internal differential is zero.

The transferred edge map is its genuine cohomology transport.
Write $p_j=p(v_j)$ and $i_j=i(v_j)$ for the contraction matrices
at the endpoints. The first positive-degree term is
\[
 b'_1=p_0\bigl[(i_1-i_0)+(r_1-1)i_1\bigr]
       =p_0r_1i_1-1.
\]
The first bracketed term is $d_{\mathrm{full}}i$; it cannot
be omitted. Thus $r'_1=p_0r_1i_1$, an isomorphism on the
zero-differential cohomology fibers, equal to $1$ on a unit.
Moreover $D_H-d_{\mathrm{full}}$ is right-linear, so it is
left multiplication by a normalized cochain $b'$.
Its degree-zero part vanishes. These facts give an actual
minimal object of the raw representation category.

Now assume $r_0=0$ and denote the edge map by $R(g,m)$,
from the fiber at $gm$ to the fiber at $m$.
The degree-two MC equation says
\begin{equation}\label{eq:affine-minimal-edge}
 R(g_{12}g_{01}\exp(hZ_m),m)
       =R(g_{01},m)\,R(g_{12},g_{01}m).
\end{equation}
At $h=0$ this is a strict contravariant action-groupoid
representation. Units and inverse arrows make its edge maps
invertible. Setting both displayed $g_{ij}$ to the identity
gives $R(\exp(hZ_m),m)=1$.
Here $\rho(Z_m)|_m=0$, so $\exp(hZ_m)$ fixes $m$.
The triangle and the last-to-first transport are therefore
typed over the stated vertices, including the origin.

Let $E_1,E_2$ be the matrix generators whose fundamental fields
are $x\partial_x,y\partial_x$, and set
\[
 L_i(m)=\left.\frac{d}{dt}\right|_{t=0}R(\exp(tE_i),m).
\]
Differentiating the boundary identity gives
$-yL_1+xL_2=0$. Global Hadamard division yields
$L_1=xL$, $L_2=yL$ as before.
The graded vertex gauge
\[
 g_x=-Lg,\qquad g(0,y)=1
\]
is globally smooth and invertible. It changes the
infinitesimal matrices to
$g^{-1}L_i g+g^{-1}\rho_i(g)=0$.

This kills the whole edge representation, not merely its
first derivative. Along a smooth path $a(t)$ in $G$ starting
at the identity, write
$a(t+\epsilon)=c_t(\epsilon)a(t)$.
Then $c_t(0)=1$ and
$c_t'(0)=\dot a(t)a(t)^{-1}$.
The contravariant cocycle equation gives
\[
 \frac{d}{dt}R^g(a(t),m)
 =R^g(a(t),m)\,
    L^g_{\dot a(t)a(t)^{-1}}(a(t)m)=0.
\]
The group $G=\Aff_+(1)$ is connected, so $R^g(a,m)=1$
for every arrow. Here the ODE gives an ordinary change of
frame, rather than algebroid transport from the origin.
The ordinary gauge is a strict closed
isomorphism from the edge-trivialized representation to
the preceding minimal one.
\end{proof}

\subsection{Finite filtered Maurer--Cartan lifting}

We give the finite obstruction argument needed for reconstruction.
The finite filtration makes every sum and gauge below terminate.

\begin{lem}[Finite MC lifting]\label{lem:finite-mc-lifting}
Let $L,M$ be nonunital DGAs over $\R$ with finite descending
filtrations
$F^1L=L$, $F^{\nu+1}L=0$, and similarly for $M$.
Assume the differential preserves filtration and
$F^iF^j\subset F^{i+j}$.
Let $F:L\rightsquigarrow M$ be an uncurved $\ai$-map with the
arity-additive bounds
\begin{equation}\label{eq:mc-arity-filtration}
 F_n(F^{i_1}L,\ldots,F^{i_n}L)
       \subset F^{i_1+\cdots+i_n}M.
\end{equation}
Here gauges are elements $g\in1+M^0$ in the unitization, with
\[
 b^g=g^{-1}bg+g^{-1}dg.
\]
Their inverses are finite by the filtration.
If $\operatorname{gr}^\ell F_1$ is a quasi-isomorphism for
every $\ell$, every target MC element is gauge-equivalent
to $F_*(a)$ for a source MC element $a$.
\end{lem}

\begin{proof}
Use the suspended bar convention of
Section~\ref{sec:coefficient-integration}.
If $\overline F_n:(sL)^{\otimes n}\to sM$ are its bar
components, define
\[
 F_*(a)=s^{-1}\sum_{n=1}^{\nu}
                 \overline F_n((sa)^{\otimes n}).
\]
Higher terms vanish by \eqref{eq:mc-arity-filtration}.
This is the primitive part of the formal group-like bar
expression; only its finite filtered components are used.
The suspension and tensor signs follow the convention fixed in
Section~\ref{sec:coefficient-integration}.

Put $\Omega(a)=da+a^2$.
Since $|sa|=0$, the bar identities
$b_1(sa)=-s(da)$ and
$b_2(sa,sa)=-s(a^2)$ give the curvature identity
\[
 s\Omega(F_*(a))
  =\sum_{r,t\ge0}
     \overline F_{r+t+1}
       ((sa)^{\otimes r},s\Omega(a),(sa)^{\otimes t}).
\]
All terms in this primitive expression are finite by the
same filtration bound.

Suppose, at stage $\ell$, that
$\Omega(a)\in F^\ell L^2$ and that the current target MC
element $b$, after the accumulated target gauge, satisfies
$c=b-F_*(a)\in F^\ell M^1$.
This holds initially with $\ell=1$ and $a=0$.
Work modulo $F^{\ell+1}$.
The Bianchi identity makes $\Omega(a)$ closed in
$\operatorname{gr}^\ell L$. In the curvature identity,
every additional insertion of $a\in F^1$ raises filtration,
so
\[
 0=\Omega(b)\equiv F_1(\Omega(a))+dc.
\]
Indeed the omitted products $F_*(a)c,cF_*(a)$ have
filtration at least $\ell+1$, and $c^2$ has filtration
$2\ell\ge\ell+1$.
Injectivity on $H^2$ gives $u\in F^\ell L^1$ with
$du\equiv-\Omega(a)$. Replacing $a$ by $a+u$ kills its
curvature in this layer.
Recompute $c=b-F_*(a)$ after this correction; it changes
by $-F_1(u)$ modulo $F^{\ell+1}$.

The new discrepancy $c$ is closed in the target layer.
Surjectivity on $H^1$ gives a source representative
$t\in F^\ell L^1$, closed modulo $F^{\ell+1}$, whose
image has the same class. Replacing $a$ by $a+t$
preserves the improved curvature condition.
Again recompute the discrepancy, now subtracting $F_1(t)$
in this layer.
The remaining discrepancy is $dz$ modulo $F^{\ell+1}$,
for $z\in F^\ell M^0$.

In the unitization take $g=1-z$; its inverse is a finite geometric series.
Since products with $b\in F^1$ have higher filtration,
\[
 b^g=g^{-1}bg+g^{-1}dg
       \equiv b-dz\equiv F_*(a)\pmod{F^{\ell+1}}.
\]
This completes the induction. After finitely many stages,
$\Omega(a)=0$ and $b^q=F_*(a)$ for the accumulated gauge
$q=g_1g_2\cdots g_\nu$ in induction order:
$b^{g_1g_2}=(b^{g_1})^{g_2}$.
Its orientation is fixed by
\begin{equation}\label{eq:mc-gauge-orientation}
 (d+b)q=q(d+b^q).
\end{equation}
Thus $q$ is a closed invertible map from $F_*(a)$ to $b$,
giving an actual morphism between coefficient objects by
finite induction.
\end{proof}

\subsection{The global affine Riemann--Hilbert equivalence}

\begin{thm}[Global regular affine RH]\label{thm:affine-rh}
The actual integration functors
\[
 \mathsf{IN}^K:\Coh^{\mathrm{reg}}(B)
       \longrightarrow\mathcal R_K^{\mathrm{reg}},
 \qquad
 \mathsf{IN}_E:\Coh^{\mathrm{reg}}(B)
       \longrightarrow
       \Rep_{\mathrm{sm}}^{\mathrm{perf,reg}}(\Int(E))
\]
are strictly unital $\ai$ quasi-equivalences.
For each target object there is a bounded finite-projective
regular cohesive preimage and an actual closed degree-zero
homotopy equivalence from its integration to that object.
\end{thm}

\begin{proof}
Quasi-full faithfulness is Theorem~\ref{thm:affine-hom}.
It remains to construct a preimage of an arbitrary regular
raw representation $R$ of $K$.
Lemma~\ref{lem:affine-target-minimalization} gives a minimal,
edge-trivialized representation on a fixed finite graded
space $H$. Write it as $e+b$.

Let
\[
 N=\bigoplus_{q<0}\End_\R^q(H),\qquad
 F^\ell N=\bigoplus_{q\le-\ell}\End_\R^q(H).
\]
This is a finite multiplicative nilpotent filtration.
Because $r_0=0$, $r_1=1$, every remaining structure
cochain lies in $C_K^\bullet\otimes N$, in total degree one.
On the source use $B\otimes N$.
A total-degree-one element there has CE degree at least two.

Restrict the actual $\mathsf{IN}^K$ to the endomorphism
algebra of the trivial rank-one coefficient and take its
graded matrix extension
\eqref{eq:in-matrix-sign}--\eqref{eq:in-unsuspended-matrix-sign}.
Its $\ai$ property is the forms-first matrix extension
already fixed in Theorem~\ref{thm:coefficient-integration},
using \cite[Proposition~A.17]{AriasAbadSchaetz2013}.
It gives the $\ai$ map
\[
 F^N:B\otimes N\rightsquigarrow C_K^\bullet\otimes N.
\]
The matrix products make its Taylor filtration arity-additive.
In layer $\ell$, its first component is
\[
 I\otimes 1_{\End^{-\ell}(H)},
\]
with that matrix space placed in internal degree $-\ell$.
It is a quasi-isomorphism by
Theorem~\ref{thm:affine-scalar-comparison}.
Only tensoring with a finite-dimensional real vector space
occurs here.

Lemma~\ref{lem:finite-mc-lifting} produces an MC element
$a\in(B\otimes N)^1$ and a finite gauge $q$ with
$b^q=F^N_*(a)$.
Let $P=(H\otimes B,Q+a)$.
It is bounded finite projective and regular, with zero
internal differential and zero one-form part.

To identify the actual object, use
$\overline\Gamma_n=s\Gamma_n(s^{-1})^{\otimes n}$.
The tensor sign is
\[
 (s^{-1})^{\otimes n}
     (sa_1\otimes\cdots\otimes sa_n)
 =(-1)^{\sum_i(n-i)(|a_i|-1)}
     a_1\otimes\cdots\otimes a_n.
\]
Every MC input has total degree $|a_i|=1$, including its
matrix degree, so this exponent is zero. Therefore
$F^N_*(a)=\sum_n\Psi^*\Gamma_n^H(a,\ldots,a)$ literally,
with the unsuspended sign
\eqref{eq:in-unsuspended-matrix-sign}.
No extra $(-1)^{n(n-1)/2}$ occurs.
Since $\Psi^*$ is a strict map of normalized cochain DGAs,
\eqref{eq:in-object-formula} gives
\[
 \mathsf{IN}^K(P)=e+F^N_*(a).
\]
No second target conversion $T$ is applied.
In particular the scalar component still sends $\eta$ to
$-h$, as fixed by \eqref{eq:affine-surface-coordinate}.
Every object word has strictly decreasing internal degree,
so these particular twisting sums are finite.

The gauge $q$ is an actual normalized raw representation
isomorphism toward the edge-trivialized target, by
\eqref{eq:mc-gauge-orientation}. Every degree-zero
correction in $F^\ell(C_K\otimes N)$ has positive
cochain degree, so $q_0=1$.
Let $g$ denote the ordinary vertex gauge from that object
to the minimal representation, and let $i'$ be the
perturbed inclusion from the minimal representation to
$R$. Then
\[
 \epsilon_R=i'\,g\,q:\mathsf{IN}^K(P)\longrightarrow R
\]
is closed, degree zero, and a homotopy equivalence.
The two gauge factors are strictly invertible; $i'$ is the
homotopy equivalence supplied by finite transfer.
This proves essential surjectivity on $H^0$.

Finally $\Psi^*$ is the actual raw dg quasi-equivalence
of Corollary~\ref{cor:affine-representation-invariance},
and $\mathsf{IN}^K=\Psi^*\mathsf{IN}_E$.
More explicitly, that corollary supplies closed invertible
gauges $U_Z:Z\to\Phi^*\Psi^*Z$ on the full MC target.
For a reconstructed map
$v:\Psi^*\mathsf{IN}_E(P)\to\Psi^*R$, the actual full-MC map is
\[
 U_R^{-1}\,\Phi^*(v)\,U_{\mathsf{IN}_E(P)}
       :\mathsf{IN}_E(P)\longrightarrow R.
\]
It is a closed degree-zero homotopy equivalence.
The Hom quasi-isomorphisms likewise transfer through
$\Psi^*$, proving both assertions for the intrinsic functor.
\end{proof}

\begin{cor}[Other global affine resolutions]
\label{cor:affine-rh-global-resolutions}
Let $E'$ be another supplied strict global geometric resolution
of $\mathfrak f=\langle x\partial_x,y\partial_x\rangle$, with
$E,E'\in\mathscr R_{\mathfrak f}$, under the
fixed-base comparison hypotheses of
Theorem~\ref{main:geometric-invariance}.
Then $\mathsf{IN}_{E'}$ is a quasi-equivalence on the
regular cohesive and raw perfect representation categories.
\end{cor}

\begin{proof}
By Theorem~\ref{thm:coherent-resolutions}, choose a comparison
$F:E\rightsquigarrow E'$ commuting with the resolution augmentations.
Both complexes resolve $\mathfrak f$, so $F\in\Wone$.
Its cohesive pullback and actual full-MC representation
pullback are dg quasi-equivalences on the regular classes.
The naturality of $\mathsf{IN}$ identifies the resulting
square with the affine integration square. Quasi-full
faithfulness and essential surjectivity transfer through
these two equivalences.
\end{proof}

\begin{rem}[Reconstruction and choices]\label{rem:affine-rh-scope}
The preimages and maps $\epsilon_R$ were chosen objectwise.
This appendix gives a global regular theorem for the specified
affine foliation and its stated global resolutions, using the
whole-plane section $\Psi$.
Theorem~\ref{thm:affine-open-rh} separately proves RH for
nonregular coefficients on every ordinary open by a different
comparison mechanism and separately chosen range-preserving data.
The regular tools here remain inputs to the vertical-family
theorem. Effective descent through a refined holonomy quotient
is a further question, as discussed in Section~\ref{sec:programme}.
\end{rem}

\section*{Acknowledgments}

This paper grew out of part of the author's Ph.D. thesis work.
The author is grateful to Jonathan Block for his advice and support,
to David White for discussions on the homotopy theory of algebras
over operads, and to Sylvain Lavau and Jim Stasheff for discussions
on $L_\infty$-algebroids and singular foliations.

\section*{AI assistance}
This manuscript was developed with substantial assistance from
GitHub Copilot CLI and Azure OpenAI GPT-5.6 Sol, including
mathematical proof exploration, drafting and revision of the
exposition, and automated consistency checks. Preparation of
the present revision also used Azure OpenAI GPT-6 Astra.

\printbibliography
\end{document}